\documentclass[12pt]{amsart}
\usepackage{amsbsy}
\usepackage{graphicx,epsfig,subfigure,psfrag}
\begin{document}

\newtheorem{theorem}{Theorem}
\newtheorem{proposition}{Proposition}
\newtheorem{lemma}{Lemma}
\newtheorem{corollary}{Corollary}
\newtheorem{definition}{Definition}
\newtheorem{remark}{Remark}
\newcommand{\beq}{\begin{equation}}
\newcommand{\eeq}{\end{equation}}
\numberwithin{equation}{section} \numberwithin{theorem}{section}
\numberwithin{proposition}{section} \numberwithin{lemma}{section}
\numberwithin{corollary}{section}
\numberwithin{definition}{section} \numberwithin{remark}{section}
\newcommand{\ren}{\mathbb{R}^N}
\newcommand{\re}{\mathbb{R}}
\newcommand{\n}{\nabla}
\newcommand{\iy}{\infty}
\newcommand{\pa}{\partial}
\newcommand{\fp}{\noindent}
\newcommand{\ms}{\medskip\vskip-.1cm}
\newcommand{\mpb}{\medskip}
\newcommand{\BB}{{\bf B}}
\newcommand{\Am}{{\bf A}_{2m}}
\renewcommand{\a}{\alpha}
\renewcommand{\b}{\beta}
\newcommand{\g}{\gamma}
\newcommand{\G}{\Gamma}
\renewcommand{\d}{\delta}
\newcommand{\D}{\Delta}
\newcommand{\e}{\varepsilon}
\newcommand{\var}{\varphi}
\renewcommand{\l}{\lambda}
\renewcommand{\o}{\omega}
\renewcommand{\O}{\Omega}
\newcommand{\s}{\sigma}
\renewcommand{\t}{\tau}
\renewcommand{\th}{\theta}
\newcommand{\z}{\zeta}
\newcommand{\wx}{\widetilde x}
\newcommand{\wt}{\widetilde t}
\newcommand{\noi}{\noindent}
\newcommand{\inA}{\quad \mbox{in} \quad \ren \times \re_+}
\newcommand{\inB}{\quad \mbox{in} \quad}
\newcommand{\inC}{\quad \mbox{in} \quad \re \times \re_+}
\newcommand{\inD}{\quad \mbox{in} \quad \re}
\newcommand{\forA}{\quad \mbox{for} \quad}
\newcommand{\whereA}{,\quad \mbox{where} \quad}
\newcommand{\asA}{\quad \mbox{as} \quad}
\newcommand{\andA}{\quad \mbox{and} \quad}
\newcommand{\ef}{\eqref}
\newcommand{\ssk}{\smallskip}
\newcommand{\LongA}{\quad \Longrightarrow \quad}
\def\com#1{\fbox{\parbox{6in}{\texttt{#1}}}}

\title
{\bf  Shock waves and compactons for\\ fifth-order nonlinear
dispersion equations. II}

\author {
Victor A.~Galaktionov}

\address{Department of Mathematical Sciences, University of Bath,
 Bath BA2 7AY, UK}
\email{vag@maths.bath.ac.uk}



 \date{\today}

\begin{abstract}

The following {\em first problem} is posed:
 to justify that the {\em standing shock wave}
 \beq
  S_-(x) = -{\rm sign}\, x =\{
 -1 \,\,\, \mbox{for} \,\,\, x<0, \,\,\,
 \,\,\,\,1 \,\,\, \mbox{for} \,\,\, x > 0\},
 \notag
  \eeq
 is a correct ``entropy" solution of the Cauchy problem for  the fifth-order degenerate nonlinear
 dispersion equations, NDEs (same as for the classic Euler one $u_t+ u u_x=0$)
  \beq
  u_t=-(u  u_{x})_{xxxx} \quad \mbox{and}
  \quad u_t=-(u u_{xxxx})_x
  \quad \mbox{in} \quad \re \times \re_+.
  \notag
   \eeq
 These two quasilinear
degenerate PDEs
 are chosen as
typical representatives, so other $(2m+1)$th-order NDEs
 of non-divergent form admit such shocks waves.
 As a  related {\em second
problem},  the opposite initial shock $S_+(x)=-S_-(x)= {\rm
sign}\, x$ is shown to be a non-entropy solution creating  a
 {\em rarefaction wave}, which becomes $C^\iy$  for any $t>0$.
 Formation of shocks
leads to  nonuniqueness of any ``entropy solutions".
Similar phenomena are studied for a {\em fifth-order in time} NDE
 $
 u_{ttttt}=(uu_x)_{xxxx}
 $
 in {\em normal form}.
Other  NDEs,
 \beq
\mbox{e.g.,} \quad  u_t=-(|u| u_{x})_{xxxx} + |u| u_x
  \quad \mbox{in}\quad \re \times \re_+,
   \notag
   \eeq
are
  shown to admit   smooth
{\em compactons}, as oscillatory {\em travelling wave} solutions
with compact support.
 The well known nonnegative compactons, which appeared in
 various applications (first examples by Day, 1998, and Rosenau--Levy, 1999), are
 nonexistent in general  and are not robust relative small perturbations of
 parameters of the PDE.

This is more extended and detailed version of the arXiv preprint
\cite{Gal5NDEarX}. Particularly, essential novelties are available
in \S~5.4, where a family of similarity extensions after blow-up
was detected, which were mentioned but not found in
\cite[\S~5]{Gal5NDEarX}.

\end{abstract}

\maketitle

\section {Introduction: nonlinear dispersion PDEs,
 and main directions of study}
 \label{Sect1}

\subsection{Five main problems and layout:  shocks, rarefaction waves, and
compactons
 for  fifth-order NDEs}

Let us introduce our basic models, which  are five fifth-order
nonlinear dispersion equations (NDEs). These are ordered by
numbers of derivatives inside and outside the quadratic
differential operators involved on the right-hand sides:
 \begin{align}
  & u_t=-u u_{xxxxx} \quad \quad\big(\mbox{NDE--$(5,0)$}\big), \label{N50}\\
 & u_t=-(u u_{xxxx})_x \quad \big(\mbox{NDE--$(4,1)$}\big), \label{N41}\\
& u_t=-(u u_{xxx})_{xx} \quad \big(\mbox{NDE--$(3,2)$}\big),
\label{N32}\\
 & u_t=-(u u_{xx})_{xxx} \quad \big(\mbox{NDE--$(2,3)$}\big), \label{N23}\\
   & u_t=-(u u_{x})_{xxxx} \quad \big(\mbox{NDE--$(1,4)$}\big). \label{N14}
 \end{align}
 The only fully divergent operator is in the last NDE--$(1,4)$
that, being written as
 \beq
 \label{N05}
  \mbox{$
u_t=-(u u_{x})_{xxxx} \equiv -\frac 12\, (u^2)_{xxxxx}
 \quad \big(\mbox{NDE--$(1,4)$ $=$ NDE--$(0,5)$}\big),
 $}
 \eeq
 becomes also the NDE--$(0,5)$, or simply the NDE--5. This
 completes the list of such quasilinear degenerate PDEs
  under consideration.

  \ssk

The main feature of these degenerate odd-order PDEs is that they
admit {\em shock} and {\em rarefaction waves}, similarly to
 the first-order conservation laws such as
{\em Euler's equation}
 \beq
 \label{3}
  u_t + uu_x=0
\quad \mbox{in} \quad \re\times \re_+, \quad u(x,0)=u_0(x) \quad
\mbox{in} \quad \re.
 \eeq
Before explaining the physical significance of the NDEs and their
role in general PDE theory, we pose  four main problems for the
above NDEs (the same as for (\ref{3})):

\ssk

\noi{\bf (I)}  \underline{\em Problem ``Blow-up to $S_-$"}
(Section \ref{Sect2}): {\em to show that the shock of the shape
$-{\rm sign}\, x$ can be obtained by blow-up limit from  a smooth
self-similar solution $u_-(x,t)$ of $(\ref{N50})$--$(\ref{N14})$
in $\re \times (0,T)$, i.e., the following holds:}
 \beq
 \label{con32}
 u_-(x,t) \to S_-(x)= -{\rm sign}\, x =  \left\{
 \begin{matrix}
 \,\,\,\,1 \,\,\, \mbox{\em for} \,\,\, x<0,\\
-1 \,\,\, \mbox{\em for} \,\,\, x > 0,
 \end{matrix}
 \right.
 \quad \mbox{\em as} \quad t \to T^- \quad
 \mbox{\em in \,\, $L^1_{\rm loc}(\re)$}.
  \eeq

\ssk

\noi{\bf (II)} \underline{\em The Riemann Problem $S_+$ $($RP$+)$}
(Section \ref{Sect3}): {\em to show that the initial shock}
   \beq
  \label{Ri12}
  S_+(x) = {\rm sign}\, x =  \left\{
 \begin{matrix}
 -1 \,\,\, \mbox{\em for} \,\,\, x<0,\\
\,\,\,\, 1 \,\,\, \mbox{\em for} \,\,\, x > 0,
 \end{matrix}
 \right.
  \eeq
 {\em for NDEs $(\ref{N50})$--$(\ref{N14})$
 generates  a  ``rarefaction wave", which is  $C^\iy$-smooth
  for $t>0$.}

 \ssk

 \noi{\bf (III)} \underline{\em The Riemann Problem $S_-$ $($RP$-)$}
 (Section \ref{SS1}): {\em introducing a ``$\d$-entropy test" (smoothing of
 discontinuous solutions at shocks via a
``$\d$-deformation"),
    to show
that}
\beq
 \label{rr1}
 S_-(x) \,\,\, \mbox{\em is an ``$\d$-entropy" shock wave, and} \,\,\,
 S_+(x) \,\,\, \mbox{\em is not}.
  \eeq


\ssk

\noi{\bf (IV)} \underline{\em Problem: nonuniqueness/entropy}
(Section \ref{SNonU}):
  {\em to show that a single point ``gradient catastrophe" for the
  NDE $(\ref{N14})$ leads to the principal nonuniqueness of a
  shock wave extension after singularity.} This also suggests {\em nonexistence} of any
   proper
  entropy mechanism for choosing any  ``right" solution after single point blow-up.


 \ssk

In Section \ref{SCK1}, we  discuss  these problems in application
to other NDEs including the following rather unusual one:
 \beq
  \label{kk1}
  u_{ttttt}=(uu_x)_{xxxx},
   \eeq
   which indeed  can be reduced to a first-order system that, nevertheless, {\em is not hyperbolic}, so
   that  modern advanced theory of 1D hyperbolic systems (see e.g.,
    Bressan \cite{Bres} or Dafermos \cite{Daf}) does not apply.
The main convenient mathematical feature of (\ref{kk1}) is that it
is in the {\em normal form}, so it obeys the Cauchy--Kovalevskaya
theorem that guarantees local existence of a unique analytic
solution and makes easier application of our $\d$-entropy
(smoothing) test. Regardless this, (\ref{kk1}) is shown to create
in finite time shocks of the type $S_-(x)$ in (\ref{con32}) and
rarefaction waves for other discontinuous data $\sim S_+(x)$ in
(\ref{Ri12}).

\ssk

Finally, we consider the last:


\ssk

 \noi {\bf (V)} \underline{\em Problem ``Oscillatory Smooth
Compactons"} (Section \ref{Sect6}): {\em to show that
 the  perturbed version
 of the NDE $(\ref{N14})$, as a typical example},
  \beq
  \label{z1}
  u_t=-(|u|u_x)_{xxxx}+ |u| u_x \quad \mbox{in} \quad \re \times \re_+,
   \eeq
{\em admits compactly supported travelling wave (TW) solutions of
{changing sign} near finite interfaces}.
 Equation (\ref{z1})
 is written for solutions with infinitely many sign changes, by replacing
   $u^2$ by the monotone function $|u|u$.

 Nonnegative compact structures have been known since beginning of
the 1990s as {\em compactons} (Rosenau--Hyman, 1993,
\cite{RosH93}). We show that more standard in literature {\em
nonnegative compactons} of fifth-order NDEs such as (\ref{z1}) are
nonexistent in general, and, moreover,  these {\em are not robust}
(not ``{structurally stable}"), i.e., do not exhibit continuous
dependence upon the parameters of PDEs (say, arbitrarily small
perturbations of nonlinearities).

 \subsection{A link to classic entropy shocks for conservation laws}

Indeed,  the above problems {\bf (I)}--{\bf (III)} are classic for
entropy theory of 1D conservation laws from the 1950s.
 It is well recognized that shock waves first appeared in gas dynamics that led to
   mathematical theory of
 entropy solutions of
 the first-order conservation laws and
{Euler's equation}
 (\ref{3}) as a key representative.
The entropy theory for PDEs such as (\ref{3}), with arbitrary
measurable initial data $u_0$, was created by Oleinik \cite{Ol1,
Ol59} and Kruzhkov \cite{Kru2} (analogous scalar  equations in
$\ren$) in the 1950--60s; see details on the history, main
results, and modern developments in the well-known monographs
\cite{Bres, Daf, Sm}. Note that first analysis of the formation of
shocks for (\ref{3}) was performed by Riemann in 1858 \cite{Ri58};
see further details and the history in \cite{Chr07}. It is worth
mentioning that the implicitly given  solution $u=u(x,t)$ of the
Cauchy problem \ef{3}, via the characteristic
 formula
  $$
  u=u_0(x-u \, t),
   $$
containing the key
 wave ``overturning" effect, was obtained earlier by Poisson in 1808
\cite{Poi08}; see \cite{Pom08}.



According to entropy theory  for conservation laws such as
(\ref{3}), it is well-known that  (\ref{rr1}) holds.
This means that
 \beq
 \label{rr41}
 u_-(x,t) \equiv S_-(x)=-{\rm sign} \, x
  \eeq
is the unique entropy solution of the PDE (\ref{3}) with the same
initial data $S_-(x)$. On the contrary, taking $S_+$-type initial
data (\ref{Ri12}) in the Cauchy problem (\ref{3})
  yields the continuous  {\em rarefaction wave}
    with a simple similarity piece-wise linear
    structure,
   \beq
   \label{rr6}
   \mbox{$
    u_0(x)=S_+(x)={\rm sign} \, x \,\,\, \Longrightarrow \,\,\,
   u_{+}(x,t)= g(\frac xt) = \left\{ \begin{matrix}
    -1 \quad \mbox{for} \,\,\, x<-t, \cr
 \,\,\,\frac xt \quad  \mbox{for} \,\,\, |x|<t, \cr
  1 \quad \,  \mbox{for} \,\,\, x>t.
 \end{matrix}
   \right.
    $}
 \eeq

 Our first goal is to justify the same conclusions for the fifth-order NDEs,
 where, of course, the rarefaction wave in the RP$+$ is supposed to
 be
 different from that in (\ref{rr6}).




\ssk

We now return to main applications of the NDEs.

\subsection{NDEs from theory of integrable PDEs and water waves}

 Talking about odd-order PDEs under consideration,
  these naturally appear in classic theory of integrable
 PDEs from shallow water applications,  beginning
with
  the {\em KdV equation},
 \beq
 \label{31}
 u_t+uu_x=u_{xxx},
 \eeq
 the {\em fifth-order KdV equation},
 $$
 u_t + u_{xxxxx} + 30 \, u^2 u_x + 20 \, u_x u_{xx} + 10 \, u
 u_{xxx}=0,
 $$
  and
others.
 These are {\em semilinear}
dispersion equations, which being endowed with smooth
semigroups (groups), generate smooth flows, so discontinuous weak
solutions
 are unlikely, though strong oscillatory behaviour of solutions is typical; see references in \cite[Ch.~4]{GSVR}.

The situation is changed for the quasilinear case. In particular,
consider  the quasilinear
  {\em Harry Dym}
 {\em equation}
  \beq
 \label{HD0}
  u_t = u^3 u_{xxx} \, ,
  \eeq
  which 
   is
one of the most exotic integrable soliton equations; see
\cite[\S~4.7]{GSVR} for survey and references
 therein. Here, (\ref{HD0})  indeed belongs to the NDE family,
 though it seems proper semigroups of its discontinuous solutions
 (if any)
 have never been examined.
 On the other hand,
 moving blow-up singularities and other types of complex
 singularities of the {\em modified Harry Dym equation},
  $$
   \mbox{$
   u_t= u^3 u_{xxx}- u_x- \frac 12\, u^3,
 $}
  $$
  have been described in \cite{Cost06} by delicate asymptotic expansion
  techniques.

    In addition,
integrable equation theory produced various hierarchies of
quasilinear higher-order NDEs, such as
 the fifth-order {\em
Kawamoto equation} \cite{Kaw85}, as a typical example
  \beq
  \label{Kaw111}
 u_t = u^5 u_{xxxxx}+ 5 \,u^4 u_x u_{xxxx}+ 10 \,u^5 u_{xx}
 u_{xxx}.
 \eeq

 We can enlarge  this list talking about possible
   quasilinear  extensions of the integrable
 {\em Lax's
 seventh-order KdV equation}
 $$
 u_t+[35 u^4 + 70(u^2 u_{xx}+ u(u_x)^2)+7(2 u u_{xxxx}+3 (u_{xx})^2 +
 4 u_x u_{xxx})+u_{xxxxxx}]_x=0,
  $$
  and the {\em seventh-order Sawada--Kotara
  equation}
 $$
 u_t+[63 u^4 + 63(2u^2 u_{xx}+ u(u_x)^2)+21( u u_{xxxx}+ (u_{xx})^2 +
  u_x u_{xxx})+u_{xxxxxx}]_x=0;
  $$
see references in \cite[p.~234]{GSVR}.

  The  modern  mathematical theory of odd-order quasilinear PDEs is partially
originated and continues to be strongly connected with the class
of integrable equations. Special advantages of integrability by
using the
inverse scattering transform method, Lax pairs, Liouville
transformations,
 and other explicit algebraic manipulations have made it
possible to create a rather complete theory for some of these
difficult quasilinear PDEs. Nowadays,  well-developed theory
 and  most of rigorous results on existence,
uniqueness, and various singularity and non-differentiability
properties are associated with NDE-type integrable models such as
{\em Fuchssteiner--Fokas--Camassa--Holm} (FFCH) {\em
 equation}
 \beq
  \label{R1}
  \mbox{$
  (I-D_x^2)u_t=
 - 3u u_x+ 2u_x u_{xx} + u u_{xxx}
  \equiv -(I-D_x^2)(u u_x)- \big[u^2 + \frac 12\, (u_x)^2\big].
  $}
\eeq
  Equation (\ref{R1})
 is   an asymptotic model  describing  the wave dynamics
at the free surface of fluids under gravity. It is derived from
Euler equations for inviscid fluids under the long wave
asymptotics of shallow water behaviour  (where the function $u$ is
the height of the water above a flat bottom).
 Applying to (\ref{R1}) the integral operator $(I-D_x^2)^{-1}$
 with the $L^2$-kernel $\o(s)= \frac 12 \, {\mathrm
 e}^{-|s|}>0$, reduces it, for a class of solutions, to the
conservation law (\ref{3}) with a compact {\em first-order}
perturbation,
 \beq
 \label{Con.eq}
  \mbox{$
 u_t  + u u_x=  - \big[ \o*\big(u^2 + \frac 12  (u_x)^2 \big)
 \big]_x.
  $}
 \eeq
 Almost all mathematical results (including entropy
 inequalities and Oleinik's condition (E)) have been obtained by
 using this integral representation of the FFCH equation;
 see the long list of references given in \cite[p.~232]{GSVR}.

There is  another  integrable PDE from the family with third-order
quadratic operators,
 \beq
 \label{ff1SS}
 u_t-u_{xxt}= \a u u_x + \b u_x u_{xx} + u u_{xxx} \quad (\a, \,
 \b \in \re),
  \eeq
  where $\a=-3$ and $\b=2$ yields the FFCH equation (\ref{R1}).
This is  the {\em Degasperis--Procesi} (DP) {\em equation} for
another choice $\a=-4$ and $\b=3$:
 \beq
 \label{DP0}
  \mbox{$
 u_t - u_{xxt}= -4 u u_x+ 3u_x u_{xx} + u u_{xxx},
 \quad \mbox{or} \quad
u_t+u u_x= - \big[ \o*\big(\frac 32 \, u^2 \big)
 \big]_x.
  $}
 \eeq
On existence, uniqueness (of entropy solutions in $L^1 \cap BV$),
parabolic $\e$-regularization, Oleinik's entropy estimate, and
generalized PDEs, see \cite{Coc06}.

 Note that, since the non-local
term in the DP equation (\ref{DP0}) does not contain $u_x$, the
differential properties of its solutions are distinct from those
for the FFCH one (\ref{Con.eq}). Namely, the solutions are less
regular, and  (\ref{DP0}) admits {\em shock waves}, e.g., of the
form
 $$
  \mbox{$
  u_{\rm shock}(x,t)=- \frac 1t \, {\rm sign} \, x \, {\mathrm
  e}^{-|x|},
   $}
   $$
 with rather standard (induced by (\ref{3}))
 but more involved entropy theory; see \cite{Land07, Esch07}.

 Besides (\ref{R1}) and (\ref{DP0}), the family
(\ref{ff1SS}) does not contain other integrable entries.
 A  list of
more applied papers related to various NDEs is also available in
\cite[Ch.~4]{GSVR}.


\subsection{NDEs from compacton theory}

Other important applications of odd-order PDEs are associated with
{\em compacton phenomena} for more general non-integrable models.
For instance,
  the   {\em Rosenau--Hyman} (RH)
{\em equation}
  \beq
  \label{Comp.4}
  \mbox{$
  u_t =  (u^2)_{xxx} + (u^2)_x
  $}
  \eeq
   has special important applications as a widely used model of
   the effects of {nonlinear dispersion} in the pattern
  formation in liquid drops \cite{RosH93}. It is
   the $K(2,2)$ equation from the general $K(m,n)$ family of
   the following NDEs:
  \beq
  \label{Comp.5}
   u_t =  (u^n)_{xxx} +  (u^m)_x \quad (u \ge 0),
   \eeq
   that
 describe   phenomena of compact pattern
 formation, \cite{RosCom94, Ros96}. Such PDEs also appear in curve motion and shortening
flows \cite{Ros00}.
 Similar to well-known parabolic models of
  the porous medium type, the $K(m,n)$ equation (\ref{Comp.5}) with $n>1$
  is degenerate at $u=0$, and
 therefore may exhibit finite speed of propagation and admit solutions with finite
 interfaces.
 The crucial advantage of the RH equation
 (\ref{Comp.4}) is that it possesses
  {\em explicit}  moving compactly supported
 soliton-type solutions, called {\em compactons}
 \cite{RosH93, RosCom94}, which are {\em travelling wave} (TW) solutions to
 be discussed for the PDEs under consideration.


   Various families of quasilinear
third-order KdV-type equations can be found  in \cite{6R.1}, where
further references concerning such PDEs and their exact solutions
 are given.
   Higher-order generalized KdV
equations are of increasing interest; see {e.g.,} the quintic KdV
equation in \cite{Ros98}, and also \cite{Yao7}, where the
seventh-order PDEs are studied.

 More general $B(m,k)$ equations,
 $$
 u_t+ a(u^m)_x = \mu (u^k)_{xxx},
 $$
  which coincide with the
 $K(m,k)$
  after scaling,
 also admit simple semi-compacton solutions \cite{RosK83}. The same is true for
   the $Kq(m,\o)$ nonlinear dispersion equation
    (another nonlinear extension of the KdV) \cite{RosCom94}
 $$
 u_t + (u^m)_x + [u^{1-\o}(u^\o u_x)_x]_x =0.
 $$
Setting $m=2$ and $\o=\frac 12$ yields a typical quadratic PDE
 \beq
  \label{BB1SS}
{\bf B}(u) \equiv  u_t + (u^2)_x + u u_{xxx} + 2 u_x u_{xx}=0.
  \eeq
 It is curious that (\ref{BB1SS}) admits an extended compacton-like
 dynamics
 on a standard trigonometric-exponential
   subspaces, on which
    \beq
    \label{BB2}
     \mbox{$
    u(x,t)= C_0(t) + C_1(t) \cos \l x + C_2(t)  \sin \l x \in W_3= {\rm
    Span}\{1,\cos \l x, \sin \l x\},
     $}
    \eeq
    where $  \l= \sqrt{\frac 23}$.
   This subspace is invariant under the quadratic operator ${\bf
   B}$ in the usual sense that
     ${\bf B}(W_3) \subseteq
    W_3$. Therefore substituting (\ref{BB2})
     into the PDE
    \ref{BB1SS}) yields for the expansion coefficients on $W_3$
     $\{C_0,C_1,C_2\}$  a 3D  nonlinear dynamical
   system; see further such examples of exact solutions of NDEs on invariant subspaces
   in \cite[Ch.~4]{GSVR}.

 Combining the $K(m,n)$ and $B(m,k)$ equations gives the
 dispersive-dissipativity entity $DD(k,m,n)$ \cite{Ros98DD}
  $$
  u_t + a(u^m)_x + (u^n)_{xxx} = \mu (u^k)_{xx},
  $$
  which can also admit solutions on invariant subspaces for
  some values of parameters.

For the fifth-order NDEs,
 such as
  \beq
  \label{Comp.3}
  u_t =  \a (u^2)_{xxxxx} + \b (u^2)_{xxx} + \g (u^2)_x \quad
  \mbox{in} \,\,\, \re \times \re_+,
  \eeq
compacton solutions  were first constructed in \cite{Dey98}, where
the  more general $K(m,n,p)$ family of
PDEs,
 $$ 
 u_t+ \b_1 (u^m)_x + \b_2 (u^n)_{xxx} + \b_3 D_x^5 (u^p)=0,
 $$
  with
 $m,\,n,\,p>1$,
  was introduced. Some of these equations will be treated later on.
Equation  (\ref{Comp.3}) is also
  associated with the
 family $Q(l,m,n)$ of more general quintic  evolution PDEs with
nonlinear dispersion,
  \beq
  \label{qq1.1}
  \mbox{$
u_t + a (u^{m+1})_{x} + \o \bigl[u(u^n)_{xx}\bigr]_x + \d
\bigl[u(u^l)_{xxxx}\bigr]_x=0,
 $}
 \eeq
 possessing multi-hump, compact solitary solutions \cite{Ros599}.

Concerning higher-order in time quasilinear PDEs, let us mention
a generalization of the {\em combined dissipative
double-dispersive} (CDDD) {\em equation} (see, {e.g.,}
\cite{Por02})
 \beq
 \label{Por.1}
 u_{tt}= \a u_{xxxx} + \b u_{xxtt} + 
   \g
 (u^2)_{xxxxt} + \d (u^2)_{xxt} + \e (u^2)_t,
 \eeq
and also the {\em nonlinear modified dispersive Klein--Gordon
equation} ($mKG(1,n,k)$),
 \beq
 \label{mkg1}
 u_{tt}+ a (u^n)_{xx}+ b(u^k)_{xxxx}=0, \quad n,k>1 \quad (u \ge 0);
  \eeq
 see some exact TW solutions in \cite{Inc07}.
 For $b>0$, (\ref{mkg1}) is of hyperbolic (or Boussinesq) type in
 the class of nonnegative solutions.
 We also mention related 2D {\em dispersive
Boussinesq equations} denoted by $B(m,n,k,p)$ \cite{Yan03},
 $$ 
  (u^m)_{tt} + \a (u^n)_{xx} + \b (u^k)_{xxxx} + \g
  (u^p)_{yyyy}=0 \quad \mbox{in}
  \quad \re^2 \times \re.
  $$ 
See \cite[Ch.~4-6]{GSVR} for more references and examples of exact
solutions on invariant subspaces of NDEs of various types and
orders.


\subsection{On canonical third-order NDEs}

Until recently,  quite a little was  known about proper
mathematics concerning discontinuous solutions, rarefaction waves,
and ``entropy-like" approaches, even for the simplest third-order
NDEs such as (\ref{Comp.4}) or (see
\cite{GPnde, GPndeII})
 \beq
 \label{1}
 u_t= (uu_x)_{xx}.
  \eeq
However,  the smoothing results for sufficiently regular solutions
of linear and nonlinear third-order PDEs are well know from the
1980-90s. For instance, infinite smoothing results were proved in
\cite{Cr90}  (see also \cite{Hos99}) for the general linear
equation
 \beq
 \label{lin11}
 u_t+a(x,t) u_{xxx}=0 \quad (a(x,t) \ge c>0),
 \eeq
 and in \cite{Cr92} for the corresponding fully nonlinear PDE
  \beq
  \label{Lin12}
  u_t+f(u_{xxx},u_{xx},u_x,u,x,t)=0 \quad \big(\,f_{u_{xxx}} \ge
  c>0\, \big);
   \eeq
   see also \cite{Cai97} for semilinear equations.
   Namely, for a class of such equations, it is shown that, for data
   with minimal regularity and sufficient (say, exponential) decay at infinity, there
   exists a unique
   solution $u(x,t) \in C^\infty_x$ for small $t>0$.
 Similar smoothing local in time results for unique solutions
 are available for
 equations in $\re^2$,
  \beq
  \label{lin13}
  u_t + f(D^3u,D^2u,Du,u,x,y,t)=0;
   \eeq
   see \cite{Lev01} and further  references therein.

These smoothing results have been used in \cite{GPndeII} for
developing  a kind of a $\d$-entropy test  for discontinuous
solutions by using techniques of smooth deformations. We will
follow these ideas applied now to shock and compacton solutions of
higher-order NDEs and others.

\section{{\bf (I) Problem ``Blow-up"}: existence of shock $S_-$ similarity
solutions}
 \label{Sect2}

We now show that  Problem {\bf (I)} on blowing up to the shock
$S_-(x)$ can be solved in a unified manner by constructing
self-similar solutions. As often happens in nonlinear evolution
PDEs, the refined structure of such bounded and generic shocks
 is described in a scaling-invariant  manner.

\subsection{Finite time blow-up formation of the shock wave
$S_-(x)$}
 \label{Sect3.1}

One can see that all five  NDEs
 (\ref{N50})--(\ref{N14}) admit the following
similarity substitution:
 \beq
 \label{2.1}
 u_-(x,t)=g(z), \quad z= x/(-t)^{\frac 15} \quad (t<0),
  \eeq
where,  by translation, the blow-up time in reduces to $T=0$.
Substituting (\ref{2.1}) into the NDEs yields for $g$  the
following ODEs in $\re$, respectively:
 \begin{align}
 & gg^{(5)}=- \textstyle{\frac 15}\, g'z, \label{1E}\\
& (gg^{(4)})'=- \textstyle{\frac 15}\, g'z, \label{2E}\\
 & (gg''')''=- \textstyle{\frac 15}\, g'z, \label{3E}\\
 & (gg'')'''=- \textstyle{\frac 15}\, g'z, \label{4E}\\
 & (gg')^{(4)}=-\textstyle{ \frac 15}\, g'z, \label{5E}
  \end{align}
  with the following conditions at infinity for the shocks $S_-$:
 \beq
 \label{2.2}
  \mbox{$
   g(\mp
  \infty)=\pm 1.
   $}
   \eeq
    In view of the symmetry of the ODEs,
     \beq
     \label{symm88}
      \left\{
      \begin{matrix}
     g \mapsto -g,\\
     z \mapsto -z,
      \end{matrix}
      \right.
      \eeq
it suffices to get
 odd solutions for $z<0$ posing anti-symmetry
conditions at the origin,
 \beq
 \label{2.4}
 g(0)=g''(0)=g^{(4)}(0)=0.
  \eeq

\subsection{Shock similarity profiles exist and are unique: numerical results}

Before performing a rigorous approach to Problem (I), it is
convenient and inspiring to check whether the shock similarity
profiles $g(z)$ announced in (\ref{2.1}) actually exist and are
unique for each of the ODEs (\ref{1E})--(\ref{5E}). This is done
by numerical methods that supply us with positive and convincing
conclusions. Moreover, these numerics clarify some crucial
properties of profiles, which will determine    the actual
strategy of rigorous study.

\begin{figure}
\centering
\includegraphics[scale=0.70]{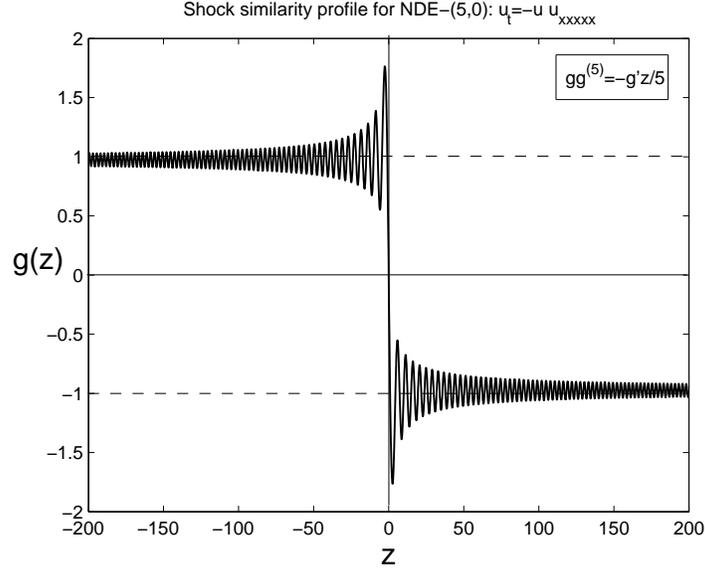}
 \vskip -.3cm
\caption{\small The shock similarity profile $g(z)$ as the unique
solution of the problem (\ref{1E}), (\ref{2.4}); $z \in
[-200,200]$.}
\label{F1}
\end{figure}

\begin{figure}
\centering
\includegraphics[scale=0.70]{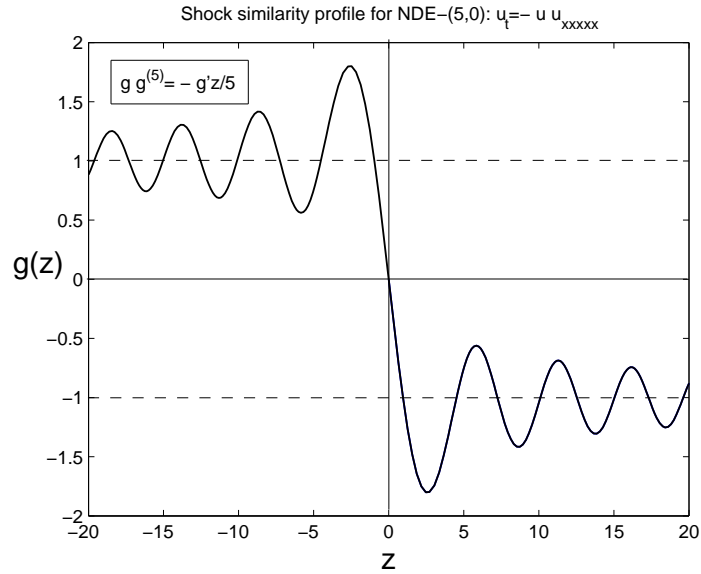}
 \vskip -.3cm
\caption{\small The shock similarity profile $g(z)$ as the unique
solution of the problem (\ref{1E}), (\ref{2.4}); $z \in
[-20,20]$.}
\label{F1NN}
\end{figure}

A typical structure of this shock similarity profile  $g(z)$
satisfying (\ref{1E}), (\ref{2.4})  is shown in Figure \ref{F1}.
As a key feature, we observe a highly oscillatory behaviour of
$g(z)$ about $\pm 1$ as $z \to \mp\infty$,  that can essentially
affect the metric of the announced  convergence (\ref{con32}).
Therefore,  we will need to describe this oscillatory behaviour in
detail. In Figure \ref{F1NN}, we show the same profile $g(z)$ for
smaller $z$. It is crucial that, in all numerical experiments, we
obtained the same profile that indicates that it is the unique
solution of
 (\ref{1E}), (\ref{2.4}).

\begin{figure}
\centering
\subfigure[equation (\ref{2E})]{
\includegraphics[scale=0.52]{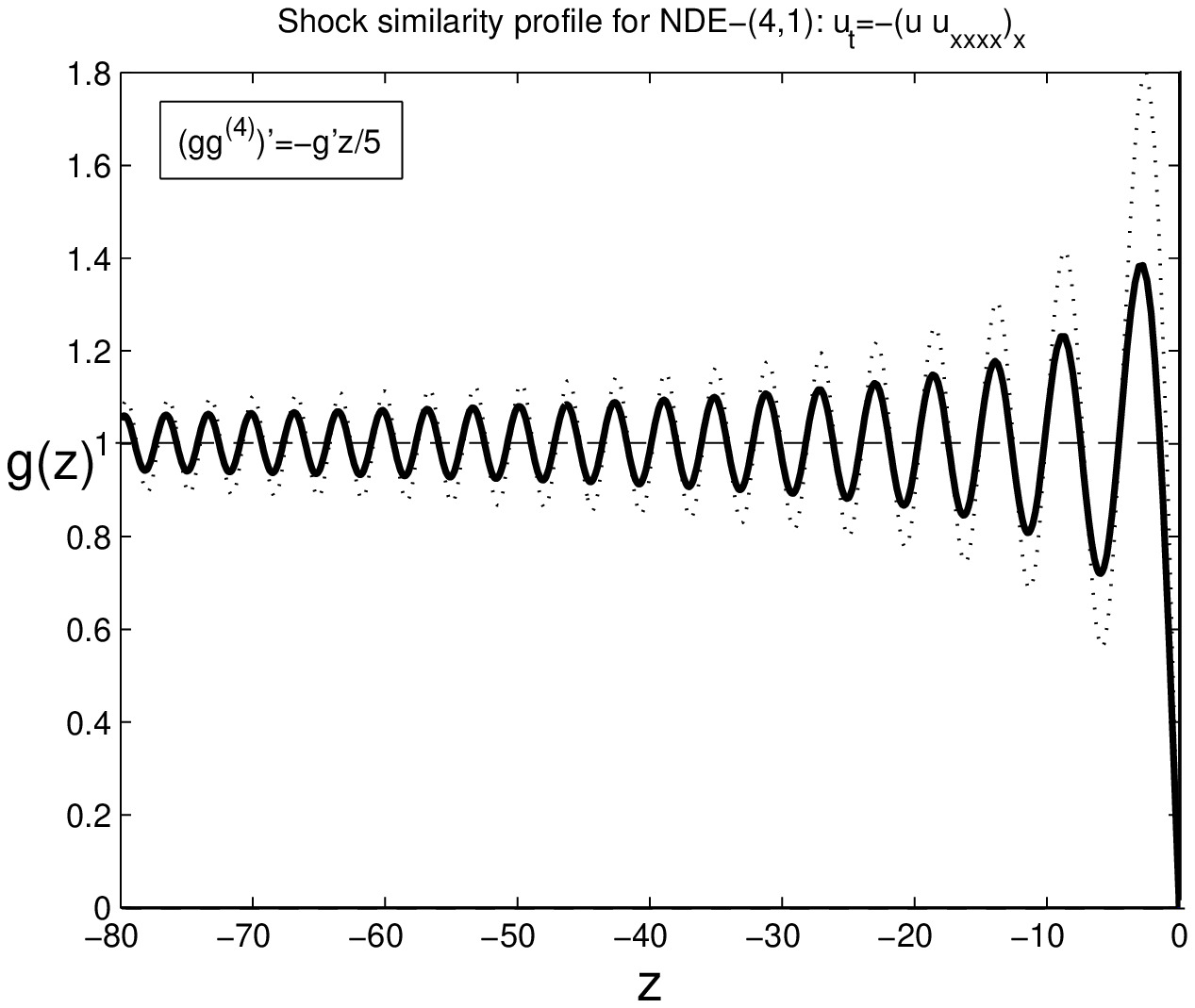} 
}
\subfigure[equation (\ref{3E})]{
\includegraphics[scale=0.52]{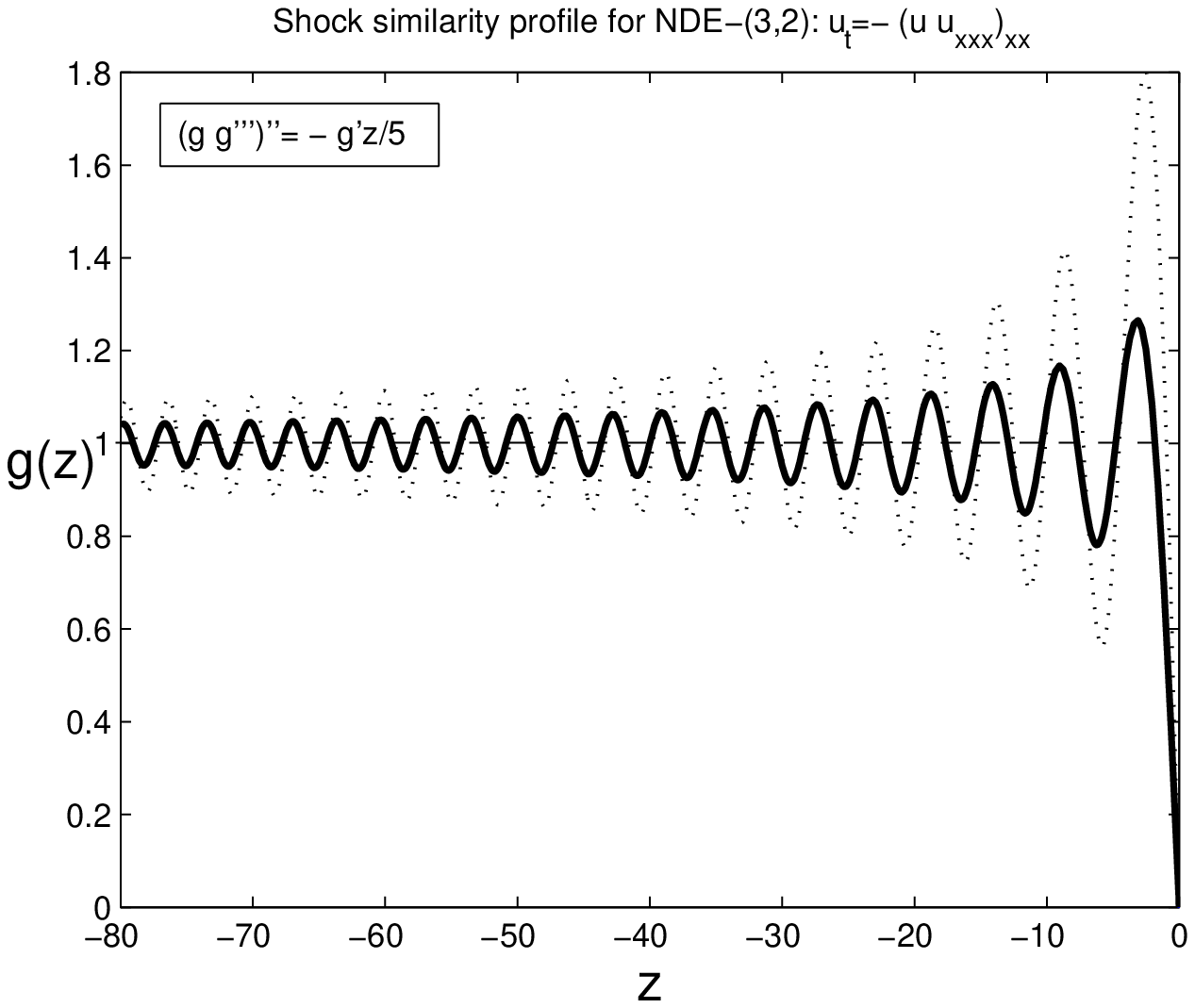} 
}
 \subfigure[equation (\ref{4E})]{
\includegraphics[scale=0.52]{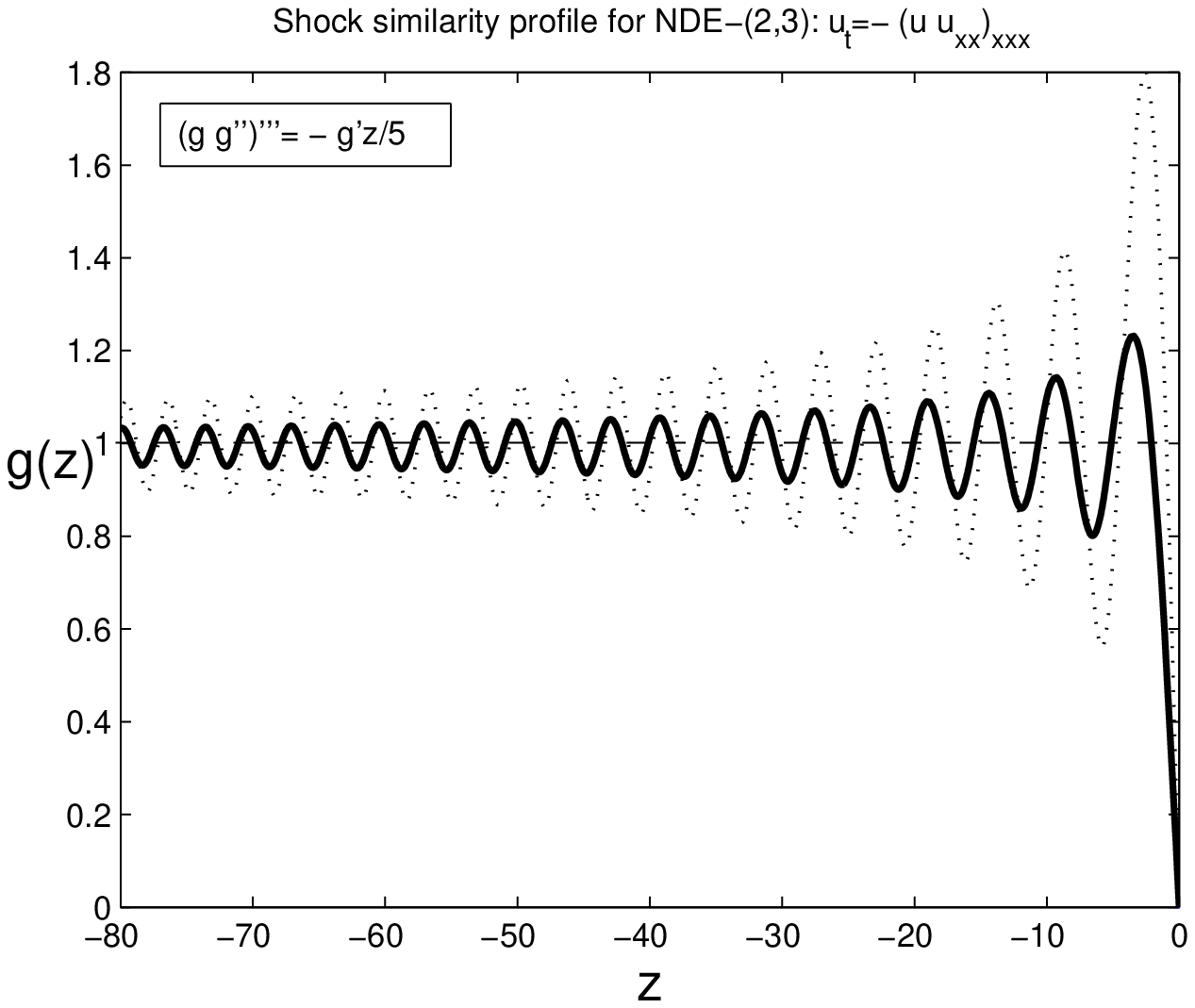} 
}
\subfigure[equation (\ref{5E})]{
\includegraphics[scale=0.52]{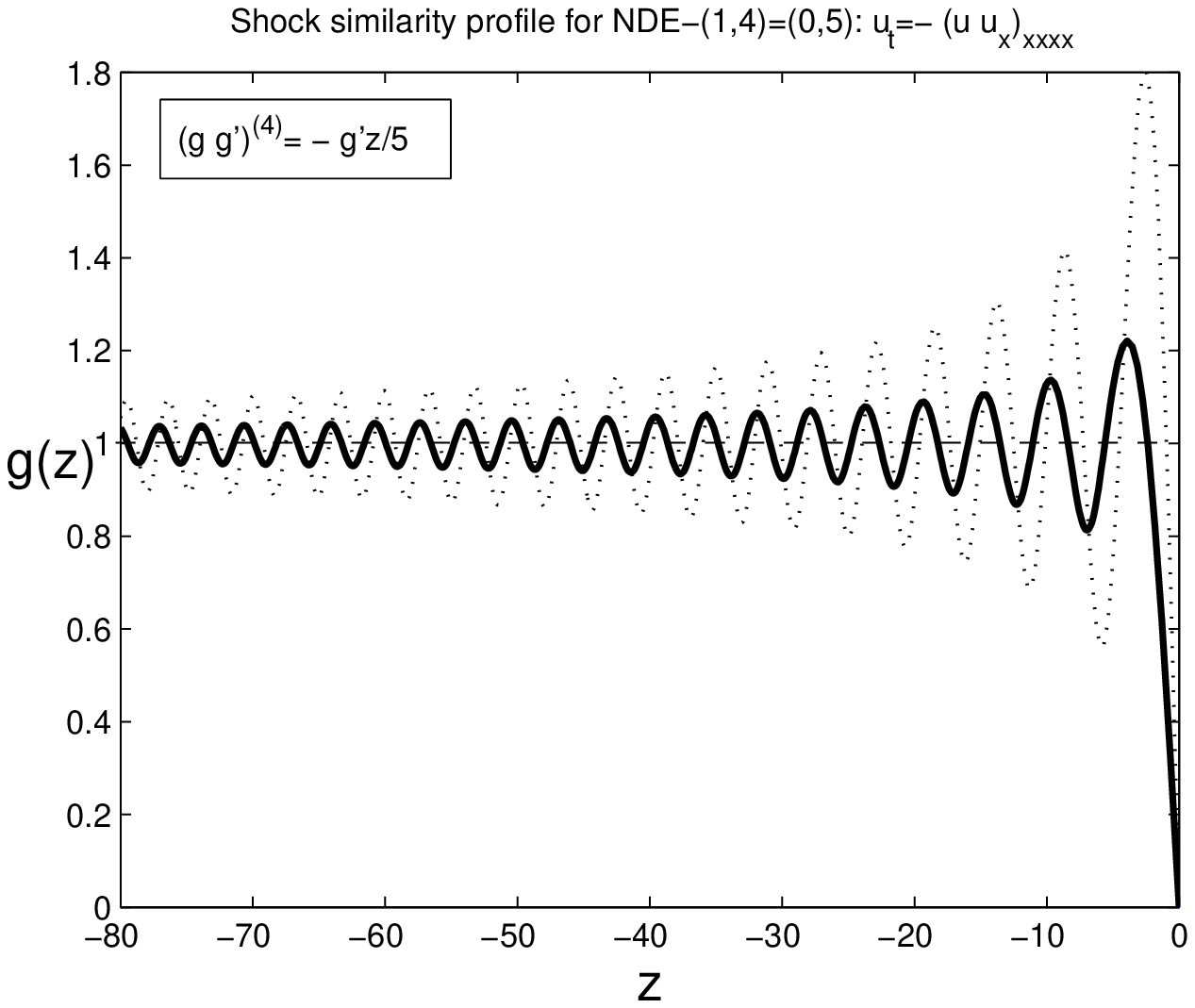} 
}
 \vskip -.3cm
\caption{\rm\small Shock similarity profiles as solutions of
(\ref{2E})--(\ref{5E}), (\ref{2.4}) respectively. For comparison,
dotted lines denote the profile from Figures
  \ref{F1} and \ref{F1NN}.}
 \label{F7}
\end{figure}

Figure \ref{F7}(a)--(d) show the shock similarity profiles for the
rest of NDEs (\ref{N41})--(\ref{N14}). They differ from each other
rather slightly.


\ssk

\noi{\bf Remark: on regularization in numerical methods.}
 For the fifth-order NDEs,
this and further numerical constructions are performed by {\tt
MatLab} by using
the {\tt bvp4c} solver. Typically, we take the relative and
absolute tolerances
 \beq
 \label{T1NN}
 {\rm Tols}=10^{-4}.
  \eeq
Instead of the degenerate ODE (\ref{1E}) (or others), we  solve
the regularized equation
 \beq
 \label{T2}
  \mbox{$
  g^{(5)} =-\frac {{\rm sign}\, g}{\sqrt{\nu^2+g^2}} \, \big(\frac 15 \, g'z\big),
 \quad \mbox{with the regularization parameter $\nu=10^{-4}$},
 $}
  \eeq
  where the choice of small $\nu$ is coherent with the
  tolerances in (\ref{T1NN}). Sometimes, we will need to use the
  enhanced parameters ${\rm Tols} = \nu = 10^{-7}$ or even
  $\sim 10^{-9}$.

\subsection{Justification of oscillatory behaviour
 about equilibria $\pm 1$ and other asymptotics}

Thus, the shock profiles $g(z)$ are oscillatory about $+1$ as $z
\to - \infty$. In order to describe these oscillations in detail,
we linearize all the ODEs (\ref{1E})--(\ref{5E}) about the regular
equilibrium $g(z) \equiv 1$ by setting $g = 1+\hat g$ to get the
linear ODE
 \beq
 \label{2.5}
  \mbox{$
   {\bf B}_5^*\, \hat g \equiv -
  \hat g^{(5)} -\frac 15 \, \hat g'z=0.
   $}
   \eeq

Note that  this equation reminds us of that for the rescaled
kernel $F(z)$ of the
 fundamental solution of the corresponding linear dispersion
 equation,
 \beq
 \label{Lin.717}
 u_t = -u_{xxxxx} \quad \mbox{in} \,\,\, \,\re \times \re_+.
 \eeq
The fundamental solution of  the corresponding linear
 operator
$\frac{\partial}{\partial t}+D_x^5$ in (\ref{Lin.717}) has the standard similarity form 
 \beq
 \label{bbb.5}
  \mbox{$
 b(x,t)= t^{-\frac 15} F(y), \,\,\,\,\mbox{with} \,\,\,\, y= x/{t^{1/5}} ,
  $}
 \eeq
where $F(y)$ is a unique solution of the ODE problem
 \beq
 \label{FODf.1}
\mbox{$ {\bf B}_5F \equiv -F^{(5)} + \frac 15\,  (F y)'=0 \,\,
\mbox{in} \,\, \re, \,\,\, \int
 F
 =1, \quad \mbox{or $F^{(4)}= \frac 15 \, F y$ on integration}.
 $}
  \eeq
However, the operator ${\bf B}_5$ in (\ref{FODf.1}) is not
identical to that in (\ref{2.5}). Moreover, this ${\bf B}^*_5$ is
adjoint to ${\bf B}_5$ in some indefinite metric and both the
operators possess countable families of eigenfunctions, which
particularly are generalized Hermite polynomials for ${\bf
B}_5^*$. We will not use   this Hermitian spectral theory later
on, so refer to \cite[\S~9]{2mSturm} and \cite[\S~8.2]{GalPet2m}
for further results and applications.

\ssk

Let us return to the linearized ODE (\ref{2.5}). Looking for
possible asymptotics as $z \to -\iy$ yields the following
exponential ones with the characteristic equation:
  \beq
  \label{a111}
   \mbox{$
 \hat  g(z) \sim {\mathrm e}^{a |z|
  ^{5/4}} \quad \Longrightarrow \quad  a^4= \frac {4^4}{5^5}.
  $}
   \eeq
 Finally,  choosing the purely imaginary root of the algebraic equation in
  (\ref{a111}) with ${\rm Re} \, a=0$ gives
 a refined WKBJ-type asymptotics of solutions of (\ref{1E}):
  \beq
  \label{as11}
 \mbox{$
 g(z) = 1+  |z|^{-\frac 58}
\bigl[A \sin\bigl( {a_0}|z|^{\frac 54}\bigr) + B \cos\bigl(
{a_0}|z|^{\frac 54}\bigr)\bigr]+... \,\,\, \mbox{as} \,\,\, z \to
-\iy,
   $}
   \eeq
    where $A$ and $B$ are some real constants satisfying $A^2+B^2 \ne 0$.

 The asymptotic behaviour (\ref{as11})  implies two important conclusions:

 \begin{proposition}
  \label{Pr.Var}
    The shock wave profiles $g(z)$ solving
   $(\ref{1E})$--$(\ref{5E})$, $(\ref{2.2})$ satisfy:  {\rm (i)}
  \beq
   \label{L11}
 \quad g(z)-1 \not \in L^1(\re_-), \quad \mbox{and}
 \eeq

 \noi {\rm (ii)}
  the
  total
 variation of $g(z)$ $($and hence of $u_-(x,t)$ for any $t<0)$
  is {infinite}.
 \end{proposition}

 \noi{\em Proof.}
   Setting $|z|^{\frac 54}=v$ in the integrals below
  yields by (\ref{as11}):
 \beq
 \label{pp1ss}
  \begin{matrix}
 \mbox{$
 {\rm (i)} \quad \int\limits_{-\iy}|g(z)-1|\, {\mathrm d}z \sim
 \int\limits^\iy \frac{|\cos z^{\frac 54}|}{z^{5/8}}\, {\mathrm d}z \sim
  \int\limits^\iy
 \frac{|\cos v|}{v^{7/10}} \, {\mathrm d}v= \iy; \quad
 \mbox{and}\qquad\qquad
 $}\ssk\ssk
 \\
  {\rm (ii)} \quad \mbox{$
 | g(\cdot)|_{\rm Tot.Var.} =
 \int\limits_{-\infty}^{+\infty} |g'(z)|\, {\mathrm d}z
 \sim \int\limits^\iy  \frac{|\cos z^{\frac 54}|}{z^{3/8}}\, {\mathrm
 d}z
 \sim \int\limits^\infty \frac
 {|\cos v|}{\sqrt v}\, {\mathrm d}v= \infty. \qed\qquad\qquad
 $}
  \end{matrix}
  \eeq

\ssk

 This is in striking contrast with the case of conservation laws
 (\ref{3}), where finite total variation approaches and  Helly's second theorem
 (compact embedding of sets of bounded functions of bounded total variations into $L^\infty$)
   used to be  key; see Oleinik's pioneering approach \cite{Ol1}.
  In view of the
  presented  properties of the similarity shock  profile $g(z)$, the
  convergence in (\ref{con32}) takes place
 for any $ x \in \re$, uniformly in $\re \setminus (-\mu,\mu)$, $\mu>0$ small, and
 in $L^p_{\rm loc}(\re)$ for  $p \in[1, \infty)$, that, for convenience, we fix in
 the following:


\begin{proposition}
 \label{Pr.Co}
 For the shock similarity profile $g(z)$
 the
  convergence $(\ref{con32})$ with $T=0$:

 {\rm (i)} does not hold in $L^1(\re)$, and

  {\rm (ii)} does hold in $L^1_{\rm loc}(\re)$, and moreover,
 for any fixed finite $l>0$,
   \beq
   \label{conMM}
   \|u_-(\cdot,t)-S_-(\cdot)\|_{L^1(-l,l)} =O((-t)^{\frac 18})\to 0 \quad \mbox{as}
   \quad t \to 0^-.
    \eeq

  \end{proposition}

\noi{\em Proof} of (\ref{conMM}) is the same as in (\ref{pp1ss})
with a finite interval of integration: for $l=1$,
 $$
 \mbox{$
\|\cdot \|_{L^1(-1,1)} \sim (-t)^{\frac 15}
\int\limits^{(-t)^{-1/5}} z^{-\frac 58}|\cos z^{\frac 54}|\,
{\mathrm
 d}z \sim (-t)^{\frac 15} \int\limits^{(-t)^{-1/4}}
v^{-\frac 7{10}}|\cos v|\, {\mathrm
 d}v \sim (-t)^{\frac 18}. \qed
 $}
 $$

 \ssk


Finally, note that each $g(z)$ has
 a regular  asymptotic expansion near
the origin. For instance, for the first ODE (\ref{1E}), there
exist solutions such that
 \beq
 \label{2.7}
  \mbox{$
 g(z) = C z +D z^3 - \frac 1{600} \, z^5 + \frac D{6300 C}\, z^7+... \, ,
  $}
  \eeq
 where $C<0$ and $D \in \re$ are some constants.
 The local
uniqueness of such  asymptotics is traced out by using Banach's
Contraction Principle
 applied to the equivalent integral equation in the metric
 of $C(-\mu,\mu)$, with $\mu>0$  small.
 Moreover,
   it can be shown that
(\ref{2.7}) is the expansion of an analytic function.
 Other ODEs admit similar local representations of solutions.

 We now need the following scaling invariance of the
ODEs (\ref{1E})--(\ref{5E}): if $g_1(z)$ is a solution, then
 \beq
 \label{2.8}
 \mbox{$
 g_a(z) = a^5 g_1\big(\frac z a\big) \quad \mbox{is a solution for any $a \not =
 0$}.
  $}
  \eeq

\subsection{Existence
of a shock similarity profile}
 \label{S3.4}

Using the asymptotics derived above, we now in a position to prove
the following:

\begin{proposition}
\label{Pr.1}
 The problem  $(\ref{2.2})$, $(\ref{2.4})$ for ODEs
$(\ref{1E})$--$(\ref{5E})$ admits a
   solution $g(z)$, which
is an odd analytic function.
 \end{proposition}

Uniqueness for such higher-order ODEs is a more difficult problem,
which is not studied here, though it has been seen numerically.
Moreover,  there are some analogous results. We refer to the paper
\cite{Gaz06} (to be used later on), where uniqueness of a
fourth-order semilinear ODE was established by an improved
shooting argument.

 Notice another difficult aspect of the problem.  Figures
\ref{F1}--\ref{F7} above, which were obtained by careful numerics,
clearly convince that the positivity holds:
 \beq
 \label{gg11}
 g(z) > 0 \quad \mbox{for} \quad z<0,
  \eeq
  which is also difficult to prove rigorously; see further
 comments below. Actually, (\ref{gg11}) is not that important for
 the key convergence (\ref{con32}), since possible sign changes (if any) disappear in the limit
  as $t \to T^-$.
  It seems that nothing prevents existence of some
 ODEs from the family
 (\ref{1E})--(\ref{5E}), with different nonlinearities, for which the shock profiles can change
 sign for $z<0$.

  \ssk

  \noi{\em Proof.}
    As above, we
  consider the first ODE (\ref{1E}) only. We use a shooting argument using the 2D
  bundle of asymptotics (\ref{2.7}).
 By scaling (\ref{2.8}), we put $C=-1$,  so, actually, we
 deal with the {\em one-parameter shooting problem} with
  the 1D family of orbits
 satisfying
 \beq
 \label{2.7D}
  \mbox{$
 g(z;D) = -z+ D \, z^3 - \frac 1{600} \, z^5 - \frac D{6300} \, z^7+... \, ,
 \quad D \in \re.
  $}
  \eeq

It is not hard to check that, besides stabilization to unstable
constant  equilibria,
 \beq
 \label{2.7D1}
 g(z) \to C_- > 0 \quad \mbox{as}
 \quad z \to - \infty,
  \eeq
  the ODE (\ref{2E}) admits an unbounded stable behaviour given by
   \beq
   \label{g11}
    \mbox{$
   g(z) \sim g_*(z)= - \frac 1{120}\, z^5 \to + \infty\quad \mbox{as}
 \quad z \to - \infty.
  $}
  \eeq
The overall asymptotic bundle about the exact solution $g_*(z)$ is
obtained by linearization: as $z \to -\iy$,
 \beq
 \label{lin12}
  \mbox{$
  g(z)=g_*(z) + Y(z) \quad \Longrightarrow \quad g_* Y^{(5)}=-
  \frac 15 Y'z+... \,\,\, \mbox{or} \,\,\, z^5 Y^{(5)}= 24 Y'
  z+...\, .
   $}
   \eeq
 This is Euler's type homogeneous equation with the characteristic equation
  \beq
  \label{lin13NN}
   Y(z)= z^m \quad \Longrightarrow \quad m_1=0 \,\,(Y_2(z) \equiv 1) \,\,\, \mbox{or}
   \,\,\,
   (m-1)(m-2)(m-3)(m-4)=24.
    \eeq
This yields another $m_2=0$ (hence there exists $Y_2(z) = \ln
|z|$), $m_3=5$ (not suitable), and a proper single  complex root
with ${\rm Re}\, m= \frac 52 <5$ yielding oscillatory
$Y_{3,4}(z)$. Thus:
  \beq
  \label{lin14}
  \mbox{as $z \to -\iy$, there exists a 4D asymptotic bundle about
  $g_*(z)=- \frac 1{120}\, z^5$.}
   \eeq

Therefore, at $z = -\iy$, we are given a 2D bundle of proper
solutions (\ref{as11}), as well as a 4D fast growing profiles from
(\ref{lin14}).
 This determines the strategy of the 1D shooting via the
$D$-family (\ref{2.7D}):

(i) obviously, for all $D \ll -1$, we have that
 $g(z;D)>0$ is {monotone decreasing} and  approaches the stable behaviour
 (\ref{g11}), (\ref{lin14}), and

(ii) on the contrary, for all $D\gg 1$, $g(z;D)$ gets non-monotone
and has a zero at some finite $z_0=z_0(D)<0$, satisfying $z_0(D)
\to 0^-$ as $ D \to + \infty$, and eventually  approaches the
bundle in  (\ref{lin14}), but in an essentially {\em non-monotone}
way.

It follows from different and opposite ``topologies" of the
behaviour announced in (i) and (ii) that there exists a constant
$D_0$ such that $g(z;D_0)$ does not belong to those two sets of
orbits (both are open)  and hence does not approach $g_*(z)$ as $z
\to - \infty$  at all. This is precisely the necessary shock
similarity profile.  $\qed$




\ssk

\begin{figure}
\centering
\includegraphics[scale=0.7]{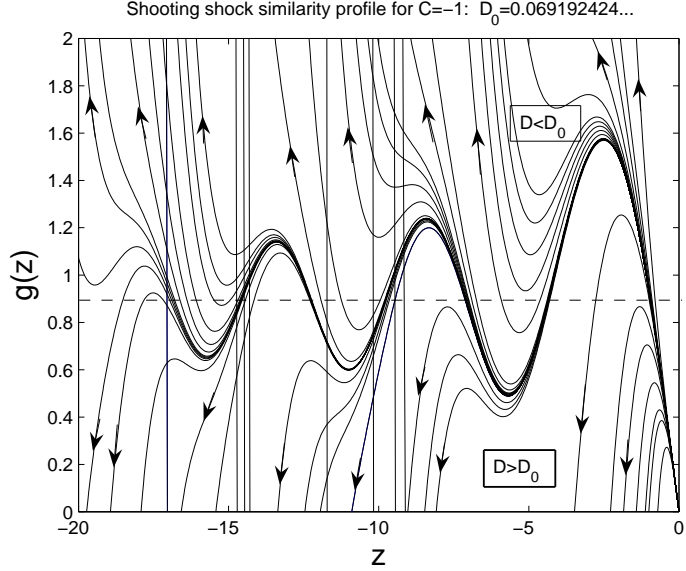}
 \vskip -.3cm
\caption{\small Shooting the shock similarity profile $g(z)$ via
the family (\ref{2.7D}); $D_0=0.069192424...\,$.}
\label{FF1n}
\end{figure}

This 1D shooting approach is explained in Figure \ref{FF1n}
obtained numerically, where
 \beq
 \label{D01}
 D_0=0.069192424... \, .
 \eeq
 It seems that as $D \to D_0^+$, the zero of $g(z;D)$ must disappear
 at infinity, i.e.,
  \beq
  \label{FF12}
  z_0(D) \to - \infty \quad \mbox{as}
  \quad D \to D_0^+,
   \eeq
   and this actually happens as Figure \ref{FF1n} shows. Then this
  would justify the positivity (\ref{gg11}). Unfortunately, in
  general (i.e., for similar ODEs with different sufficiently arbitrary
    nonlinearities),   this is not true, i.e., cannot be guaranteed by a
  topological argument. So that the actual operator structure of
  the ODEs should be involved in the study, so, theoretically, the positivity
  is difficult to guarantee in general. Note again that, if the shock
  similarity profile $g(z)$ had a few zeros for $z<0$, this
  would not affect the crucial convergence property such as
  (\ref{con32}).

\subsection{Self-similar formation of other shocks}

\begin{figure}
\centering
\includegraphics[scale=0.75]{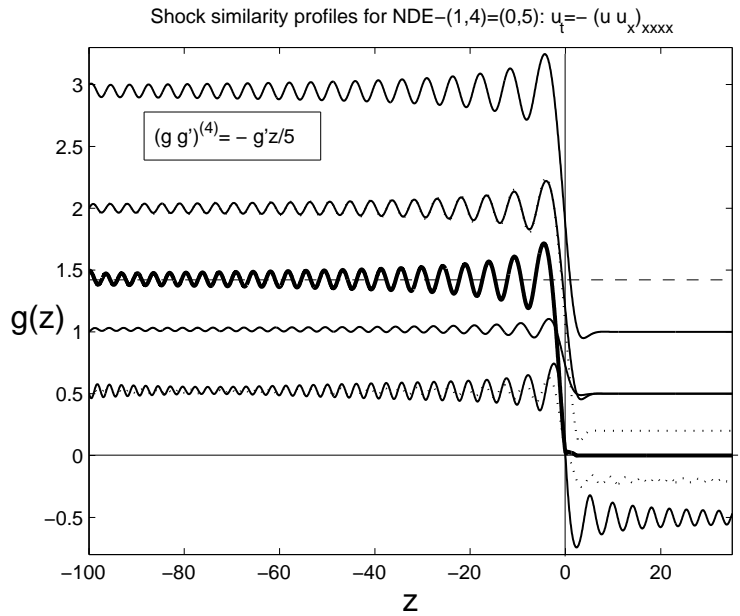}
 \vskip -.3cm
\caption{\small Various shock similarity profiles $g(z)$ as
solutions of the problem (\ref{5E}), (\ref{2.4}).}
\label{F3NNN}
\end{figure}

 \noi \underline{\em NDE--$(1,4)$}.
Let us first briefly consider the last ODE (\ref{5E}) for the
fully divergent NDE (\ref{N14}).
 Similarly,  by the same arguments, we show
that, according to (\ref{2.1}), there exist other non-symmetric
shocks as non-symmetric step-like functions, so that, as $t \to
0^-$,
 \beq
 \label{2.10}
  u_-(x,t) \to
  \left\{
  \begin{matrix}
  C_->0 \,\,\, \mbox{for} \,\,\, x<0, \\
  C_0  \,\,\,\quad \, \quad \mbox{for} \,\,\, x=0, \\
 C_+<0 \,\,\, \mbox{for} \,\,\, x>0,
  \end{matrix}
   \right.
    \eeq
where $C_- \not = -C_+$ and $C_0 \not = 0$.
 Figure \ref{F3NNN} shows a
few of such similarity profiles $g(z)$, where three of these are
strictly positive. The most interesting is the boldface one with
 $$
 C_-=1.4 \quad \mbox{and} \quad C_+=0,
  $$
   which has the finite right-hand interface at
$z=z_0 \approx 5$,
 with the expansion
 \beq
 \label{C001}
   \mbox{$
   g(z)=- \frac {z_0}{4200}\, (z_0-z)_+^4(1+o(1)) \to 0^-\,\,\,
 \mbox{as}\,\,\, z \to z_0.
  $}
   \eeq
It follows that this $g(z) <0$ near the interface so
 the function changes sign there, which
 is also seen in Figure \ref{F3NNN} by carefully checking the shape of profiles above the
 boldface one with the finite interface bearing in mind a natural continuous dependence on parameters.

\ssk

\begin{figure}
\centering
\includegraphics[scale=0.75]{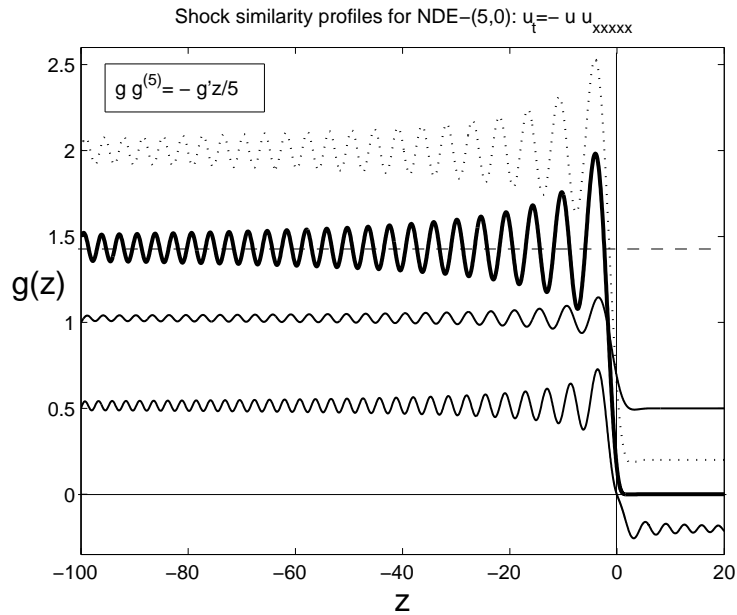}
 \vskip -.3cm
\caption{\small Various shock similarity profiles $g(z)$ as
solutions of the problem (\ref{1E}), (\ref{2.4}).}
\label{F3NNNs}
\end{figure}

 \noi \underline{\em NDE--$(5,0)$}.
Consider next the first ODE (\ref{1E}) for the fully non-divergent
NDE--$(5,0)$ (\ref{N50}). We can again  describe  formation of
shocks (\ref{2.10}); see
 Figure \ref{F3NNNs}.
  The   boldface profile with $C_-=1.4$
and $C_+=0$  has  finite right-hand interface at $z=z_0 \approx
5$,
 with a different expansion
 \beq
 \label{C001s}
   \mbox{$
   g(z)= \frac {6z_0}{5}\, (z_0-z)^4 |\ln(z_0-z)|(1+o(1)) \to 0^+\,\,\,
 \mbox{as}\,\,\, z \to z_0^-.
  $}
   \eeq

\subsection{Shock formation for a uniformly dispersive NDE: an example}

\begin{figure}
\centering
\includegraphics[scale=0.75]{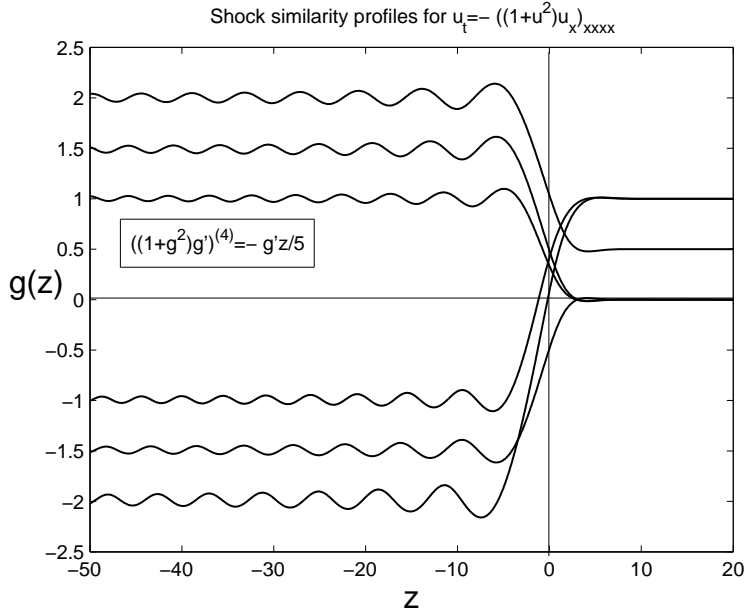}
 \vskip -.3cm
\caption{\small Various shock similarity profiles $g(z)$
satisfying the ODE (\ref{nn2}).}
\label{FMN1}
\end{figure}

Here, as a key example to be continued, we show shocks for {\em
uniform} (non-degenerate) NDEs, such as the fully divergent one,
 \beq
 \label{nn1}
 u_t= -((1+u^2)u_x)_{xxxx},
  \eeq
  where the dispersion coefficient $-(1+u^2)$ of the principal operator is an {\em
  even}
  function. Recall that, for all the previous ones
  (\ref{N50})--(\ref{N14}), the dispersion coefficient is equal to $-u$ and is an odd function of $u$.
 Equation (\ref{nn1}) is non-degenerate and represents a ``uniformly dispersive"
NDE. The ODE for self-similar solutions (\ref{2.1}) then takes the
form
 \beq
 \label{nn2}
 \mbox{$
 ((1+g^2)g')^{(4)}= - \frac 15 \, g'z \, .
  $}
   \eeq
 The mathematics of such equations is similar to that in Section
 \ref{S3.4}. In Figure \ref{FMN1}, we present a few  shock
 similarity profiles for (\ref{nn2}). Note that both shocks
 $S_\pm(x)$ are admissible, since for the ODE (\ref{nn2}) (and for the NDE
 (\ref{nn1})), we have, instead of symmetry (\ref{symm88}),
  \beq
  \label{nn3}
  -g(z) \quad \mbox{is also a solution}.
   \eeq

\section{{\bf (II) Riemann Problem $S_+$}: similarity rarefaction waves}
 \label{Sect3}

Using the reflection symmetry of all the NDEs
(\ref{N50})--(\ref{N14}),
 \beq
 \label{symm1}
 \left\{
 \begin{matrix}
u \mapsto -u,\ssk \\
 t \mapsto -t, \,\,
  \end{matrix}
  \right.
  \eeq
  we conclude that these
 admit {\em global}  similarity solutions defined for all $t>0$,
 \beq
  \label{2.14JJ}
  u_+(x,t) = g(z), \,\,\, \mbox{with} \,\, z=x/t^{\frac 15}.
   \eeq
   Then
 $g(z)$ solves the ODEs (\ref{1E})--(\ref{5E}) with the opposite
 terms
  \beq
  \label{op1}
   \mbox{$
 ... = \frac 15 \, g'z
 $}
  \eeq
  on the right-hand side.
 The conditions (\ref{2.2}) also take the opposite form
  \beq
  \label{2.2op}
  f(\pm \infty)=\pm 1.
   \eeq
 Thus, these profiles are obtained from the blow-up ones in (\ref{2.1})
by reflection, i.e.,
 \beq
 \label{2.15}
\mbox{if  $g(z)$ is a shock profile in (\ref{2.1}), then $g(-z)$
is a rarefaction one in (\ref{2.14JJ}).}
 \eeq

  These are sufficiently regular similarity solutions of
  NDEs that
have  the necessary initial data: by Proposition \ref{Pr.Co}(ii),
in $L^1_{\rm loc}$,
 \beq
 \label{2.3R}
  u_+(x,t) \to S_+(x) \quad \mbox{as}
  \quad t \to 0^+.
   \eeq
Other profiles $g(-z)$ from shock wave similarity patterns
generate further rarefaction solutions including those with finite
left-hand interfaces.


\section{{\bf (III) Riemann Problem $S_-$}: towards $\d$-entropy test}
 \label{SS1}

\subsection{Uniform NDEs}

In this section, for definiteness, we consider the fully
non-divergent NDE (\ref{N50}),
 \beq
 \label{ss1}
 u_t = {\bf A}(u) \equiv  - u u_{xxxxx} \quad \mbox{in} \quad
 \re \times (0,T), \quad u(x,0)=u_0(x) \in C_0^\infty(\re).
  \eeq
In order to concentrate on shocks and to avoid difficulties with
 finite interfaces or transversal zeros at which $u=0$ (these are
weak discontinuities via non-uniformity of the PDE), we deal with
strictly positive solutions satisfying
 \beq
 \label{ss2}
  \mbox{$
  \frac 1C \le u \le C, \quad \mbox{where \,\, $C > 1$ \,\, is a
  constant.}
   $}
   \eeq

   \ssk

\begin{figure}
\centering
\includegraphics[scale=0.75]{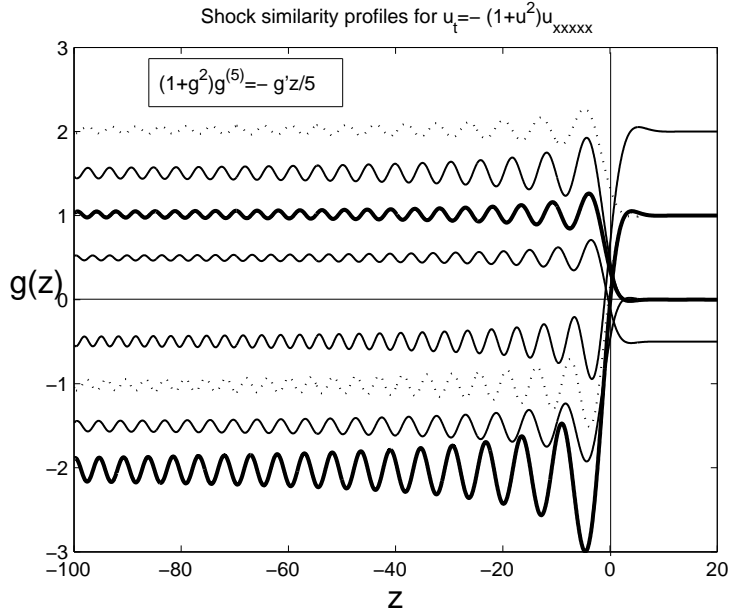}
 \vskip -.3cm
\caption{\small Various shock similarity profiles $g(z)$
satisfying the ODE (\ref{ss4}).}
\label{FNG1}
\end{figure}

   \noi{\bf Remark: uniformly non-degenerate  NDEs.} Alternatively,
  in order to avoid the assumptions like (\ref{ss2}), we can consider
the uniform equations such as (cf. (\ref{nn1}))
 \beq
 \label{ss3}
 u_t = -(1+u^2)u_{xxxxx},
  \eeq
for which no finite interfaces are available. Of course,
(\ref{ss3}) admits analogous blow-up similarity formation of
shocks by (\ref{2.1}). In Figure \ref{FNG1}, we show a few
profiles satisfying
 \beq
 \label{ss4}
  \mbox{$
  (1+g^2)g^{(5)}= - \frac 15 \, g'z, \quad z \in \re.
    $}
    \eeq
Recall that, for (\ref{ss4}), (\ref{nn3}) holds, so both
$S_\pm(x)$ are admissible and entropy (see below).

\subsection{On uniqueness, continuous dependence, and {\em a priori} bounds
 for smooth solutions}

Actually, in our $\d$-entropy construction, we will need just a
local semigroup of smooth solutions that is continuous is
$L^1_{\rm loc}$. The fact that such results are true for
fifth-order (or other odd-order NDEs) is easy to illustrated as
follows.  One can see that, since (\ref{ss1}) is a dispersive
equation, which contains no dissipative terms, the uniqueness
follows as for parabolic equations such as
 $$
  \mbox{$
 u_t =- u u_{xxxx} \quad \mbox{or} \quad u_t = u u_{xxxxxx}
 \quad \big(\mbox{in the class} \,\,\,\big\{\frac 1C \le u \le C\big\}\big).
  $}
  $$

 Thus, we assume that $u(x,t)$ solves (\ref{ss1})
with initial data $u_0(x)\in H^{10}(\re)$, satisfies (\ref{ss2}),
  and  is
  sufficiently smooth, $u \in L^\infty([0,T], H^{10}(\re))$,
  $u_t \in L^\infty([0,T], H^{5}(\re))$, etc.
Assuming that $v(x,t)$ is  the second smooth solution, we subtract
 equations to obtain for the difference $w=u-v$ the PDE
 \beq
 \label{wv1}
 w_t=-u w_{xxxxx} - v_{xxxxx} w.
  \eeq
  We next divide by  $u \ge \frac 1C>0$ and multiply by $w$ in
  $L^2$, so integrating by parts that vanishes the
  dispersive term $w_{xxxxx}$ yields
  \beq
  \label{mm3}
  \mbox{$
  \int \frac{w w_t} u \equiv \frac 12\, \frac{\mathrm d}{{\mathrm d}t}
   \int \frac {w^2}u+ \frac 12 \, \int \frac{u_t}{u^2} \, w^2
=
  - \int \frac{v_{xxxxx}}u \, w^2.
   $}
   \eeq
Therefore, using (\ref{ss2}) and the assumed regularity yields
 \beq
 \label{mm4}
 \mbox{$
 \frac 12\, \frac{\mathrm d}{{\mathrm d}t}
   \int \frac {w^2}u= \int \big(-\frac 12 \, \frac{u_t}{u^2} -\frac{v_{xxxxx}}u\big)  \, w^2
 \le C_1 \int \frac {w^2}u,
  $}
  \eeq
  where the derivatives $u_t(\cdot,t)$ and $v_{xxxxx}(\cdot,t)$ are from
  $L^\infty([0,T])$.
  By Gronwall's inequality, (\ref{mm4}) yields $w(t) \equiv 0$.
Obviously, these estimates can be translated to the continuous
dependence result in $L^2$ and hence in $L^1_{\rm loc}$.

 Other
{\em a priori} bounds on solutions can be also derived along the
lines of computations in \cite[\S\S~2,\,3]{Cr92} that lead to
rather technical manipulations. The principal fact is the same as
seen from (\ref{mm4}): differentiating $\a$ times in $x$ equation
(\ref{ss1}) and setting $v=D_x^\a u$ yields the equations with the
same  principal part as in (\ref{wv1}):
 \beq
 \label{wv2}
 v_t = - u v_{xxxxx}+... \, .
  \eeq
Multiplying this by $\zeta \frac vu$, with $\zeta$ being a cut-off
function, and using various interpolation inequalities makes it
possible to derive necessary {\em a priori} bounds and hence to
observe the corresponding smoothing phenomenon for exponentially
decaying initial data.

\subsection{On local semigroup of smooth solutions of uniform NDEs and linear operator theory}

We recall that local $C^\infty$-smoothing phenomena are known for
third-order linear and fully nonlinear dispersive PDEs; see
  \cite{Cr90, Cr92, Hos99, Lev01} and earlier references therein.
 We claim that, having obtained {\em a priori} bounds, a smooth
 local solution can be constructed by the iteration techniques as
 in \cite[\S~3]{Cr92} by using a standard scheme of iteration of the
  equivalent integral equation for spatial derivatives.
We present further comments concerning other approaches to local
existence, where we return to integral equations.



\ssk

We then need a detailed spectral theory of fifth-order  operators
such as
 \beq
 \label{bb1S}
 \textstyle{
 {\bf P}_5= a(x) D_x^5 + b(x) D_x^4 +... \, ,
 \quad x \in
 (-L,L) \quad \big(a(x) \ge \frac 1C>0\big),
  }
  \eeq
  with bounded coefficients. This theory can be found in, e.g.,  Naimark's
   book \cite[Ch.~2]{NaiPartI}. For {\em regular
  boundary conditions} (e.g., for periodic ones that are regular for any order, which suits us
  well), operators  (\ref{bb1S}) admit a
  discrete spectrum $\{\l_k\}$, where the eigenvalues $\l_k$ are
  all simple for all large $k$.

  It is crucial for further use of eigenfunction expansion techniques
   that  the complete in $L^2$ subset of
  eigenfunctions $\{\psi_k\}$  creates a {\em Riesz  basis}, i.e.,
  for any $f \in L^2$,
 \beq
 \label{zz1}
  \mbox{$
  \sum |\langle f,  \psi_k \rangle|^2 < \infty,
  \quad \mbox{where}
  \quad \langle f,  \psi_k \rangle = \int   f \, \overline
  \psi_k,
   $}
   \eeq
    and, for any  $\{c_k\} \in l^2$ $\big($i.e., $\sum|c_k|^2 <
    \infty\big)$, there exists a function $f \in L^2$ such that
     \beq
     \label{zz2}
      \langle f, \psi_k \rangle = c_k.
       \eeq
Then there exists a unique set of ``adjoint" generalized
eigenfunctions $\{\psi_k^*\}$  (attributed to the``adjoint"
operator ${\bf P}_5^*$) being also a Riesz basis that is
bi-orthonormal to $\{\psi_k\}$:
 \beq
 \label{zz3}
  \langle \psi_k, \psi_l^* \rangle = \d_{kl} \quad\mbox{(Kronecker's delta)}.
   \eeq
 Hence, for any $f \in L^2$, in the sense of the mean convergence,
  \beq
  \label{zz4}
   \mbox{$
  f= \sum c_k \psi_k, \quad \mbox{with} \quad c_k= \langle f,
   \psi_k^*\rangle.
   $}
   \eeq
See further details in \cite[\S~5]{NaiPartI}.


The eigenvalues of (\ref{bb1S})
 have the asymptotics
 \beq
 \label{bb2}
 \l_k \sim (\pm 2\pi k {\rm i})^5 \quad \mbox{for all \,\, $k \gg 1$}.
 \eeq
In particular, it is known that ${\bf P}_5$ has compact resolvent,
which  makes it possible to use it in the integral representation
 of the NDEs; cf. \cite[\S~3]{Cr92}, where integral equations are
 used to construct a unique smooth solution of third-order NDEs.

On the other hand, this means that ${\bf P}_5 - aI$ for any $a \gg
1$ is not a sectorial operator, which makes suspicious using
advanced theory of analytic semigroups \cite{PG, Esch06, Lun}, as
is natural for even-order parabolic flows; see further discussion
below. Analytic smoothing effects for higher-order dispersive
equations were studied in \cite{Tak06}. Concerning unique
continuation and continuous dependence properties for dispersive
equations, see  \cite{Daw07} and references therein, and also
\cite{Tao00} for various estimates.


\subsection{Hermitian spectral theory and analytic semigroups}

Let us continue to discuss related spectral issues for odd-order
operators. For the linear dispersion equation with constant
coefficients (\ref{Lin.717}),
  the Cauchy problem with integrable data $u_0(x)$ admits the
  unique solution
 \beq
  \label{bbb1}
 u(x,t)= b(x-\cdot,t)*u_0(\cdot),
  \eeq
 where $b(x,t)$ is the fundamental solution (\ref{bbb.5}).
Analyticity of solutions in  $t$ (and $x$) can be associated with
the rescaled operator
 \beq
 \label{BB1}
  \mbox{$
 {\bf B}_5= - D_z^5 + \frac 15 \, z D_z + \frac 15\, I
 \quad \mbox{in}
 \quad L^2_\rho(\re),
 $}
  \quad  \mbox{where} \,\,\,\rho(z)=
 \left\{
  \begin{matrix}
 {\mathrm e}^{a|z|^{5/4}},
 \,\,\,z<0, \\ {\mathrm e}^{-a z^{5/4}}, \,\, z>0,
  \end{matrix}
   \right.
 \eeq
  and $a>0$ is a sufficiently small constant.
 Here, ${\bf B}_5$ in (\ref{BB1}) is the operator in (\ref{FODf.1}) that generates the rescaled
 kernel $F$ of the fundamental solution in (\ref{bbb.5}).

  Next,
using in (\ref{Lin.717}) the same rescaling as in (\ref{bbb.5}),
we set
 \beq
 \label{sz1}
 u(x,t)= t^{-\frac 15} v(y,\t), \quad y= x/t^{\frac 15}, \quad \t
 = \ln t,
  \eeq
  to get the rescaled PDE with the operator (\ref{BB1}),
   \beq
   \label{sz2}
   v_\t= {\bf B}_5 v.
    \eeq
  Next,  on Taylor  expansion of the kernel in (\ref{bbb1}) yields
     \beq
     \label{sz3}
      \mbox{$
  v(y,\t)= \int F(y-z{\mathrm e}^{-\t/5}) \, u_0(z)\, {\mathrm d}z
= \sum\limits_{(k)} \frac{(-1)^k}{\sqrt {k!}} \, F^{(k)}(y)
{\mathrm e}^{-\frac k 5\, \t} \frac 1{\sqrt{k!}}\, \int z^k u_0(z)
\, {\mathrm d}z,
 $}
 \eeq
 where the series converges uniformly on compact subsets,
 defining an analytic solution, and also in the mean in $L^2_\rho$.
According to the eigenfunctions expansion (\ref{sz3}) of the
semigroup, there is a proper definition of the operator
 (\ref{BB1}) with a real spectrum and  eigenfunctions (see details in \cite[\S~9]{2mSturm},
  \cite[\S~8.2]{GalPet2m})
  $$
   \mbox{$
   \s({\bf B}_5)=\big\{- \frac
 k 5, \, k=0,1,2,...\big\} \quad \mbox{and}
 \quad \psi_k(y)= \frac{(-1)^k}{\sqrt {k!}} \, F^{(k)}(y), \,\,\, k \ge 0.
  $}
  $$
  The  basis of the ``adjoint" operator (cf. (\ref{2.5})),
  in a space with an indefinite metric,
   $$
    \mbox{$
    \BB_5^*=- D^5_y- \frac 15\, y D_y
    \quad \mbox{in}
 \quad L^2_{\rho^*}(\re), \quad \rho^*(z)={\mathrm e}^{-a|z|^{5/4}}
 \,\, \mbox{in} \,\, \re,
     $}
     $$
 has the same point spectrum and eigenfunctions $\{\psi_k^*\}$,
 which are generalized Hermite polynomials.
   Cf. a full  ``parabolic" version of
 such a Hermitian spectral theory in \cite{Eg4, 2mSturm}.
   This implies  that ${\bf B}_5-a I$ is sectorial for
 $a \ge 0$ ($\l_0=0$ is simple), and this justifies the fact that
 (\ref{bbb1}) is an analytic (in $t$) flow. Let us mention again
 that
analytic smoothing effects are well known for higher-order
dispersive equations with operators of principal type,
\cite{Tak06}.

  Actually, this also suggests
  to treat (\ref{ss1}), (\ref{ss2}) by a classic approach as in Da
 Prato--Grisvard \cite{PG} by linearizing about a sufficiently
 smooth $u_0=u(t_0)$, $t_0 \ge 0$, by setting $u(t)=u_0+v(t)$ giving the
 linearized equation
  \beq
  \label{leq1}
  v_t= {\bf A}'(u_0)v + {\bf A}(u_0)+ g(v), \quad t > t_0; \quad
  v(t_0)=0,
   \eeq
   where $g(v)$ is a quadratic perturbation. Using the good  semigroup
    ${\mathrm e}^{{\bf A}'(u_0)t}$, this makes it possible to
    study local regularity properties of the corresponding
    integral equation
     \beq
     \label{leq2}
     \mbox{$
     v(t) =\int\limits_{t_0}^t{\mathrm e}^{{\bf A}'(u_0)(t-s)}({\bf A}(u_0)+
     g(v(s)))\, {\mathrm d}s.
      $}
       \eeq
 Note that this smoothing approach
demands a fast exponential decay of solutions $v(x,t)$ as $x \to
\infty$, since one needs that $v(\cdot,t) \in L^2_\rho$; cf.
\cite{Lev01}, where $C^\infty$-smoothing for third-order NDEs was
also established under the exponential decay.
Equation (\ref{leq2}) can be used to guarantee local existence of
smooth solutions of a wide class of odd-order NDEs.






\ssk

 Thus, we state the following conclusion to be used later on:
 \beq
 \label{cc11}
  \begin{matrix}
 \mbox{any sufficiently smooth solution $u(x,t)$ of (\ref{ss1}),
 (\ref{ss2}) at $t=t_0$,}\ssk\\
 \mbox{can be uniquely extended to some interval
 $t \in(t_0,t_0+\nu)$,  $\nu>0$.}
  \end{matrix}
  \eeq

\subsection{Smooth deformations and $\delta$-entropy test for  solutions with shocks}
 \label{SD1}

The situation dramatically changes if we want to treat solutions
with shocks. Namely, it is known that even for the NDE--3
(\ref{1}), the similarity formation mechanism of shocks
immediately shows  nonunique extensions of solutions after a
typical ``gradient" catastrophe \cite{Gal3NDENew}. Therefore, we
do not have  a chance to get, in such  an easy (or any) manner, a
uniqueness/entropy result for more complicated NDEs such as
(\ref{N14}) by using the $\delta$-deformation (evolutionary
smoothing) approach. However, we will continue using these  ideas,
 turned out to be fruitful, in order to develop a much weaker {\em
``$\d$-entropy test"} for distinguishing some simple shock and
rarefaction waves.

 Thus, given a
small $\d>0$ and a sufficiently small bounded continuous (and,
possibly, compactly supported)  solution $u(x,t)$ of the Cauchy
problem (\ref{ss1}),
 satisfying (\ref{ss2}), we construct its  smooth {\em $\d$-deformation}, aiming get smoothing in
 a small neighbourhood of bounded shocks,
   as follows. Note that we deal here with simple shock
   configurations (mainly, with 1-shock structures), and do not
   aim to cover more general shock geometry, which can be very
  complicated; especially since we do not know all types of simple
  single-point moving
  shocks.

\ssk

(i) We perform a smooth $\d$-deformation of initial data $u_0(x)$
 by introducing a suitable $C^1$ function $u_{0\d}(x)$ such that
 \beq
 \label{p2}
 \mbox{$
\int |u_0-u_{0\d}| < \d.
 $}
 \eeq
   If $u_0$ is already sufficiently smooth, this step {\em must be abandoned} (now and always later on).
 By $u_{1\d}(x,t)$, we denote the unique local smooth solution of the Cauchy
 problem with data $u_{0\d}$, so that, by (\ref{cc11}),  the continuous
 function
 $u_{1\d}(x,t)$ is defined on the maximal interval $t \in
 [t_0,t_1(\d))$, where we denote $t_0=0$ and  $t_1(\d)=\D_{1\d}$.
 At this step, we are able to eliminate non-evolution (evolutionary unstable)
 initially posed shocks, which then create corresponding smooth rarefaction waves.

 \ssk

 (ii) At $t= \D_{1\d}$, a shock-type discontinuity
  (or possibly infinitely many shocks) is supposed to
 occur, since otherwise we extend the continuous solution by (\ref{cc11}), so we
 perform another suitable $\d$-deformation of the ``data"
 $u_{1 \d}(x,\D_{1\d})$ to get a unique continuous solution $u_{2\d}(x,t)$
 on the maximal interval $t \in [t_1(\d),t_2(\d))$, with
 $t_2(\d)=\D_{1\d}+\D_{2\d}$, etc. Here and in what follows, we always mean a
 ``$\d$-smoothing" performed
 in a small neighbourhood of occurring singularities {\sc only} as
 discontinuous
 shocks.


$\dots$

\ssk

 We continue  in this manner with  suitable
choices of each $\d$-deformations of ``data" at the moments
$t=t_j(\d)$, when $u_{j\d}(x,t)$ has a shock, there exists a
$t_{k}(\d)
> 1$ for some finite $k=k(\d)$, where $k(\d) \to +\infty$ as $\d
\to 0$. It is easy to see that, for bounded solutions, $k(\d)$ is
always finite. A contradiction is obtained by assuming that
$t_j(\d) \to \bar t<1$ as $j \to \infty$ for arbitrarily small
$\d>0$ meaning a kind of ``complete blow-up" that
 was excluded by assumption of smallness of the data.


\ssk

This gives  a {\em global $\d$-deformation} in $\re \times [0,1]$
of the solution $u(x,t)$, which is the discontinuous orbit denoted
by
 \beq
 \label{p3}
 \mbox{$
u^\d(x,t)= \{u_{j\d}(x,t) \,\,\, \mbox{for} \,\,\, t \in
[t_{j-1}(\d),t_j(\d)), \quad j=1,2,...,k(\d)\}.
 $}
 \eeq
One can see that this $\d$-deformation construction aims at
checking a kind of {\em evolution stability} of possible shock
wave singularities and therefore, to exclude those that are not
entropy and evolutionarily generate smooth
 rarefaction waves.

  Finally,  by an
arbitrary {\em smooth $\d$-deformation}, we will  mean the
function (\ref{p3}) constructed by any sufficiently refined finite
partition $\{t_j(\d)\}$ of $[0,1]$,  without reaching a shock of
$S_-$-type at some or all intermediate points $t=t_{j}^-(\d)$.

\ssk

We next say that, given a solution $u(x,t)$, it is {\em stable
relative smooth deformations}, or simply $\d$-stable ({\em
$\d$eformation-stable}), if for any $\e>0$, there exists
$\d=\d(\e)>0$ such that, for any finite $\d$-deformation of $u$
given by (\ref{p3}),
 \beq
 \label{p4}
 \mbox{$
\iint |u-u^\d| < \e.
 $}
 \eeq
 Recall that (\ref{p3}) is an $\d$-orbit, and, in general, is not
 and cannot be aimed to
 represent a fixed solution in the limit $\d \to 0$; see below.


\subsection{On $\d$-entropy solutions}

  Having checked that the local smooth solvability problem above is well-posed,
we now present the corresponding definition that will be
 applied to  particular weak solutions. Recall that the metric
 of convergence, $L^1_{\rm loc}$ under present consideration, for (\ref{1}) was justified
by a similarity analysis presented in Proposition \ref{Pr.Co}. For
other types of shocks and/or  NDEs, the metric  may be different.


 Thus, under the given hypotheses, a function $u(x,t)$ is called
a $\d$-entropy solution of
 the Cauchy problem $(\ref{ss1})$,
  if there exists a sequence of its smooth $\d$-deformations
 $\{u^{\d_k}, \, k=1,2,...\}$, where $\d_k \to 0$, which converges in $L^1_{\rm loc}$ to
 $u$ as $k \to \infty$.


 This is slightly weaker (but equivalent) to the condition of
$\d$-stability.

\ssk

\noi{\bf Remark: $\d$-entropy solution is unique for 1D
conservation law.} Consider, as a typical example, (\ref{3}) for
general measurable $L^1$-data. The classical Oleinik--Kruzhkov's
entropy theory for (\ref{3}) defines the unique semigroup of
contractions in $L^1$ (see \cite{Sm}), i.e., for an arbitrary pair
of entropy solutions $u(\cdot,t)$ and $v(\cdot,t)$, in the sense
of distributions,
 \beq
  \label{SM1}
 \mbox{$
  \frac {{\mathrm d}}{{\mathrm d}t}\, \|u(t)-v(t)\|_{L^1} \le 0
   \quad \mbox{for a.a. $t \ge 0$}.
   $}
   \eeq

Consider now the above $\d$-deformation construction of an orbit
$\{u_\d\}$ in the case, when the entropy solution $u(x,t)$ is
continuous a.e. for all $t \ge 0$, i.e., shocks  have zero
measure. It means that $u_\d(x,t)$ for $t  \ge 0$ is smooth and
essentially differs from $u(x,t)$ on a set of arbitrarily small
measure $\sim \d \to 0$. Therefore, under these (possibly,
non-constructive) assumptions, (\ref{SM1}) implies that any smooth
$\d$-deformations in $L^1$ inevitably lead  to the {\em unique}
entropy solution of (\ref{3}) as $\d \to 0$. In other words,
 \beq
 \label{SM2}
  \fbox{$
 \mbox{for Euler's equation (\ref{3}), classic entropy solutions $=$ $\d$-entropy ones.}
  $}
  \eeq
Of course, this is just the trivial consequence of the
$L^1$-contractivity (\ref{SM1}), which, in its turn, is induced by
the Maximum Principle. It is also worth mentioning that, somehow,
(\ref{SM1}) reflects the fact that the conservation laws such as
(\ref{3}) admit the direct algebraic solution via characteristics.
Indeed, the characteristic method guarantees the unique
solvability in the regularity domain, while the ``shocks cut off"
can be performed at necessary points by the corresponding
Rankine--Hugoniot relations. Thus, the entropy conditions just
describe the correct evolution from initially posed singularities
(evolutionary, such ``rarefaction waves" cannot appear by
characteristics).


Therefore, the absence of the Maximum Principle and absence of any
characteristic-based approaches for higher-order NDEs  recall that
a result such as (\ref{SM2}) cannot be expected  in principle
here. The situation is even more terrible: we will show that any
uniqueness/entropy results for such NDEs fail always and anyway.




\subsection{$\d$-entropy test and  nonexistent uniqueness}



Since, for obvious reasons, the $\d$-deformation construction gets
rid of non-evolutionary shocks (leading to non-singular
rarefaction waves), a first consequence of the construction is
that it defines the {\em $\d$-entropy test} for solutions, which
allows one, at least, to distinguish the true simple isolated
shocks from smooth rarefaction waves.


In Section \ref{SNonU}, we show that it is completely unrealistic
to expect from this construction something essentially stronger
in the direction of uniqueness and/or entropy-like selection of
proper solutions.
Though these expectations   correspond well to previous classical
PDE entropy-like theories, these are  excessive for higher-order
models, where such a universal property is not achievable at all
any more. Even proving convergence for a fixed special
$\d$-deformation is not easy at all.
 Thus, for particular cases, we will use the above notions with
 convergence along a subsequence of $\d$'s
 to classify and distinguish  shocks and rarefaction waves of simple geometric
  configurations:

\subsection{First easy  conclusions of $\d$-entropy test}

   As a first
application, we have:

\begin{proposition}
\label{Pr.E} Shocks of the type $S_-(x)$  are $\d$-entropy for
$(\ref{ss1})$.
 \end{proposition}

 The result follows from the properties of similarity solutions
 (\ref{2.1}), with $-t \mapsto T_t$, which,
  by varying the blow-up time $T \mapsto T+\d$, can be used as their local smooth
 $\d$-deformations
 at any point $t  \ge 0$.


\begin{proposition}
\label{Pr.NE} Shocks of the type $S_+(x)$  are not $\d$-entropy
for $(\ref{ss1})$.
 \end{proposition}

Indeed, taking initial data $S_+(x)$ and constructing its smooth
$\d$-deformation via the self-similar solution (\ref{2.14JJ}) with
shifting $t \mapsto t+\d$, we obtain the global $\d$-deformation
$\{u^\d=u_+(x,t+\d)\}$, which goes away from $S_+$.


Thus, the idea of smooth $\d$-deformations allows   us to
distinguish basic $\d$-entropy and non-entropy shocks {\em
without} any use of mathematical manipulations associated with
standard entropy inequalities, which, indeed, are illusive for
higher-order NDEs; cf. \cite{Gal3NDENew}. We believe that
successful  applications of the $\d$-entropy test can be extended
to any configuration with a {\em finite} number of {\em isolated}
shocks. However, it is completely illusive to think that such a
simple procedure could be applied to general solutions, especially
since the uniqueness after singularity formation cannot be
achieved in principle, as we show next.

In other words, the $\d$-entropy test allows us to {\em prohibit}
formation of non $\d$-deformation stable shocks of type $S_+$ and
proposes a smooth rarefaction wave instead. However, this approach
cannot detect a unique shock of the opposite geometry $S_-$, since
such a formation is principally nonunique.

 \section{{\bf (IV) Nonuniqueness} after shock formation}
 \label{SNonU}

 Here we mainly follow the ideas from \cite{Gal3NDENew} applied there to the NDE--3 \ef{1}, so we
 will omit some technical data and present more convincing
 analytic and numerical results concerning the nonuniqueness. For the hard 5D dynamical
   systems under consideration, numerics becomes more and more essential and unavoidable for understanding the nature
    of such  nonunique extensions of solutions.
  Without loss of generality, we
 always deal with the  NDE--5 (\ref{N14}) of the fully divergent form.

\subsection{Main strategy towards  nonunique continuation: pessimistic
conclusions}
 \label{S5.1}

We begin with the study of new shock patterns, which are induced
by other (cf. (\ref{2.1})) similarity solutions of (\ref{N14}):
 \beq
 \label{s1a}
  \mbox{$
  u_{-}(x,t)=(-t)^\a f(y), \quad y = \frac x{(-t)^\b}, \quad
  \b= \frac {1+\a}5, \quad \mbox{where $\a \in \big(0, \frac 14\big)$ and}
   $}
   \eeq
   \beq
   \label{s2a}
     \left\{
     \begin{matrix}
    -(ff')^{(4)}- \b f'y + \a f=0 \inB \re_-, \quad
    f(0)=f''(0)=f^{(4)}(0)=0, \ssk\\
   f(y) = C_0
    |y|^{\frac \a \b}(1+o(1)) \asA y \to - \iy, \quad C_0>0.
    \qquad\qquad\qquad\quad\,\,\,
     \end{matrix}
     \right.
     \eeq
    In this section, in order to match the key results in
    \cite{Gal3NDENew}, in
      \ef{2.1} and later on, we change the variables $\{g,z\} \mapsto
     \{f,y\}$. In the next Section \ref{SCK1}, we return to the
     original notation.
The anti-symmetry conditions in (\ref{s2a}) allow us  to extend
the solution to the positive semi-axis  $\{y>0\}$ by the
reflection  $-f(-y)$ to get a global pattern.

Obviously, the solutions (\ref{2.1}), which are suitable for
Riemann problems, correspond to the simple case $\a=0$ in
(\ref{s1a}). It is easy to see that, for positive $\a$,
 the asymptotics in \ef{s2a} ensures
getting first {\em gradient blow-up} at $x=0$ as $t \to 0^-$, as a
weak discontinuity, where the final time profile remains locally
bounded and continuous:
 \beq
 \label{s3a}
  u_{-}(x,0^-)= \left\{
   \begin{matrix}
   C_0 |x|^{\frac \a \b}\,\,\,\, \forA x<0, \ssk\\
- C_0 |x|^{\frac \a \b} \forA x>0,
 \end{matrix}
 \right.
  \eeq
  where $C_0>0$ is an arbitrary constant.
  Note that the standard ``gradient catastrophe", $u_x(0,0^-)=
  - \iy$, then occurs in the range, which we will deal within,
 \beq
 \label{Ran1}
  \mbox{$
 \frac \a \b<1 \quad \mbox{provided that} \quad \a < \frac 14.
 $}
  \eeq

   Thus, the wave braking
  (or ``overturning") begins at $t=0$, and next we show that it is
  performed again in a self-similar manner and is described by
  similarity solutions
\beq
 \label{s1Na}
  \mbox{$
  u_{+}(x,t)=t^\a F(y), \quad y = \frac x{t^\b}, \quad
  \b= \frac {1+\a}5, \quad \mbox{where}
   $}
   \eeq
   \beq
   \label{s2Na}
    \left\{
    \begin{matrix}
    -(F F')^{(4)}+ \b F'y - \a F=0 \inB
    \re_-,\qquad\qquad\qquad\qquad\quad\,\,
  \ssk\\
    F(0)=F_0>0, \,\,\,   F(y) = C_0
    |y|^{\frac \a \b}(1+o(1)) \asA y \to - \iy,
     \end{matrix}
     \right.
     \eeq
     where the constant $C_0>0$ is fixed by blow-up data
     (\ref{s3a}).
The asymptotic behaviour as $y \to -\iy$ in (\ref{s2Na})
guarantees the continuity of the global discontinuous pattern
(with $F(-y) \equiv -F(y)$) at the singularity blow-up instant
$t=0$, so that
 \beq
 \label{s3Na}
 u_-(x,0^-)= u_+(x,0^+) \quad \mbox{in} \quad  \re.
  \eeq

Then any suitable couple $\{f,F\}$ defines a global solution
$u_\pm(x,t)$, which is continuous at $t=0$, and then it is called
an {\em extension pair}. It was shown in \cite{Gal3NDENew} that,
for the typical NDEs--3, the pair is not uniquely determined and
there exist infinitely many shock-type extensions of the solution
after blow-up at $t=0$. We are going to describe a similar
nonuniqueness phenomenon for the NDEs--5 such as \ef{N14}.

It is worth mentioning that, for conservation laws such as
(\ref{3}), such an extension pair $\{f,F\}$ {\em is always
unique}; see similarity analysis in \cite[\S~4]{Gal3NDENew}. Of
course, this is not surprising due to existing Oleinik--Kruzhkov's
classic uniqueness-entropy theory \cite{Ol59,Kru2}. Note again
that any sufficient multiplicity of extension pairs $\{f,F\}$,
obtained via small micro-scale blow-up analysis of the PDEs, would
always lead to a principle nonuniqueness, so this approach could
be referred to as a {\rm ``uniqueness test"}.


A first immediate consequence of our similarity blow-up/extension
analysis is as follows:
 \beq
 \label{s5N}
 \mbox{in the CP, formation of shocks for the NDE (\ref{N14}) can lead to
 nonuniqueness.}
 \eeq
The second conclusion is more subtle and is based on
 the fact that,
 for some initial data at $t=0$ (i.e., created by single point {\em gradient blow-up}
 as $t \to 0^-$), the whole admitted solution set for $t>0$
    does not contain any ``minimal",
``maximal", ``extremal" in any reasonable sense, or any isolated
points, which might
play a role of a unique ``entropy" one chosen by introducing a
hypothetical entropy inequalities, conditions, or otherwise.
If this is true for the whole set of such weak solutions of
\ef{N14} with initial data (\ref{s3a}), then, for the Cauchy
problem,
 \beq
 \label{s6N}
 \mbox{there exists
  no  general ``entropy mechanisms" to choose a unique
 solution.}
  \eeq
Actually,  overall, \ef{s5N} and \ef{s6N} show that the problem of
uniqueness of weak solutions for the NDEs such as \ef{N14} {\em
cannot be solved in principal}.

On the other hand, in a FBP setting by adding an extra suitable
 condition on shock lines, the problem might be well-posed with a
 unique solution, though proofs can be very difficult.
 We refer again to a more detailed discussion of these issues for
 the NDE--3 \ef{1} in \cite{Gal3NDENew}. Though we must admit that, for
 the NDE--5 \ef{N14}, which induces 5D dynamical systems for the
 similarity profiles (and hence 5D phase spaces), those
 nonuniqueness and non-entropy conclusions are
 more difficult and not that clear as for the NDEs--3,
 so some of their aspects do unavoidably  remain questionable and
 even open.

Hence, the {\em nonuniqueness} in the CP is a non-removable issue
of PDE theory for higher-order degenerate nonlinear odd-order
equations (and possibly not only for those).  The nonuniqueness of
solutions of \ef{N14} has some pure  dimensional natural features,
and, more precisely, is associated with the dimensions of ``good"
and ``bad" asymptotic bundles of orbits in the 5D phase space of
the ODE \ef{s2Na}.



\subsection{Infinite shock similarity solutions for $\a<0$}

\begin{figure}
\centering
\includegraphics[scale=0.60]{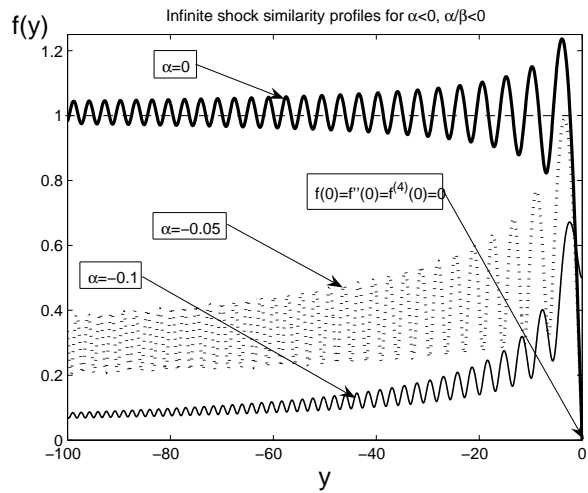}
\caption{\small Infinite shock similarity profiles as solutions of
(\ref{s2a}) for $\a<0$: $\a=-0.05$ and $\a=-0.1$.}
\label{FalN}
\end{figure}

  Let us first note that the blow-up solutions (\ref{s1a})
  represent an effective way to describe other types of
  singularities with {\em infinite shocks}. Namely, assuming that
   \beq
   \label{ra11}
   \mbox{$
   \a<0 \andA \frac \a \b <0,
    $}
    \eeq
  we again obtain the same ``data" \ef{s3a} but now $u_-(0,0^-)= \iy$.
 We do not study in any detail such interesting new singularity
 phenomena and present Figure \ref{FalN} showing that such
 infinite shock similarity profiles do exist.
 For comparison, we indicate the standard $S_-$-type profile for
 $\a=0$, which coincides with that in Figure \ref{F7}(d).

For the NDE--3 such as (\ref{1}), the infinite shock similarity
solutions in the range (\ref{ra11}) were studied  in
\cite[\S~4]{GPnde} in sufficient detail.

\subsection{Gradient blow-up similarity solutions}
 \label{S2a}

 Consider the blow-up ODE problem (\ref{s2a}), which is a
 difficult one, with a {\em 5D phase space}. Note that, by invariant
 scaling (\ref{2.8}), it can be reduced to a 4th-order ODE with
a also even more  complicated nonlinear operator composed from too
many polynomial terms, so we do not rely on that and
 work in the original phase space.
Therefore,  some more delicate issues on, say, uniqueness of
certain orbits, become very difficult or even remain open, though
some more robust properties can be detected rigorously.
 We will also  use
numerical methods for illustrating and even justifying some of our
conclusions. As before, for the fifth-order equations such as
(\ref{s2a}), this and further numerical constructions are
performed by the {\tt MatLab} with the standard {\tt ode45} solver
therein.

 Let us describe the necessary properties of orbits $\{f(y)$\} we are
 interested in.
Firstly, it follows from the conditions in \ef{s2a} that, for $y
\approx 0^-$,
 \beq
 \label{par55}
 \mbox{the set of proper orbits is 2D parameterized by
 $f_1=f'(0)<0$ and $f_3=f'''(0)$.}
  \eeq

Secondly and on the other hand,
 the necessary behaviour at infinity is as follows:
 \beq
 \label{in1a}
  \mbox{$
 f(y)=C_0 |y|^{\frac {5\a}{1+\a}}(1+o(1)) \asA y \to - \iy \quad \big( \frac{5\a}{1+\a}
 =\frac \a \b \big),
 $}
  \eeq
  where $C_0>0$ is an arbitrary constant by
  scaling (\ref{2.8}).
 It is key to derive the whole 4D
 bundle of solutions satisfying (\ref{in1a}). This is done
 by the linearization as $y \to -\iy$:
  \beq
  \label{in12a}
   \begin{split}
   & f(y)=f_0(y)+Y(y) \whereA f_0(y)=C_0(-y)^{\frac \a \b} \ssk\\
   \LongA
    &
     \mbox{$
     -  C_0 ((-y)^{\frac \a \b}Y)^{(5)}+ \b Y'(-y)+ \a Y + \frac 12\,
   (f_0^2(y))^{(5)}+...=0.
   $}
    \end{split}
   \eeq
   By WKBJ-type asymptotic techniques in ODE
   theory,  solutions of (\ref{in12a}) have a standard exponential
   form with the characteristic equation:
 \beq
 \label{in22a}
  \mbox{$
  Y(y) \sim {\mathrm e}^{a(-y)^\g}, \,\,\, \g=1+ \frac 14\,\big(1-
  \frac \a \b\big)>1
 \LongA  C_0(\g a)^4=\b,
  $}
  \eeq
which has three roots with  non-positive real parts, ${\rm Re}\,
a_k \le 0$, where $a_1<0$ is real and conjugate  $a_{2,3}\in {\rm
i}\,\re$. Hence, we conclude that:
 \beq
 \label{in23a}
 \mbox{as $y \to -\iy$, the bundle (\ref{in1a}) is four-dimensional (including $C_0$).}
  \eeq

   The behaviour corresponding the bundle \ef{in23a} gives the desired
   asymptotics.
   Indeed,
   by
  \ef{in1a}, we have the gradient blow-up behaviour at a single
  point: for any fixed $x<0$, as $t \to 0^-$, where $y=x/(-t)^\b
  \to - \iy$,  uniformly on compact subsets,
   \beq
   \label{in2a}
    \mbox{$
    u_-(x,t)= (-t)^\a f(y) = (-t)^\a C_0 \big| \frac
    x{(-t)^\b}\big|^{\frac \a \b}(1+o(1)) \to C_0
    |x|^{\frac{5\a}{1+\a}}.
    $}
    \eeq

\begin{figure}
\centering
\subfigure[$\a= \frac 19$]{
\includegraphics[scale=0.52]{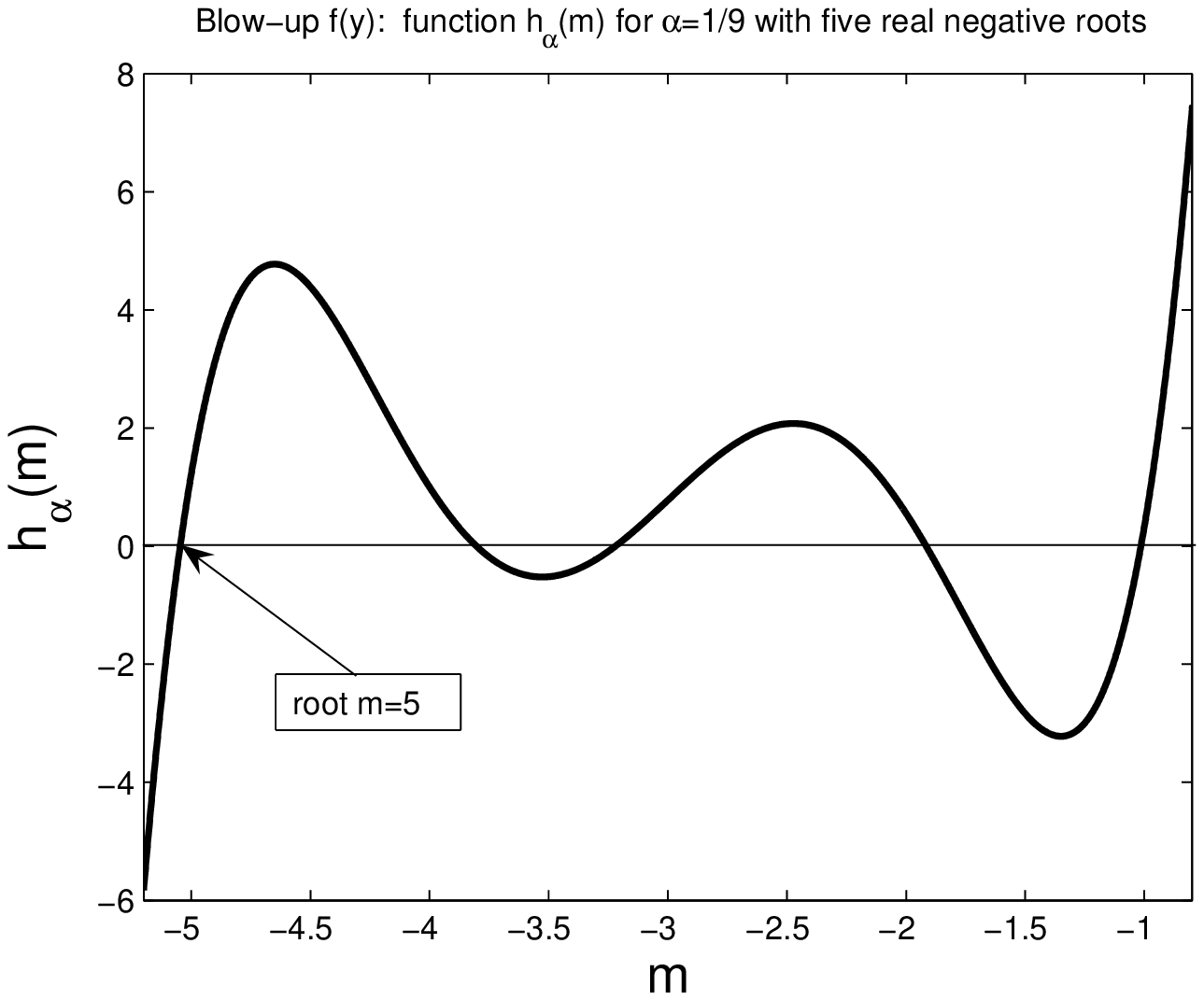}
}
\subfigure[various $\a$]{
\includegraphics[scale=0.52]{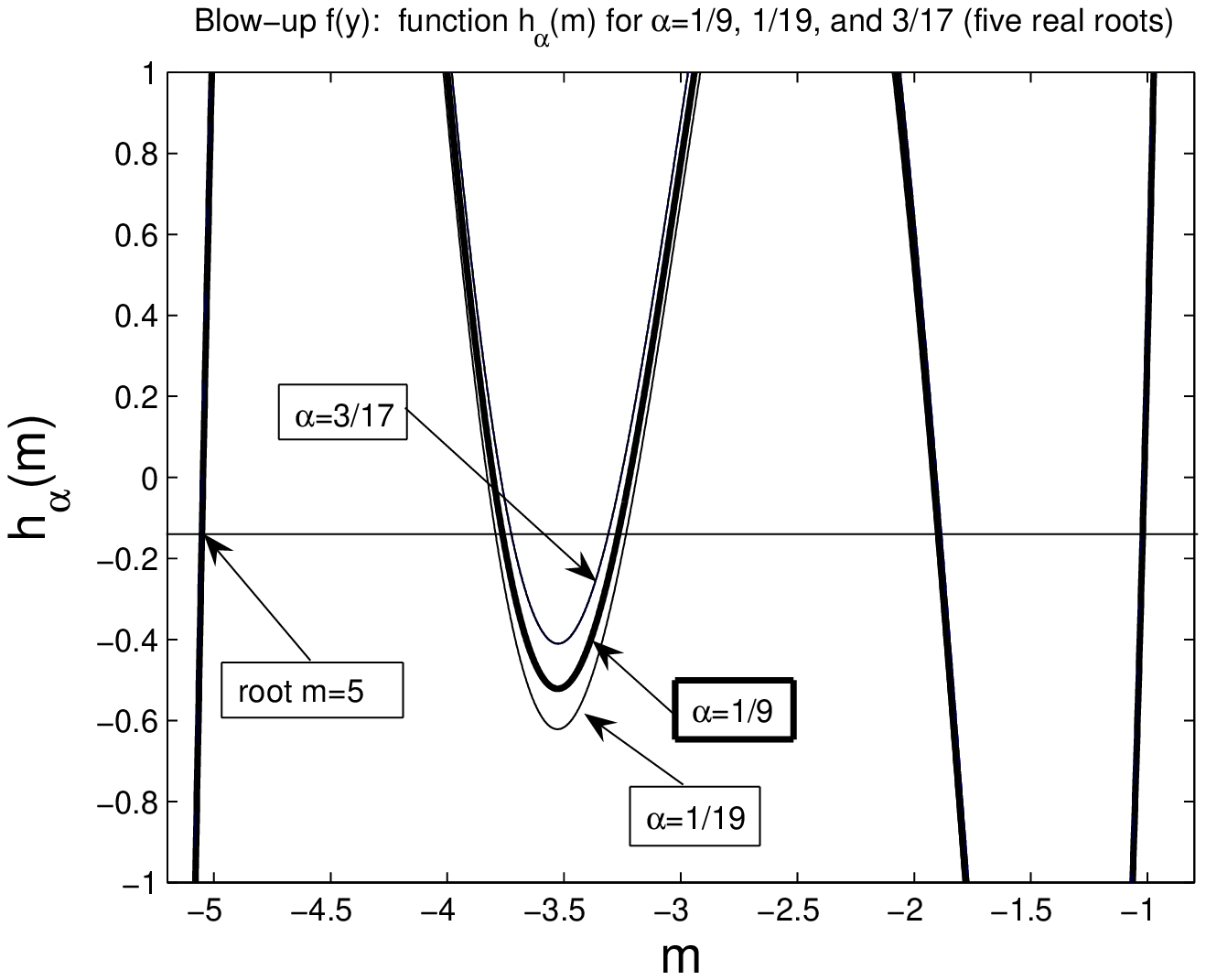}
}
 \vskip -.2cm
\caption{\rm\small The polynomial $h_\a(m)$ in (\ref{eu11}) for
various $\a \in \big(0, \frac 14\big)$: five negative roots.}
 \label{Falph1}
\end{figure}

Let us explain some other crucial properties of the phase space,
now meaning  ``bad bundles" of orbits. First, these are
  the fast growing solutions according
to the explicit solution
 \beq
 \label{2.9a}
  \mbox{$
  f_{\rm *}(y)=-  \frac {y^5}{15120} > 0 \quad \mbox{for} \quad y \ll -1.
 $}
   \eeq
 Analogously to \ef{in12a}, we compute the whole bundle
 about \ef{2.9a}:
\beq
  \label{in12aN}
    f(y)=f_*(y)+Y(y) \LongA
     \mbox{$
       \frac 1{15120}\, (y^5 Y)^{(5)}- \b Y' y+ \a Y +...=0.
   $}
   \eeq
This Euler equation has the following solutions with the
characteristic polynomial:
 \beq
 \label{eu11}
  \mbox{$
  Y(y)=y^m \LongA h_\a(m) \equiv
  \frac{(m+1)(m+2)(m+3)(m+4)(m+5)}{15120}- \b m +\a=0.
   $}
    \eeq
    One root $m=-5$ is obvious that gives the solution \ef{2.9a}.
It turns our that this algebraic equation has precisely five
negative real roots for $\a$ from the range \ef{Ran1}, as Figure
\ref{Falph1} shows. Actually, (b) explains that the graphs are
rather slightly dependent on $\alpha$. Thus:
 \beq
 \label{eu12}
 \mbox{the bundle about (\ref{2.9a}) is five-dimensional.}
  \eeq


Second, there exists a bundle of positive solutions vanishing at
some finite $y \to y_0^+<0$ with the behaviour (this bundle occurs
from both sides, as $y \to y_0^\pm$ to be also used)
 \beq
 \label{bb1}
  \mbox{$
 f_1(y)=A \sqrt{|y-y_0|}\,(1+o(1)), \quad A >0,
 $}
  \eeq
is 4D, which also can be shown by linearization about \ef{bb1}.
Indeed, the linearized operator contains the leading term
 \beq
 \label{bb2a}
  \mbox{$
  -A^2(\sqrt{|y-y_0|}\,\,Y)^{(5)}+...=0 \LongA
  Y(y) \sim |y-y_0|^{\frac 32}, \,\,\, |y-y_0|^{\frac 52}, \,\,\,|y-y_0|^{\frac
  72},
   $}
   \eeq
   which together with the parameter $y_0<0$ yields
 \beq
 \label{bb3}
 \mbox{the bundle about (\ref{bb1}) is four-dimensional.}
  \eeq

 Thus, \ef{par55}, \ef{in23a}, \ef{eu12}, and \ef{bb3} prescribe
key aspects of the 5D phase space we are dealing with. To get a
global orbit $\{f(y),\, y \in \re_-\}$ as a connection of the
proper bundles  \ef{par55} and \ef{in23a}, it is natural to follow
the strategy of ``shooting from below" by avoiding the bundle
\ef{bb1}, \ef{bb3}, i.e.,  using the parameters $f_{1,3}$ in
\ef{par55}, to obtain
 \beq
 \label{y01}
 y_0=-\iy.
  \eeq
  It is not difficult to see that this profile $f(y)$ will belong
  to the bundle \ef{in23a}. The proof of such a 2D shooting strategy
  can be done by standard arguments.
  By scaling \ef{2.8}, we always can reduce the problem to a 1D
  shooting (recall that $f_0=f_2=f_4=0$ already):
  \beq
  \label{y02}
  f_1 \equiv f'(0)=-1 \andA f_3 \equiv f'''(0) \,\,\, \mbox{is a parameter.}
   \eeq
 By the above asymptotic analysis of the 5D phase space, it
 follows that:

 {\bf (I)} for $f_3 \ll -1$ the orbit belongs to the bundle about
 \ef{2.9a}, and

 {\bf (II)} for $f_3 \gg 1$, the orbit vanishes at finite $y_0$ along
 \ef{bb1}.

 Hence, by continuous dependence, we obtain a solution $f(y)$ by
 the min-max principle (plus some usual technical details that can be omitted).
 Before stating the result, for convenience, in Figure \ref{Fa33}, obtained by the {\tt ode45} solver,
  we
  explain how we are going to justify existence of a proper blow-up
 shock profile $f(y)$; cf. Figure \ref{FF1n}.

  Thus, we fix the above speculations as follows:

\begin{figure}
\centering
\includegraphics[scale=0.70]{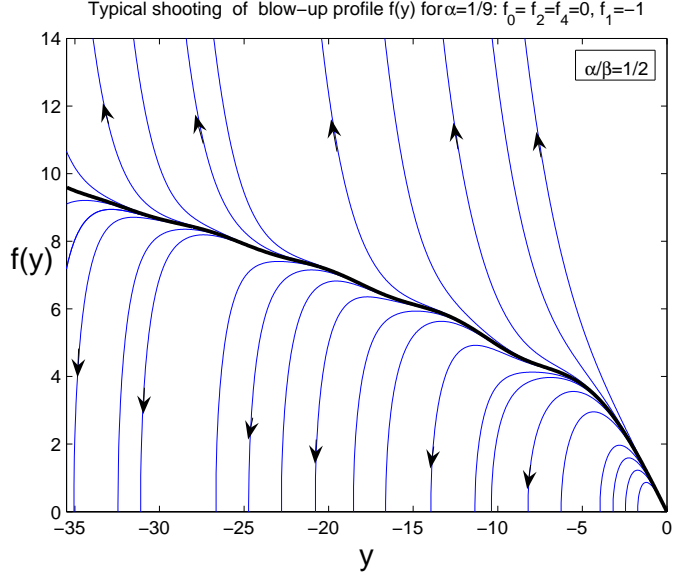}
\caption{\small The shooting strategy of a blow-up similarity
profile $f(y)$ for $\a= \frac 19$, with data
$f(0)=f''(0)=f^{(4)}(0)=0$ and $f'(0)=1$; the shooting parameter
is $f'''(0)=0.0718040128557...$\,.}
\label{Fa33}
\end{figure}



\begin{proposition}
 \label{Pr.1a}
 {\rm (i)} In the range $(\ref{Ran1})$, the problem $(\ref{s2a})$
admits a shock profile $f(y)$.

 \end{proposition}

  We have  the following expectation:
 \noi{\rm (ii)} {\em  $f(y)$ is unique up to scaling $(\ref{2.8})$
 and is positive
for $y<0$.}
This remains an open problem that was confirmed numerically.
  In \cite{Gal3NDENew}, for the NDE--3 \ef{1}, the phase space
is 3D  and a full proof is available.


 In fact, this is a rather typical result for higher-order dynamical
 systems. E.g., we refer
to  a similar and not less complicated study of a 4th-order ODE
\cite{Gaz06}, where existence and {\em uniqueness} of a positive
solution of the radial bi-harmonic equation with source:
 \beq
 \label{rrr1a}
 \D^2_r u= u^p \forA r=|x|>0, \quad u(0)=1, \quad u'(0)=u'''(0)=0, \quad
 u(\iy)=0,
  \eeq
  was proved in the supercritical Sobolev range
   $
   p > p_{\rm Sob}= \frac{N+4}{N-4}$, $N>4$.
   Here, analogously, there exists a single shooting parameter
   being the second derivative at the origin $u_2=u''(0)$; the value $u_0=u(0)=1$ is
    fixed by a scaling symmetry.
 Proving uniqueness of such a solution in \cite{Gaz06} is not easy and
 leads
 to essential technicalities, which the attentive reader can
 consult in case of necessity. Fortunately, we are not interested
 in any uniqueness of such kind.
 Instead of the global behaviour such as \ef{2.9a}, the equation
 \ef{rrr1a} admits the blow-up one governed by the principal
 operator
  $
  u^{(4)}+...= u^p \quad (u \to +\iy).
   $
The solutions vanishing at finite point otherwise can be treated
as in the family {\bf (I)}.

\ssk

\begin{figure}
\centering
\includegraphics[scale=0.70]{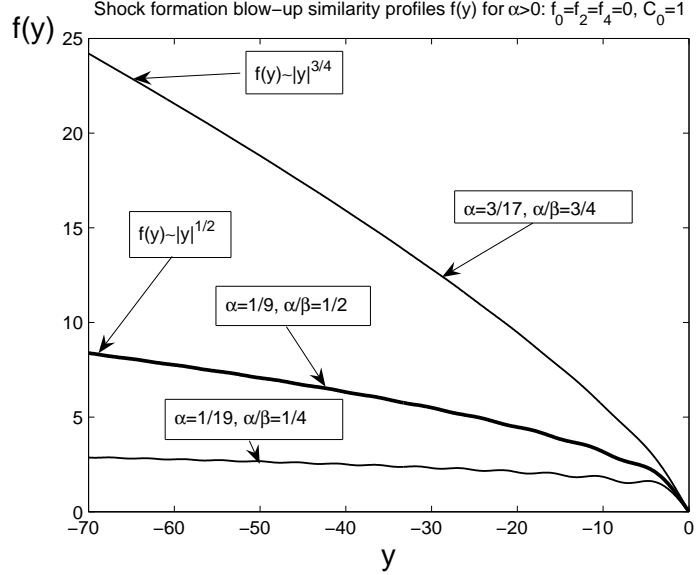}
\caption{\small The odd blow-up similarity profiles $f(y)$ in
$\re_-$ with $\a= \frac 19$,  $\frac 1{19}$, and $\frac 3{17}$.}
\label{Fa2}
\end{figure}

\noi\underline{\sc More numerics by {\tt bvp4c}}. We next use more
advanced and enhanced numerical methods towards existence (and
uniqueness-positivity, see (ii)) of $f(y)$. Figure \ref{Fa2} shows
 blow-up profiles,  with $f(0)=0$, constructed by a
 different method (via the solver {\tt bvp4c}) for convenient values $\a= \frac 19$,
 $\frac 1{19}$, and $\frac 3{17}$. Note the clear oscillatory
 behaviour of such patterns that is induces by complex roots
 of the characteristic equation (\ref{in22a}).

\ssk

\begin{figure}
\centering
\includegraphics[scale=0.70]{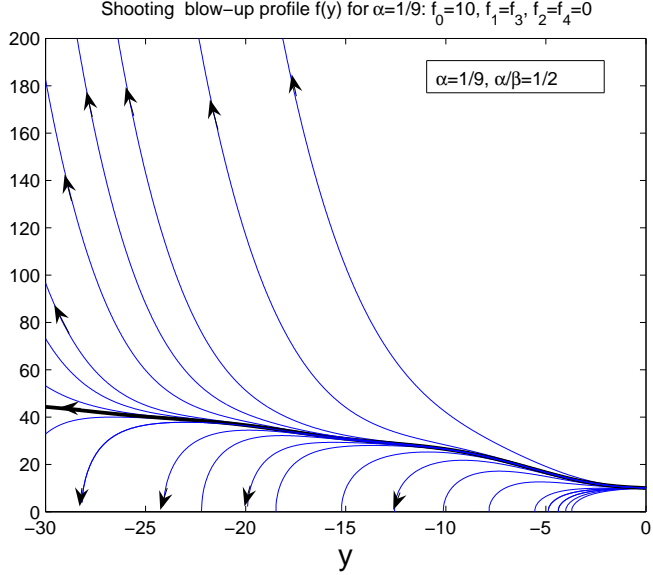}
\caption{\small Shooting the blow-up profile $f(y)$ for $\a= \frac
19$: $f(0)=10$, $f''(0)=f^{(4)}(0)=0$ and the shooting parameter
is $f'(0)=f'''(0)$.}
\label{Fa1}
\end{figure}

\noi\underline{\sc Collapse of shocks: ``backward nonuniqueness"}.
This new phenomenon is presented in Figure \ref{Fa1}, which shows
the shooting from $y=0^-$
 for
 \beq
 \label{al12}
  \mbox{$
 \a= \frac 19 \LongA \frac \a \b = \frac 12.
  $}
  \eeq
 This again illustrates the actual strategy in proving  Proposition
\ref{Pr.1a}. However, though the phase space looks similar, note
that here, as an illustration of another important evolution
  phenomenon, we solve the problem with $f(0) \ne 0$, so that
  there exists a non-zero jump of $u_-(x,t)$ at $x=0$ denoted by
  $[\cdot]$:
  \beq
  \label{pr11}
  f(0)=f_0=10 \LongA  [u_-(0,t)]= 2 f_0(-t)^\a \to 0 \asA t \to
  0^-.
   \eeq
     Therefore, this similarity solution describes {\em collapse of a shock wave} as $t \to 0^-$.

\begin{figure}
\centering
\includegraphics[scale=0.9]{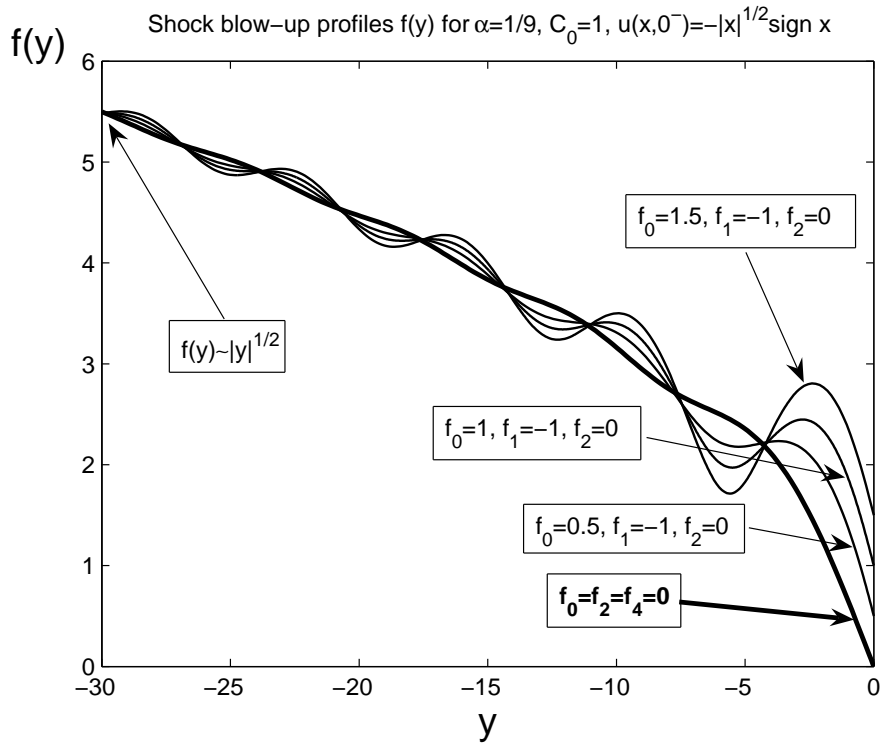}
\caption{\small Blow-up profiles $f(y)$ for $\a= \frac 19$, with
$f'(0)=-1$, $f''(0)=0$; $f(0) \in [0,1.5]$ being a parameter.}
\label{Fa1N}
\end{figure}

\begin{figure}
\centering
\subfigure[$f'(0)=f''(0)=0$]{
\includegraphics[scale=0.52]{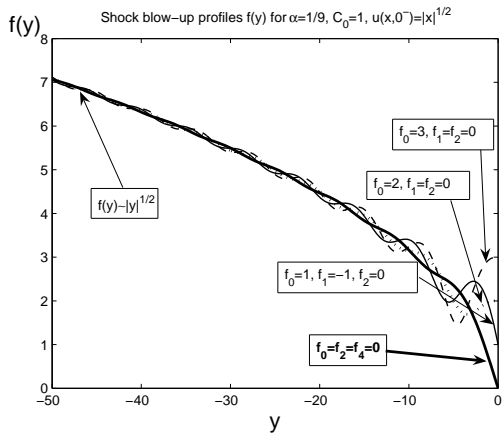}
}
\subfigure[$f''(0)=\frac 12$, $f^{(4)}(0)=0$]{
\includegraphics[scale=0.52]{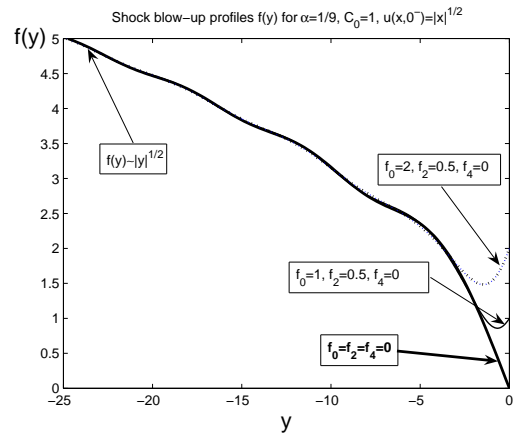}
}
 \vskip -.2cm
\caption{\rm\small Blow-up profiles $f(y)$ for $\a= \frac 19$, for
$f'(0)=f''(0)=0$ (a) and $f''(0)=\frac 12$, $f^{(4)}(0)=0$ (b);
$f(0)$ is a parameter.}
 \label{Fa2N}
\end{figure}

More numerical results of such types are presented in Figures
\ref{Fa1N} and \ref{Fa2N}, where we use other boundary conditions
at $y=0$. Note that, being extended for $y>0$ in the
anti-symmetric way, by $-f(-y)$, this will give a proper shock
wave solution with the nil speed of propagation (see the R--H
condition (\ref{jj2}) below).

In a whole, since all these blow-up profiles satisfy the necessary
behaviour as $y \to -\iy$ as indicated in (\ref{s2a}), these
create as $t \to 0^-$ the same initial data (\ref{s3a}). This
confirms the following phenomenon of ``{\bf backward
nonuniqueness"}: {\em initial data $(\ref{s3a})$ with gradient
blow-up at $x=0$ can be created by an infinite number $($in fact,
by a 2D subset parameterized, say,  by $\{f_0,f_1\})$ of various
self-similar solutions $(\ref{s1a})$.}

Indeed, such a nonuniqueness is directly associated with the fact
that, due to (\ref{in23a}), the proper asymptotic bundle as $y \to
-\iy$ is 3D (for a fixed $C_0>0$, we have to subtract the
dimension via the scaling invariance (\ref{2.8})). Therefore,
roughly speaking, shooting from $y=0^-$ with 5 parameters
$f_0=f(0)$, ...\,, $f_4=f^{(4)}(0)$ allows a 2D ($2=5-3$) subset
of solutions $f(y)$ with shocks at $y=0$. A full justification of
such a conclusion requires a more careful analysis of the phase
space including geometry of two ``bad" bundles, which we do not
perform here concentrating on other more important solutions and
true nonuniqueness phenomena.

\ssk

 \noi\underline{\sc Stationary solutions with a ``weak
shock"}. The ODE in \ef{s2a} and hence the PDE \ef{N14} admit a
number of simple continuous ``stationary" solutions. E.g.,
consider
 \beq
 \label{st1}
 \mbox{$
 \a= \frac 19, \,\, \frac \a \b= \frac 12: \quad
 \hat f(y)= \sqrt{|y|}\,\,{\rm sign} \, y \andA \hat u(x,t) \equiv   \pm \sqrt{|x|}\,\,{\rm sign} \,
 x.
  $}
  \eeq
 Note that these {\em are not} weak solutions of the stationary
 equation
 \beq
 \label{st11N}
  \mbox{$
 \frac 12\, (u^2)_{xxxxx}=0 \quad \mbox{in} \quad {\mathcal D}'.
 $}
 \eeq
 The classic stationary solution of (\ref{st11N}) $\hat u(x,t)= \pm x^2$ is
  smoother at $x=0$.

 We will show that such ``weak stationary shocks" as in (\ref{st1}) also lead to nonuniqueness.

\ssk

\noi{\bf Remark: an exact solution for a critical $\alpha$.} One
can see that the quadratic operator $\BB(f)=(f f')^{(4)}$ in
\ef{s2a} admits the following polynomial invariant subspace:
 $$
 W_6={\rm Span}\{1,y,y^2,y^3,y^4,y^5\} \LongA \BB(W_6) \subseteq
 W_6.
  $$
  Restricting the ODE \ef{s2a} to $W_6$ yields an algebraic
 system, which admits an exact solution for the following value
 of the {\em critical} $\a_{\rm c}$:
  \beq
  \label{cr11}
   \mbox{$
   \a=\a_{\rm c}=\frac{17}{84}=0.202381... \LongA \exists \,\, f(y)= C y- \frac
   {4!}{9!}\, y^5, \,\,\, C \in \re.
   $}
   \eeq
Since $\a_{\rm c}>0$, it does not deliver a ``saw"-type blow-up
profile (having infinite number of positive humps) as it used to
be for the NDE--3 \ef{1} for $\a_{\rm c}=-\frac 1{10}$; see
\cite[\S~4]{GPnde}.

\subsection{Nonuniqueness of similarity extensions beyond blow-up}
 \label{S3aa}


As in \cite{Gal3NDENew} for the NDEs--3,  a discontinuous shock
wave extension of blow-up similarity solutions \ef{s1a}, \ef{s2a}
is assumed to be done by using the global ones \ef{s1Na},
\ef{s2Na}. Actually, this leads to watching a whole 5D family of
solutions parameterized by their Cauchy values at the origin:
 \beq
 \label{N1a}
 F(0)=F_0>0, \,\,\, F'(0)=F_1<0,\,\,\, F''(0)=F_2,\,\,\, F'''(0)=F_3, \,\,\, F^{(4)}(0)=F_4.
  \eeq
 Thus, unlike \ef{par55}, the proper bundle in \ef{N1a} is 5D.
 Note that at $y=-\iy$, the solution must have the form
 \beq
 \label{in1aGG}
  \mbox{$
 F(y)=C_0 |y|^{\frac {5\a}{1+\a}}(1+o(1)) \asA y \to - \iy \quad
 (C_0>0).
 $}
  \eeq

\ssk

As above,
 the 5D phase space for the ODE in \ef{s2Na} has two stable
 ``bad" bundles:

\ssk

 \noi{\bf (I)} Positive solutions with ``singular extinction" in finite $y$, where
 $F(y) \to 0$ as $y \to y_0^+<0$. This is an unavoidable
 singularity following from the degeneracy of the equations with
 the principal term $F F^{(5)}$ leading to the singular potential $\sim \frac
 1F$. As in  \ef{bb2a}, this bundle is 4D, and

\ssk

 \noi{\bf (II)} Negative solutions with the fast growth (cf. (\ref{2.9a})):
  \beq
  \label{f1a}
   \mbox{$
  F_*(y) =  \frac {y^5}{15120}\,(1+o(1)) \to -\iy \asA y \to - \iy.
  $}
  \eeq
   The characteristic polynomial is the same as in \ef{eu11}, so
   that the bundle is 5D; cf. \ef{eu12}.

\ssk

 Both sets of such solutions are {\em open} by the standard
continuous dependence of solutions of ODEs on parameters.
 The whole bundle of solutions satisfying (\ref{in1a})
 is obtained by  linearization as $y \to -\iy$ in \ef{s2Na}:
  \beq
  \label{in12Ga}
   \begin{split}
   & f(y)=F_0(y)+Y(y) \whereA F_0(y)=C_0(-y)^{\frac \a \b} \ssk\\
   \LongA
    &
     \mbox{$
     - C_0 ((-y)^{\frac \a \b}Y)^{(5)}- \b Y'(-y)- \a Y + \frac 12\,
   (F_0^2(y))^{(5)}+...=0.
   $}
    \end{split}
   \eeq
   The WKBJ method now leads to a different characteristic equation:
 \beq
 \label{in22Ga}
  \mbox{$
  Y(y) \sim {\mathrm e}^{a(-y)^\g}, \,\,\, \g=1+ \frac 14\,\big(1-
  \frac \a \b\big)>1
 \LongA  C_0(\g a)^4=-\b,
  $}
  \eeq
 so that there exist just two complex conjugate roots with Re$\,\le 0$,
 and hence, unlike \ef{in23a},
\beq
 \label{in23GaG}
 \mbox{the bundle (\ref{in1a}) of global orbits $\{F(y)\}$ is three-dimensional.}
  \eeq

\begin{figure}
\centering
\subfigure[$F_1=-1, \, F_2=F_3=0$]{
\includegraphics[scale=0.52]{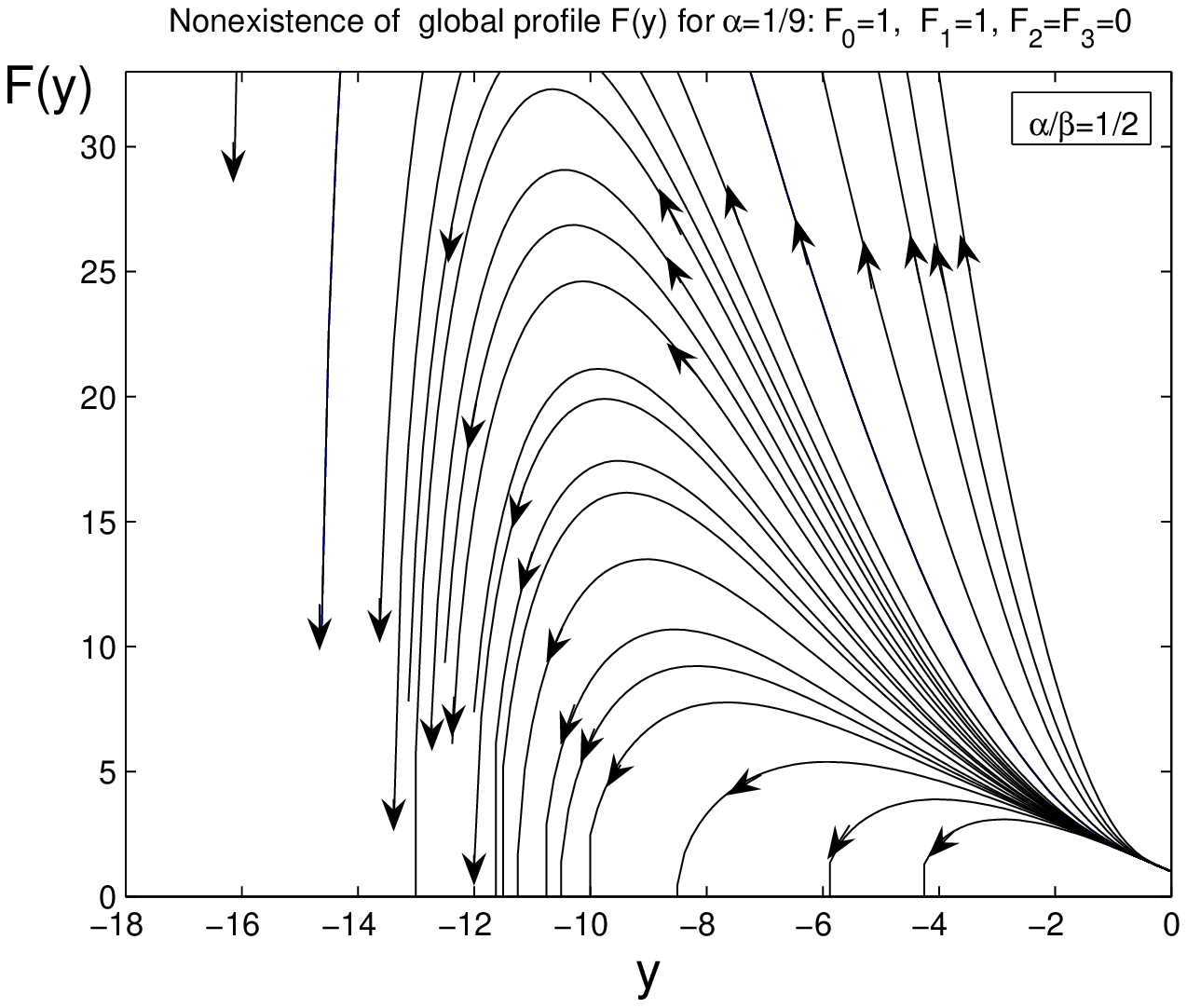}
}
\subfigure[$F_1=F_2=F_3=F_4$]{
\includegraphics[scale=0.52]{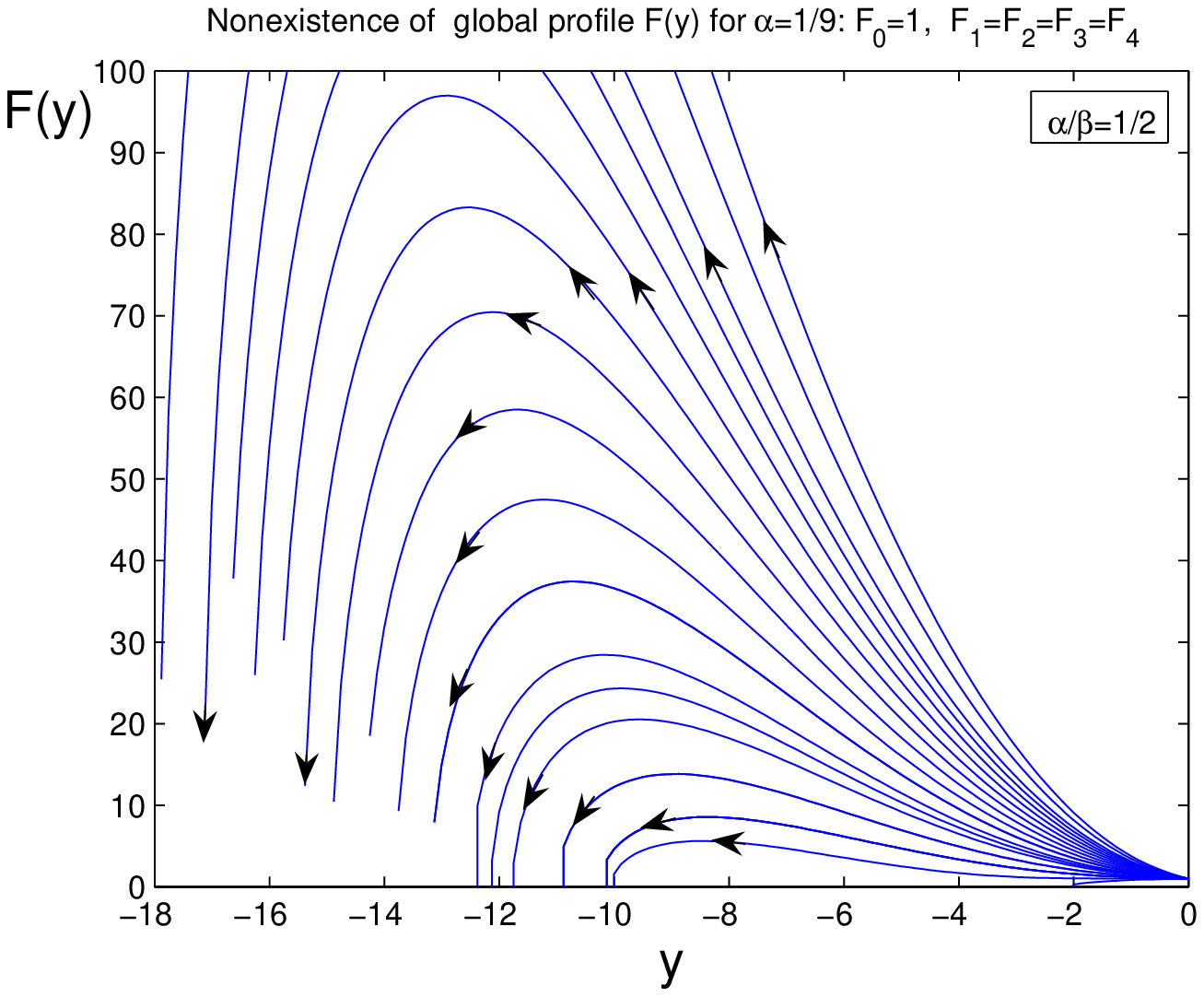}
}
 \vskip -.2cm
\caption{\rm\small Unsuccessful examples of 1D shooting of $F(y)$
of (\ref{s2Na})
 from $y=0^-$.}
 \label{FNon88}
\end{figure}
\begin{figure}
\centering
\includegraphics[scale=0.55]{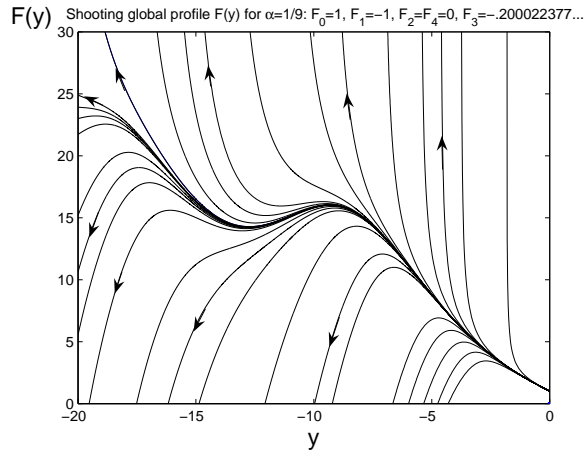}
\caption{\small Unsuccessful  1D shooting of $F(y)$ satisfying
(\ref{s2Na})
 from $y=0^-$, with conditions $F(0)=1$, $F'(0)=-1$, $F''(0)=F^{(4)}(0)=0$,
 and $F'''(0)=-0.2000223777...$ being a parameter.}
\label{Fbv3}
\end{figure}

However, the geometry of the whole phase space and the structure
of key asymptotic bundles  change dramatically in comparison with
the blow-up cases, so that the standard shooting of {\em positive}
global profiles $F(y)$ by the {\tt ode45} solver yields no
encouraging results. We refer to Figure \ref{FNon88}, which
illustrates typical negative results of a standard shooting.
 Figure \ref{Fbv3} looks better and presents shooting a kind of
 ``separatrix", which however does not belong to the necessary
 family as in (\ref{s2Na}).
Actually, this means that a 1D shooting is not possible, and, as
we will see, there occurs a more complicated  2D one, i.e., using
two parameters.



\begin{figure}
\centering
\includegraphics[scale=0.70]{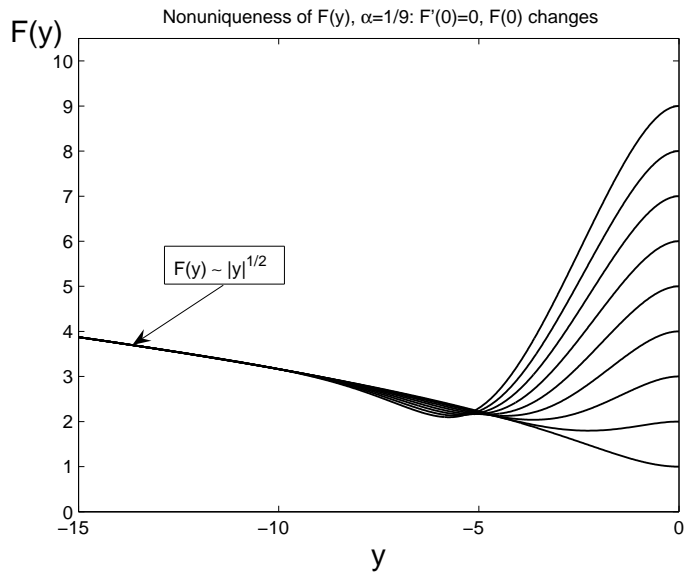}
\caption{\small Global profiles $F(y)$ of (\ref{s2Na}) for $\a=
\frac 19$, $C_0=1$, $F'(0)=0$;  $F(0)\in [1,9]$ being a
parameter.}
\label{Fnon1}
\end{figure}

\begin{figure}
\centering
\subfigure[$F'(0)=+1$]{
\includegraphics[scale=0.52]{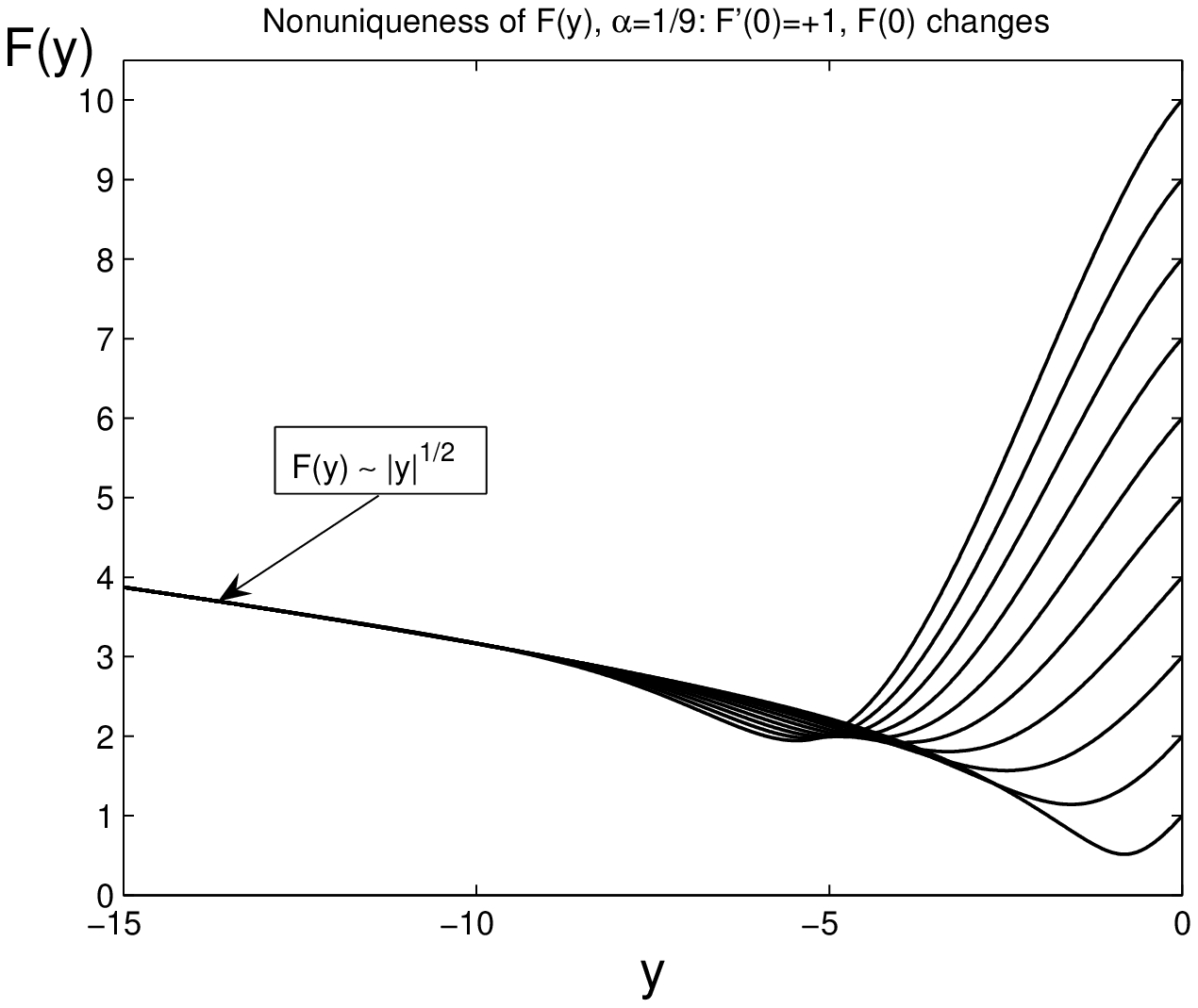}
}
\subfigure[$F'(0)=-1$]{
\includegraphics[scale=0.52]{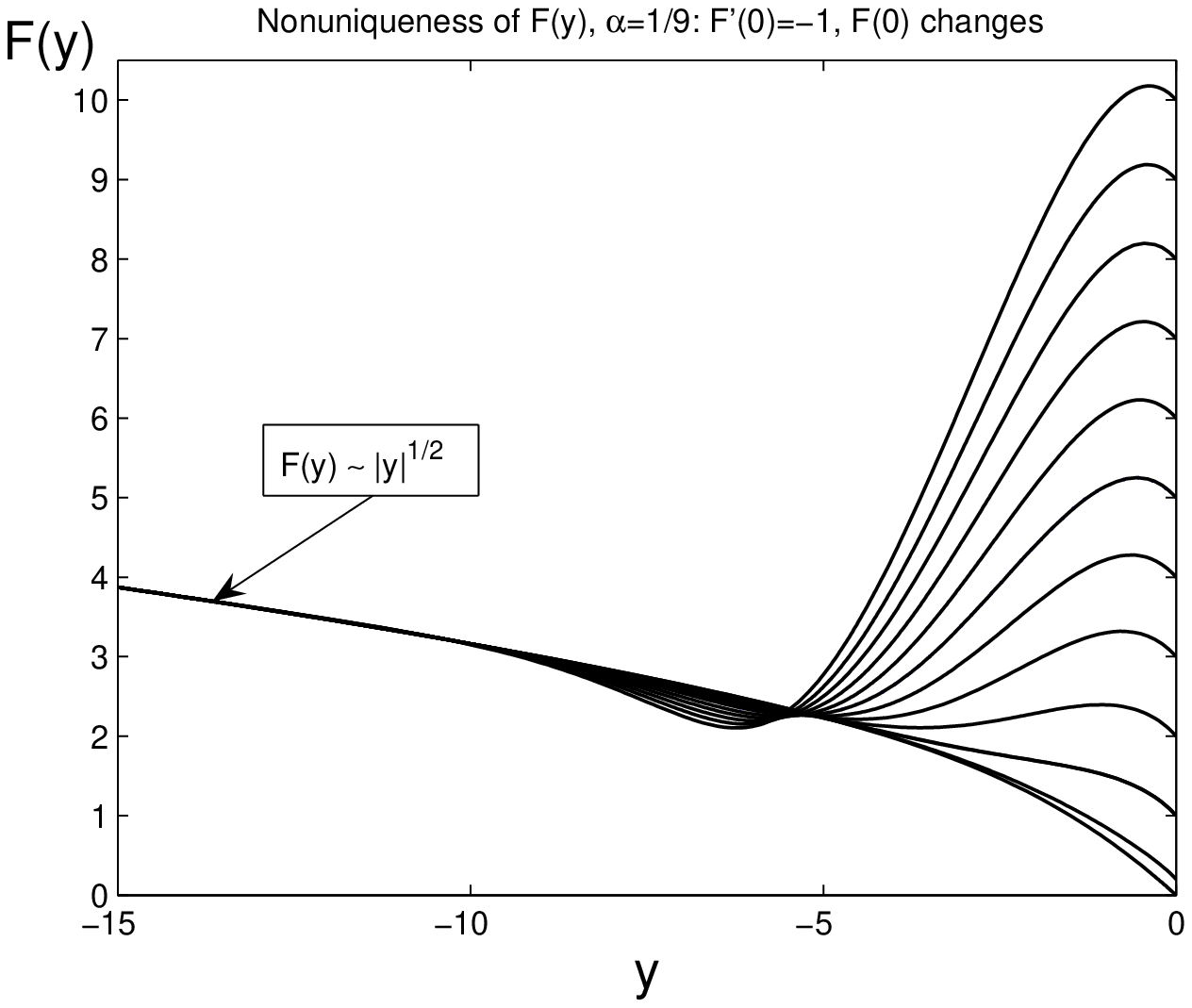}
}
 \vskip -.2cm
\caption{\rm\small Global profiles $F(y)$ of (\ref{s2Na}) for $\a=
\frac 19$, $C_0=1$ and $F'(0)=+1$ (a), $F'(0)=-1$ (b); $F(0)\in
[0,10]$ being a parameter.}
 \label{Fnon2}
\end{figure}

Therefore, we now use the {\tt bvp4c} solver, and this gives the
following results for the case (\ref{al12}), with $C_0=1$, as
usual.
Namely, we show that there are two parameters, say,
 \beq
 \label{par22}
  F_0=F(0) \andA F_1=F'(0),
   \eeq
   such that, for their arbitrary  values from some connected subset in
   $\re^2$, including all points with $F_0>0$ and $F_1\le 0$, the
   problem (\ref{s2Na}) admits a solution. This is confirmed in
   Figure \ref{Fnon1} for the case $F'(0)=0$ and in Figure
   \ref{Fnon2} for the cases $F'(0)=+1$ (a) and $F'(0)=-1$ (b).
   Obviously, all these profiles are different and exhibit fast and
   ``non-oscillatory"
    convergence as $y \to -\iy$
   to the ``good" bundle as in (\ref{s2Na}) with $C_0=1$.

Finally, carefully analyzing the dimensions of all the ``bad"
 and ``good" asymptotic bundles indicated in (i) and (ii) above, plus
 \ef{in23GaG}, unlike the result for blow-up profiles in
 Proposition \ref{Pr.1a}, we arrive at even  stronger nonuniqueness:

\begin{proposition}
\label{Pr.2}
 In the range $(\ref{Ran1})$ and any fixed $C_0>0$,
 the problem
$(\ref{s2Na})$ admits a $2$D family of solutions, which can be
parameterized by $F_0$ and $F_1$.
 \end{proposition}

Recall again that, for any hope of uniqueness, the extension pair
$\{f,F\}$ {\em must be unique} (or at least their subset should
contain some ``minimal" and/or isolated points as proper
candidates for unique
 entropy solutions) for any fixed constant
$C_0>0$, which defines the ``initial data" (\ref{s3a}) at the
blow-up time $t=0^-$. This actually happens for the Euler equation
(\ref{3}); see \cite[\S~4]{Gal3NDENew}, where the similarity
analysis is indeed easier and is reduced to algebraic
manipulations, but not that straightforward anyway even for such a
``first-order NDE".





 \subsection{``Initial nonuniqueness"}

\begin{figure}
\centering
\subfigure[$F'(0)=F'''(0)=-0.115526...$]{
\includegraphics[scale=0.52]{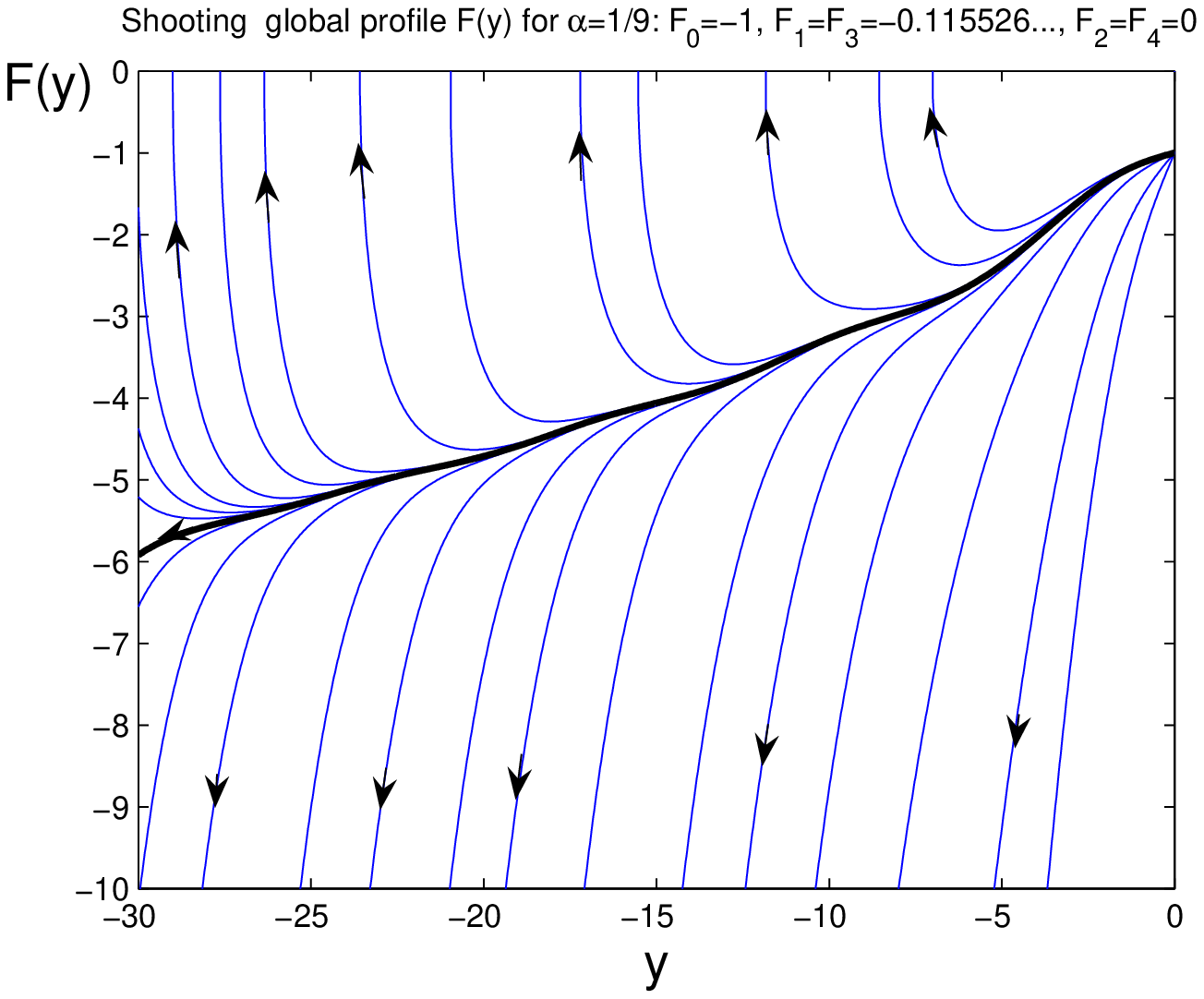}
}
\subfigure[$F'(0)=0, \,\, F''(0)=-0.16648...$]{
\includegraphics[scale=0.52]{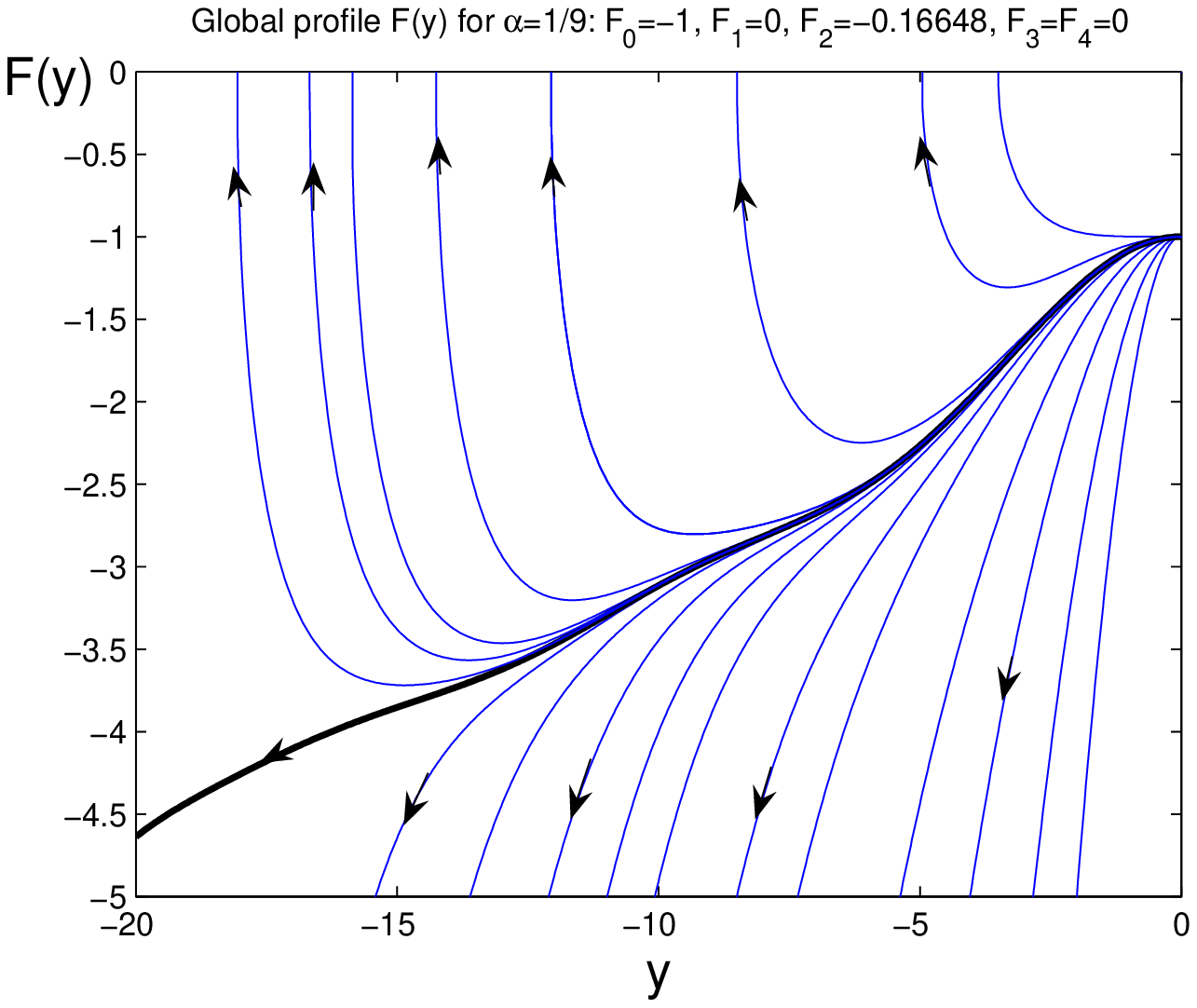}
}
 \vskip -.2cm
\caption{\rm\small  Shooting a proper solution  $F(y)$ of
\ef{s2Na} for $\a=\frac 19$ with data $F(0)=-1$,
 $F_4=0$, and $F'(0)=F'''(0)=-0.115526...$ (shooting parameter), $F_2=0$ (a), and
 $F(0)=-1$,
 $F'(0)=0$, $F_2=-0.16648...\,$ (shooting parameter), $F_3=F_4=0$ (b).}
 \label{F6a}
\end{figure}

A new ``nonuniqueness" phenomenon is achieved
for the values of parameters
  \beq
  \label{FF66}
 F(0)=F_0<0 \quad \mbox{and} \quad F'(0)=F_1 \ge 0.
   \eeq
  Figure \ref{F6a}(a), (b) shows such shock profiles leading to the nonuniqueness,
  obtained by a standard 1D shooting via the
   {\tt ode45} solver. Here, two similarity profiles $F(y)$
 are obtained via distinct types of shooting: relative to the
 parameter $F'(0)=F'''(0)$ in (a), and relative $F_2(0)$ in (b).

 The proof of existence of such profiles $F$ is based on the same
 geometric arguments as that of Proposition \ref{Pr.1} (with the
 evident change of the geometry of the phase space).
  These two different profiles posed into the similarity solutions (\ref{s1Na}) show
  a nonunique way to get  solutions with  initial data
  ($C_0=1$ by scaling) at $t=0^+$:
   \beq
   \label{FF77}
   u_0(x)=  |x|^{\frac \a \b} {\rm sign}\, x \quad \mbox{in} \quad
   \re,
    \eeq
    which already have a gradient blow-up singularity at $x=0$.
This is another potential type of nonuniqueness in the Cauchy
problem for \ef{N14}, showing the nonunique way of formation of
shocks from weak discontinuities, including the stationary ones as
in \ef{st1}.


\begin{figure}
\centering
\includegraphics[scale=0.70]{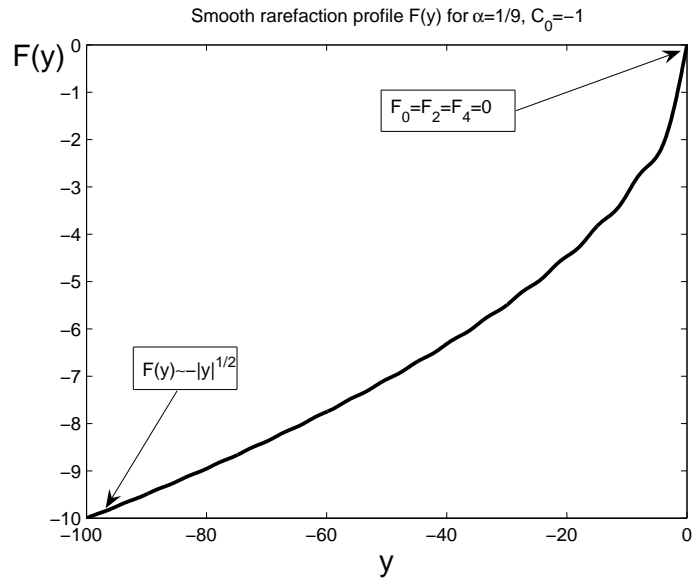}
\caption{\small Global rarefaction profile $F(y)$ of (\ref{s2Na})
for $\a= \frac 19$, $C_0=-1$; $F(0)=F''(0)=F^{(4)}(0)=0$.}
\label{Fbv10}
\end{figure}

 However, bearing in mind Proposition \ref{Pr.NE} saying that the
 shocks of $S_+$-type are not $\d$-entropy (i.e., not stable
 relative small smooth deformations), one can expect that the
 shocks as in (\ref{FF66}) are also unstable. Indeed, smooth
 extensions of weak pointwise shocks (\ref{FF77}) via rarefaction
 self-similar waves given by (\ref{s2Na}) are $\d$-entropy.
 In Figure \ref{Fbv10}, we show such a global rarefaction profile $F(y)$
for $\a=\frac 19$, which describes smooth collapse of the ``weak
equilibrium" (\ref{st1}). One can see that such rarefaction
profiles satisfy $F(y) \equiv -f(y)$, where $f$ are the
corresponding blow-up ones, as shown in Figure \ref{Fa2} for
various $\a$.

Overall, it seems that the $\d$-entropy test rules out such an
``initial nonuniqueness" with data of type $S_+$ as in
(\ref{FF77}), where a unique smooth rarefaction extension is
available. On the other hand, for other classes of data of
$S_-$-shape (according to Proposition \ref{Pr.E}), such a
nonuniqueness can take place; see Section \ref{S3aa}.

\subsection{More on nonuniqueness and well-posedness of FBPs}

The nonuniqueness \ef{s5N} in the Cauchy problem \ef{N14},
\ef{s3a} is as follows:  any  $F(y)$ yields the self-similar
continuation \ef{s1Na}, with the behaviour of the jump at $x=0$
(profiles $F(y)$ as in Figure \ref{F6a})
 \beq
 \label{jj1a}
  \mbox{$
 -[u_+(x,t)]\big|_{x=0}
 \equiv - \big(u_+(0^+,t)-u_+(0^-,t))
 = 2F_0 t^\a < 0 \forA t>0.
  $}
 \eeq
In the similarity  ODE representation, this nonuniqueness has a
pure geometric-dimensional origin associated with the dimension
and mutual
 geometry of the good and bad asymptotic bundles of the
5D phase spaces of both blow-up and global equations. Since these
shocks are stationary, the corresponding Rankine--Hugoniot (R--H)
condition on the speed $\l$ of the shock propagation:
 \beq
 \label{jj2}
 \mbox{$
 \l= \frac {[(uu_x)_{xxx}]}{[u]}\big|_{x=0}
 \equiv \frac {[(u^2)_{xxxx}]}{2[u]}\big|_{x=0}
 =  \frac {[(f^2)^{(4)}]}{2[f]}\big|_{y=0}
 =0
  $}
  \eeq
 is valid by anti-symmetry. As usual, \ef{jj2} is obtained
 by integration of the equation \ef{1} in a small neighbourhood of
 the shock.
  The R--H condition does
 not assume any novelty and is a corollary of integrating  the PDE
about the line of discontinuity.

Moreover, the R--H condition \ef{jj2} also indicates another
origin of nonuniqueness: a {\em symmetry breaking}. Indeed, the
solution for $t>0$ is not obliged to be an odd function of $x$, so
 the self-similar solution \ef{s1Na} for $x<0$ and $x>0$ can be defined using
 ten
different parameters $\{F_0^\pm,...,F_4^\pm\}$, and the only extra
condition one needs is the R--H one:
 \beq
 \label{kk1a}
 [(F F')'''](0)=0, \quad \mbox{i.e.,}
 \quad F_0^- F_4^- + 4 F_1^- F_3^-+ 3(F_2^-)^2=  F_0^+ F_4^+ + 4 F_1^+ F_3^+ +
 3(F_2^+)^2.
  \eeq
 This algebraic equations with {\em ten} unknowns admit
 many
 other solutions rather than the obvious anti-symmetric
 one:
  $$
  F_0^-=-F_0^+, \quad F_1^-=F_1^+, \quad F_2^-=-F_2^+, \quad F_3^-=F_3^+,
  \,\,\,
  \mbox{and} \,\,\, F_4^-=-F_4^+.
   $$

Finally, we note that the uniqueness can be restored by posing
specially designed conditions on moving shocks, which, overall
guarantee the unique solvability of the algebraic equation in
\ef{kk1a} and hence the unique continuation of the solution beyond
blow-up. This construction is analytically similar to that for the
NDEs--3 \ef{1} in \cite{Gal3NDENew}.

\section{Shocks for an NDE obeying the Cauchy--Kovalevskaya
theorem}
 \label{SCK1}

 In this short section, we touch on the problem of formation of
 shocks for NDEs that are higher-order in time. Instead of
 studying the PDEs such as (cf. \cite{GPndeII, GPnde})
  \beq
  \label{kk121}
  u_{tt}=-(uu_x)_{xxxx}, \quad  u_{ttt}=-(uu_x)_{xxxx}, \quad
  \mbox{etc.},
   \eeq
 we consider the fifth-order in time NDE (\ref{kk1}),
   which exhibits certain simple and, at the same time, exceptional properties.
Writing it for $W=(u,v,w,g,h)^T$ as
 \beq
 \label{jj1gg}
 \left\{
 \begin{matrix}
 u_t=v_x, \,\,\,\,\,\,\\
 v_t=w_x, \,\,\,\,\,\\
 w_t=g_x, \,\,\,\,\,\\
 g_t=h_x, \,\,\,\,\,\\
 h_t=uu_x,
  \end{matrix}
  \right.
  \quad \mbox{or} \quad W_t= A W_x, \quad \mbox{with the matrix}
  \quad
  A= \left[
   \begin{matrix}
   0\,\,\,1\,\,\,0\,\,\,0\,\,\,0\\
 0\,\,\,0\,\,\,1\,\,\,0\,\,\,0\\
  0\,\,\,0\,\,\,0\,\,\,1\,\,\,0\\
   0\,\,\,0\,\,\,0\,\,\,0\,\,\,1\\
    u
    \,\,\,0\,\,\,0\,\,\,0\,\,\,0
    \end{matrix}
    \right] ,
   \eeq
(\ref{kk1}) becomes a first-order system with the characteristic
equation for eigenvalues
 $$
 -\l^5 + u=0.
  $$
  Hence, for any $u \not = 0$, there exist complex roots,
   so that advanced results on hyperbolic systems
\cite{Bres, Daf} cannot be applied.

   \subsection{Evolution formation of shocks}

   For (\ref{kk1}),
  the blow-up  similarity
   solution is
 \beq
 \label{2.1tt}
 u_-(x,t)=g(z), \quad z= x/(-t), \quad \mbox{where}
  \eeq
 \beq
 \label{2.2tt}
 \begin{matrix}
  (g g')^{(4)}=(z^5g')^{(4)}\equiv 120 g'z+240 g'' z^2+ 120 g''' z^3 \ssk\ssk\ssk\\
    +
  20 g^{(4)} z^4+  g^{(5)} z^5
  \quad \mbox{in} \quad \re, \quad f(\mp
  \infty)=\pm 1.
   \end{matrix}
   \eeq
Integrating (\ref{2.2tt}) four times yields
 \beq
 \label{pp1}
  \mbox{$
 gg'=z^5 g' +Az+Bz^3, \quad \mbox{with constants} \quad A=(g'(0))^2>0, \,\, B=
  \frac 23\, g'(0)g'''(0),
   $}
  \eeq
 so that the necessary similarity profile $g(z)$ solves
 the first-order ODE
 \beq
 \label{ss1GG}
 \mbox{$
  \frac{{\mathrm d}g}{{\mathrm d}z}= \frac {A z + Bz^3}{g-z^5}.
  $}
  \eeq
By the phase-plane analysis of (\ref{ss1GG}) with $A>0$ and $B=0$,
we easily get the following:

\begin{proposition}
 \label{Pr.3}
 The problem $(\ref{2.2tt})$ admits a  solution $g(z)$
 satisfying the anti-symmetry conditions $(\ref{2.4})$ that is
 positive  for $z<0$, monotone decreasing, and is real analytic.
  \end{proposition}

\begin{figure}
\centering
\includegraphics[scale=0.70]{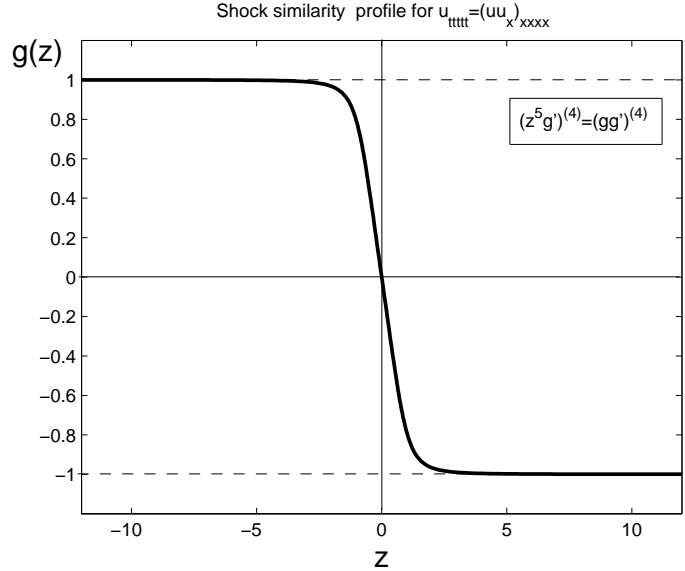}
 \vskip -.3cm
\caption{\small The shock similarity profile satisfying
(\ref{2.2tt}).}
\label{F1tt}
\end{figure}

Actually, involving the second parameter $B>0$ yields that there
exist infinitely many shock similarity profiles.
The boldface profile $g(z)$ in Figure \ref{F1tt} (by
(\ref{2.1tt}), it gives $S_-(x)$ as $t \to 0^-$) is
non-oscillatory about $\pm 1$, with the following  algebraic rate
of convergence
 to the equilibrium as $ z \to - \infty$:
 $$
 g(z) =
 \left\{
 \begin{matrix}
  1 + \frac A{3 z^{5}} +... \,\,\, \mbox{for} \,\,\, B=0, \ssk\\
  1+ \frac B z+... \,\,\,\,\,\, \mbox{for} \,\,\, B>0.
  \end{matrix}
  \right.
     $$


Note that the fundamental solutions of the corresponding linear
PDE
 \beq
 \label{dd1}
 u_{ttttt}=u_{xxxxx}
 \eeq
 is also not oscillatory as $x \to \pm \infty$.
 This has the form
  $$
  b(x,t) = t^3 F(y), \quad y=x/t, \quad \mbox{so that}
  \,\,\, b(x,0)=...=b_{ttt}(x,0)=0, \,\,\, b_{tttt}(x,0)=\d(x).
   $$
  The linear equation
  (\ref{dd1}) exhibits some features of  finite propagation via
 TWs, since
 $$
 u(x,t)=f(x- \l t) \quad \Longrightarrow \quad
  -\l^5 f^{(5)}=f^{(5)}, \,\,\,\, \mbox{i.e.,} \,\,\, \l=-1,
   $$
   since the profile $f(y)$ disappears from the ODE. This is similar to
   some canonical equations of
   mathematical physics such as
    $$
    u_t=u_x \,\,(\mbox{dispersion, $\l=-1$}) \quad \mbox{and} \quad
  u_{tt}=u_{xx} \,\,(\mbox{wave equation, $\l=\pm 1$}).
   $$
 The blow-up  solution (\ref{2.1tt}) gives
 in the limit $t \to 0^-$ the shock $S_-(x)$, and (\ref{con32})
 holds.

Since (\ref{kk1}) has the same symmetry (\ref{symm1}) as
(\ref{1}), similarity solutions (\ref{2.1tt}),
 with $-t \mapsto t$ and $g(z) \mapsto g(-z)$ according to (\ref{2.15}),
 also give the
rarefaction waves for $S_+(x)$ as well as other types of collapse
of initial non-entropy discontinuities.


\subsection{Analytic $\d$-deformations by Cauchy--Kovalevskaya
theorem}

The great advantage of the equation (\ref{kk1}) is that it is in
the {\em normal form}, so it obeys the Cauchy--Kovalevskaya
theorem \cite[p.~387]{Tay}. Hence, for any analytic initial data
$u(x,0)$, $u_t(x,0)$, $u_{tt}(x,0)$, $u_{ttt}(x,0)$, and
$u_{tttt}(x,0)$, there exists a unique local in time analytic
solution $u(x,t)$. Thus, (\ref{kk1}) generates a local semigroup
of
 analytic solutions, and this makes it easier to deal with smooth
 $\d$-deformations that are chosen to be analytic.
This defines a special analytic $\d$-entropy test for
 shock/rarefaction waves.
  On the other hand, such nonlinear PDEs can
 admit other (say, weak) solutions that are not analytic.
Actually, Proposition \ref{Pr.3} shows that the shock $S_-(x)$
 is a $\d$-entropy solution of (\ref{kk1}), which is obtained by finite-time blow-up
  as $t \to 0^-$
 from the analytic similarity solution (\ref{2.1tt}).

\subsection{On formation of single-point shocks and extension
nonuniqueness}

Similar to the analysis in Section \ref{SNonU}, for the model
\ef{kk1} (and \ef{kk121}), these assume studying extension
similarity  pairs $\{f,F\}$ induced by the easy derived analogies
of the blow-up \ef{s2a} and global \ef{s2Na},
 with
  $$
  \mbox{$
  \b= \frac {5+\a}5,
  $}
  $$
 5D dynamical
systems. These are very difficult, so that checking three types
(standard, backward, and initial) of possible nonuniqueness and
non-entropy of such flows with strong and weak shocks becomes a
hard open problem, though some auxiliary analytic steps towards
nonuniqueness are doable. Overall, in view of complicated
multi-dimensional phase spaces involved, we do not have any reason
for having a unique continuation after singularity. In other
words, for such higher-order NDEs, uniqueness can occur {\em
accidentally} only for very special phase spaces, and hence, at
least, is not robust (in a natural ODE--PDE sense) anyway.

 \section{{\bf (V)} {\bf Problem ``Oscillatory Smooth Compactons"} of fifth-order NDEs}
 \label{Sect6}

We begin with an easier explicit  example of {\em nonnegative}
compactons for a third-order NDE.

\subsection{Third-order NDEs: $\d$-entropy compactons}

Compactons as compactly supported TW solutions of the $K(2,2)$
equation (\ref{Comp.4}) were introduced in 1993, \cite{RosH93},
 as
 \beq
 \label{co1}
 u_{\rm c}(x,t)=f_{\rm c}(y), \,\,\, y=x + t
 \quad \Longrightarrow \quad f_{\rm c}: \quad
  f=(f^2)''+f^2.
  \eeq
  Integrating yields the following explicit
 compacton profile:
 \beq
 \label{co11}
 f_{\rm c}(y)= \left\{
 \begin{matrix}
  \frac 43 \cos^2\big(\frac y4\big) \quad \mbox{for} \quad |y| \le
  2\pi, \\
  \qquad\,\,\, 0 \qquad\,  \mbox{for} \quad |y| \ge
  2\pi.
   \end{matrix}
   \right.
   \eeq
The corresponding compacton (\ref{co1}), (\ref{co11}) is
G-admissible in the sense of
$\fbox{$\mbox{Gel'fand}$}$\footnote{I.M.~Gel'fand, 2.09.1913--5.10.2009.} (1959)
\cite[\S\S~2,\,8]{Gel},
 and is a
$\d$-entropy solution  \cite[\S~4]{GPndeII}, i.e., can be
constructed by smooth (and moreover analytic) approximation via
strictly positive solutions of the full third-order ODE for
$f_{\rm c}(y)$ ,
 $$
 f'=(f^2)'''+(f^2)'.
  $$
 Since the PDE is not involved unlike Section \ref{SD1},
 the $\d$-entropy notion coincides with the G-admissability.

It is curious that the same compactly supported blow-up patterns
occur in the combustion problem for the related reaction-diffusion
parabolic equation
 \beq
 \label{co2}
 u_t=(u^2)_{xx}+u^2.
 \eeq
Then the standing-wave blow-up (as $t \to T^-$) solution of
S-regime leads to the same ODE:
 \beq
 \label{co3}
 u_{\rm S}(x,t)=(T-t)^{-1} f(x) \quad \Longrightarrow \quad
  f=(f^2)''+f^2.
  \eeq
This yields the {\em Zmitrenko--Kurdyumov blow-up localized
solution}, which has been known since 1975; see more historical
details in \cite[\S~4.2]{GSVR}.

\subsection{Examples of $C^3$-smooth nonnegative compacton for higher-order NDEs}

Such an example was given in  \cite[p.~4734]{Dey98}. Following
\cite[p.~189]{GSVR}, we construct this explicit solution as
follows.
 The operator ${\bf F}_5(u)$ of the
{\em quintic NDE}
 \beq
 \label{qq1}
 u_t= {\bf F}_5(u) \equiv (u^2)_{xxxxx}+ 25 (u^2)_{xxx} + 144 (u^2)_x
 \eeq
 is shown to
  preserve
   the 5D invariant subspace
 \beq
\label{555.2II}
 W_5= {\rm Span}\{1, \cos x, \sin x, \cos 2x,
\sin 2x\},
 \eeq
 i.e., ${\bf F}_5(W_5) \subseteq W_5$. Therefore, (\ref{qq1}) restricted
 to the invariant subspace $W_5$ is a 5D dynamical system for the
 expansion coefficients of the solution
 $$ 
u(x,t)= C_1(t)+C_2(t) \cos x+C_3(t) \sin x +C_4(t) \cos 2x +
C_5(t) \sin 2x \in W_5.
 $$ 
 Solving this yields the explicit compacton TW
 \beq
 \label{CC.191}
 u_{\rm c}(x,t) = f_{\rm c}(x+t), \quad \mbox{where} \quad f_{\rm c}(y)=
 \left\{
  \begin{matrix}
  \frac 1{105} \, \cos^4 \bigl(\frac y2\bigr) \,\,\,
  \mbox{for} \,\,\,|y| \le \pi,\ssk \\
  \qquad \quad 0 \qquad  \mbox{for} \,\,\,|y| \ge \pi.
  \end{matrix}
  \right.
  \eeq
This $C^{3}_x$ solution can be attributed to the Cauchy problem
for (\ref{qq1}) since smooth solutions are not oscillatory near
interfaces; see a discussion around \cite[p.~184]{GSVR}.

\ssk


The above invariant subspace analysis applies also to the
7th-order PDE
 \beq
\label{777.1} \mbox{$
  u_t ={\bf F}_7(u) \equiv D_x^7(u^2) +  \b
D_x^5(u^2) + \g (u^2)_{xxx} + \nu (u^2)_x.
 $}
 \eeq
 Here ${\bf F}_7$ admits $W_5$ if
 $$
 \mbox{
 $\b=25$, \,\,$\g = 144$,\,\, and $\nu=0$}.
 $$
 Moreover \cite[p.~190]{GSVR},
the only operator ${\bf F}_7$ in (\ref{777.1}) preserving the 7D
subspace
 \beq
\label{777.2} W_7 = {\mathcal L}\{1, \cos x, \sin x, \cos 2x, \sin
2x, \cos 3x, \sin 3x\}
 \eeq
is in the following NDE--7:
 \beq
\label{555.3}
 u_t= {\bf F}_7(u) \equiv D_x^7(u^2) +  77 D_x^5(u^2) +
1876 (u^2)_{xxx} + 14400 (u^2)_x.
 \eeq
This makes it possible to reduce (\ref{555.3}) on $W_7$ to a
complicated dynamical system.

 \subsection{Why nonnegative compactons for fifth-order NDEs are
 not robust: a saddle-saddle homoclinic}

 Recall that, as usual in dynamical system theory,
  by robustness of trajectories  we mean that these are
 stable with respect to small perturbations of the parameters
 entering the NDE or  the corresponding ODEs.
 In other words, the dynamical systems (ODEs) admitting such non-negative ``heteroclinic"
 {\em saddle-like}
 orbits $0 \to 0$ are not {\em structurally stable} in a natural sense.
 This reminds the classic Andronov--Pontriagin--Peixoto theorem,
 where one of the four conditions for the structural stability of
 dynamical systems in $\re^2$ reads as follows
 \cite[p.~301]{Perko}:
  \beq
  \label{het11}
  \mbox{``(ii) there are no trajectories connecting saddle
  points...\,."}
   \eeq
   Actually, nonnegative compactons, such as  \ef{CC.191},
    are special homoclinics of the origin, and
   we will show that the nature of their non-robustness is in the
   fact that these represent a stable-unstable manifold of the
   origin consisting of a {\em single orbit}. Therefore, in consistency
   with (\ref{het11}), the origin is indeed a saddle in $\re^4$ in the plane
   $\{f,f',f'',f'''\}$, obtained after integration once; see
   below.

 In order to illustrate the lack of such a robustness in view of a sole heteroclinic involved,
  consider
 the NDE
  (\ref{qq1}), for which, substituting the TW solution,
   on integration, we obtain the following ODE:
  \beq
  \label{gg1s}
  u_{\rm c}(x,t)=f_{\rm c}(x+t)
  \quad \Longrightarrow \quad f_{\rm c}: \quad  2f= (f^2)^{(4)}+... \, ,
   \eeq
   where we omit the lower-order terms as $f \to 0$. Looking for the compacton
  profile $f \ge 0$, we set $f^2=F$ to get
   \beq
   \label{gg2s}
   F^{(4)}= 2 \sqrt F+...  \quad \mbox{for}
   \quad  y>0, \quad F'(0)=F'''(0)=0.
    \eeq
 As usual, we  look  for a symmetric  $F(y)$ by
 putting  two symmetry conditions at the origin.

 Let $y=y_0>0$ be the interface point of $F(y)$. Then, looking for
 the expansion as $y \to y_0^-$ in the form
  \beq
  \label{gg3}
   \mbox{$
  F(y)= \frac 1{840^2} \, (y_0-y)^8 + \e(y), \quad
  \mbox{with}
  \quad \e(y) = o((y_0-y)^8),
 $}
   \eeq
   we obtain Euler's equation for the perturbation $\e(y)$,
 \beq
 \label{gg4}
 \mbox{$
 \frac 1{840}(y_0-y)^4 \e^{(4)}-\e=0.
  $}
  \eeq
Hence, $\e(y)=(y_0-y)^m$, with the characteristic equation
 \beq
 \label{gg5}
  \mbox{$
 m(m-1)(m-2)(m-3)-840=0
 \, \Longrightarrow
 \, m_1=-4, \,\,\, m_{2,3}= \frac{3 \pm {\rm i} \sqrt{111}}{2},
 \quad m_4=7.
 $}
  \eeq
  Hence, ${\rm Re}\, m_i <8$, and,
in other words, (\ref{gg4}) does not admit any nontrivial solution
satisfying the condition in (\ref{gg3}); see further comments in
\cite[p.~142]{GSVR}.
 In fact, it is easy to see that (\ref{gg4}) with $\e=0$ is the
 unique positive smooth solution of $F^{(4)}=2 \sqrt F$. Thus,
 \beq
 \label{gg6}
 \mbox{the asymptotic bundle of solutions (\ref{gg3}) is 1D},
  \eeq
  where the only parameter is the position of the interface
  $y_0>0$.

  Obviously, as a typical property, this 1D bundle
 is not sufficient to satisfy (by shooting)
 {\sc two} conditions at the origin in (\ref{gg2s}), so such TW
 profiles $F(y) \ge 0$ are nonexistent for almost all NDEs like
 that.
 In other words, the condition of positivity of the solution,
  \beq
  \label{pos11}
  \mbox{to look for a nontrivial solution $F \ge 0$ for the ODE in
  (\ref{gg2s})}
   \eeq
 creates a free-boundary  ``obstacle" problem that, in general,
is inconsistent. Skipping the obstacle condition (\ref{pos11})
will return such ODEs (or elliptic equations), with a special
extension, into the consistent variety, as we will illustrate
below.

\ssk

 Thus, nonnegative TW compactons are not generic (robust) solutions
 of $(2m+1)$th-order quadratic NDEs with $m=2$, and also for larger
 $m$'s, where some kind of (\ref{gg6}), as a ``dimensional defect"
 (the bundle dimension is smaller than the number of conditions at $y=0$ to shoot), remains valid.

\subsection{Nonnegative compactons are robust for third-order NDEs only}

 The
 third-order case $m=1$, i.e., NDEs such as (\ref{Comp.4}), is the
 only one where propagation of perturbations via nonnegative
 TW compactons is structurally stable, i.e.,  with respect to small perturbation of the
 parameters (and nonlinearities) of equations. Mathematically speaking, then the
 1D bundle in (\ref{gg6}) perfectly matches with the {\sc single}
 symmetry condition at the origin,
  $$
 F''= 2 \sqrt F+... \quad \mbox{and} \quad  F'(0)=0.
  $$

 \subsection{Compactons of changing sign are robust and $\d$-entropy for the NDEs--5}

As a typical example, we consider the  perturbed version
(\ref{z1})
 of the NDE--(1,4) (\ref{N14}). As we have mentioned, this is
    is written for solutions of changing sign,
    since nonnegative compactons do not exist in general.
   Looking for the TW compacton (\ref{gg1s}) yields the ODE
\beq
 \label{z2}
 \mbox{$
 f= - \frac 12 \, (|f|f)^{(4)}+ \frac 12\, |f|f
 \quad \Longrightarrow
 \quad F^{(4)}=F- 2|F|^{-\frac 12}F \,\,\, \mbox{for} \,\,\, F=|f|f.
  $}
   \eeq
Such ODEs with non-Lipschitz nonlinearities are known to admit
countable sets of compactly supported solutions, which are studied
by a combination of Lusternik--Schnirel'man and Pohozaev's
fibering theory; see \cite{GMPSob}.

\begin{figure}
\centering
\includegraphics[scale=0.70]{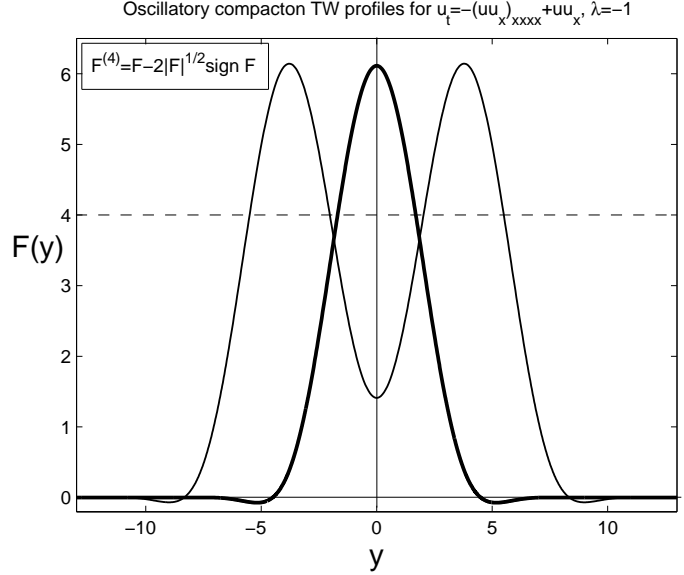}
 \vskip -.3cm
\caption{\small First two compacton TW profiles $F(y)$ satisfying
the ODE in (\ref{z2}).}
\label{FOsc1}
\end{figure}

 In Figure \ref{FOsc1}, we
present the first TW compacton patterns (the boldface line) and
the second one that is essentially non-monotone. These look like
standard compacton profiles but careful analysis of the behaviour
near the finite interface at $y=y_0$ shows that $F(y)$ changes
sign infinitely many times according to the asymptotics
 \beq
 \label{z4}
 F(y)=(y_0-y)^8 [\varphi(s+s_0)+o(1)], \,\,\, s= \ln(y_0-y)
  \quad \mbox{as} \quad y \to y_0^-.
  \eeq
Here,  the oscillatory component $\varphi(s)$ is a periodic
solution of a certain nonlinear ODE and $s_0$ is an arbitrary
phase shift; see \cite[\S~4.3]{GSVR} and \cite[\S~4]{GMPSob} for
further details. Thus, unlike (\ref{gg6}),
 \beq
 \label{zzz1}
\mbox{the asymptotic bundle of solutions (\ref{z4}) is 2D
(parameters are $y_0$ and $s_0$)},
 \eeq
 and exhibits some features of a ``nonlinear focus" (not a saddle
 as above) on some manifold.
 Hence, this is enough to match also two symmetry boundary conditions given  in
 (\ref{gg2s}). Such a robust solvability is confirmed by
 variational techniques that apply to
 rather arbitrary equations such as in (\ref{z2}) with similar singular
non-Lipschitz nonlinearities.

Let us also note that such oscillatory compactons are also
$\d$-entropy in the sense that can be approximated by analytic TW
solutions of the same ODE but having finite number of zeros (i.e.,
admit smooth analytic $\d$-deformations); see
\cite[\S~8.3]{GPndeII} for a related NDE--5.

\ssk

Regardless the existence of such sufficiently smooth compacton
solutions,  it is worth recalling again
 that, for the NDE (\ref{z1}), as well as (\ref{nn1}) and
 (\ref{ss3}),
 both containing monotone nonlinearities,
 the generic
behaviour, for other initial data, can include formation of shocks
in finite time, with the local similarity mechanism as in Section
\ref{Sect2} and in  Section \ref{S5.1}, representing   more
generic single point shock pattern formations.


\section{Final conclusions}

 The  fifth-order nonlinear dispersion equations
(\ref{N50})--(\ref{N14}) (NDEs--5), which are associated with a
number of important applications,
 are considered. The main achieved properties of such 1D degenerate
 nonlinear PDEs are as follows:

\ssk

 \noi(i)  \underline{\sc Section \ref{Sect2}}: these NDEs
admit blow-up self-similar formation from smooth solutions of {\em
shock waves} of a specific oscillatory structure, which correspond
to initial data $S_-(x)=-{\rm sign}\, x$, i.e., the same as for
the first-order conservation law \ef{1};

\ssk

 \noi(ii) \underline{\sc Section \ref{Sect3}}: as customary,
{\em self-similar rarefaction waves}, which get smooth for any
$t>0$, are created by the reversed data $S_+(x)={\rm sign}\, x$;


\ssk

 \noi(iii)
Unlike the classical theory of first-order conservation laws
developed in the 1950s--60s, the entropy-type techniques are no
longer applied for distinguishing general proper (unique)
solutions. A {\em $\d$-deformation test} via smoothing the
solutions is developed in  \underline{\sc Section \ref{SS1}},
which is able to separate shocks and rarefaction waves for
particular classes of initial data $\sim S_\pm(x)$ with a simple
geometry of initial shocks;

\ssk

\noi(iv)  \underline{\sc Section \ref{SNonU}}: by studying more
general self-similar solutions of NDEs, it was shown that
uniqueness of the solutions after formation of a shock is
principally impossible. Namely, there exist single point gradient
blow-up similarity solutions, which admit an infinite  number of
 self-similar extensions beyond. This 2D set of shock wave extensions
after singularity does not have any distinguished solutions (say,
maximal, minimal, isolated, etc.). This also suggests that  any
entropy-like mechanisms for a unique continuation do not exist
either. However, using a proper {\em free-boundary setting}, i.e.,
posing special conditions on shocks, can restore uniqueness; and

\ssk

 \noi(v)  \underline{\sc Section \ref{Sect6}}: {\em
nonnegative compacton solutions} for some NDEs--5 are shown to be
non-robust (not ``structurally stable"), i.e., these disappear
after a.a. arbitrarily small perturbations of the parameters
(nonlinearities) of the equations. However, {\em oscillatory
compactons} of changing sign near finite interface are shown to
exist and to be robust.


\ssk

{\bf Acknowledgements.}  The authors would like to thank the
anonymous Referee from the European Journal of Applied Mathematics
for an extremely useful discussion concerning $\d$-deformation
issues for shock/rarefaction waves.


\enddocument